%
%

%
%

\def\fmtversion{2.0}
\catcode`\@=11
\ifx\amstexloaded@\relax\catcode`\@=\active
 \endinput\else\let\amstexloaded@\relax\fi
\def\W@{\immediate\write\sixt@@n}
\def\CR@{\W@{}\W@{AmS-TeX - Version \fmtversion}\W@{}
\W@{COPYRIGHT 1985, 1990 - AMERICAN MATHEMATICAL SOCIETY}
\W@{Use of this macro package is not restricted provided}
\W@{each use is acknowledged upon publication.}\W@{}}
\CR@
\everyjob{\CR@}
\toksdef\toks@@=2
\long\def\rightappend@#1\to#2{\toks@{\\{#1}}\toks@@
 =\expandafter{#2}\xdef#2{\the\toks@@\the\toks@}\toks@{}\toks@@{}}
\def\alloclist@{}
\newif\ifalloc@
\def\showallocations{{\def\\{\immediate\write\m@ne}\alloclist@}\alloc@true}
\def\alloc@#1#2#3#4#5{\global\advance\count1#1by\@ne
 \ch@ck#1#4#2\allocationnumber=\count1#1
 \global#3#5=\allocationnumber
 \edef\next@{\string#5=\string#2\the\allocationnumber}%
 \expandafter\rightappend@\next@\to\alloclist@}
\newcount\count@@
\newcount\count@@@
\def\FN@{\futurelet\next}
\def\DN@{\def\next@}
\def\DNii@{\def\nextii@}
\def\RIfM@{\relax\ifmmode}
\def\RIfMIfI@{\relax\ifmmode\ifinner}
\def\setboxz@h{\setbox\z@\hbox}
\def\wdz@{\wd\z@}
\def\boxz@{\box\z@}
\def\setbox@ne{\setbox\@ne}
\def\wd@ne{\wd\@ne}
\def\iterate{\body\expandafter\iterate\else\fi}
\newlinechar=`\^^J
\def\err@#1{\errmessage{AmS-TeX error: #1}}
\newhelp\defaulthelp@{Sorry, I already gave what help I could...^^J
Maybe you should try asking a human?^^J
An error might have occurred before I noticed any problems.^^J
``If all else fails, read the instructions.''}
\def\Err@#1{\errhelp\defaulthelp@\errmessage{AmS-TeX error: #1}}
\def\eat@#1{}
\def\in@#1#2{\def\in@@##1#1##2##3\in@@{\ifx\in@##2\in@false\else\in@true\fi}%
 \in@@#2#1\in@\in@@}
\newif\ifin@
\def\space@.{\futurelet\space@\relax}
\space@. %
\newhelp\athelp@
{Only certain combinations beginning with @ make sense to me.^^J
Perhaps you wanted \string\@\space for a printed @?^^J
I've ignored the character or group after @.}
\def\futureletnextat@{\futurelet\next\at@}
{\catcode`\@=\active
\lccode`\Z=`\@ \lccode`\I=`\I \lowercase{%
\gdef@{\csname futureletnextatZ\endcsname}
\expandafter\gdef\csname atZ\endcsname                                      
 {\ifcat\noexpand\next a\def\next{\csname atZZ\endcsname}\else
 \ifcat\noexpand\next0\def\next{\csname atZZ\endcsname}\else
 \ifcat\noexpand\next\relax\def\next{\csname atZZZ\endcsname}\else
 \def\next{\csname atZZZZ\endcsname}\fi\fi\fi\next}
\expandafter\gdef\csname atZZ\endcsname#1{\expandafter                      
 \ifx\csname #1Zat\endcsname\relax\def\next
 {\errhelp\expandafter=\csname athelpZ\endcsname
 \csname errZ\endcsname{Invalid use of \string@}}\else
 \def\next{\csname #1Zat\endcsname}\fi\next}
\expandafter\gdef\csname atZZZ\endcsname#1{\expandafter                     
 \ifx\csname \string#1ZZat\endcsname\relax\def\next
 {\errhelp\expandafter=\csname athelpZ\endcsname
 \csname errZ\endcsname{Invalid use of \string@}}\else
 \def\next{\csname \string#1ZZat\endcsname}\fi\next}
\expandafter\gdef\csname atZZZZ\endcsname#1{\errhelp                        
 \expandafter=\csname athelpZ\endcsname
 \csname errZ\endcsname{Invalid use of \string@}}}}                         
\def\atdef@#1{\expandafter\def\csname #1@at\endcsname}
\def\atdef@@#1{\expandafter\def\csname \string#1@@at\endcsname}
\newhelp\defahelp@{If you typed \string\define\space cs instead of
\string\define\string\cs\space^^J
I've substituted an inaccessible control sequence so that your^^J
definition will be completed without mixing me up too badly.^^J
If you typed \string\define{\string\cs} the inaccessible control sequence^^J
was defined to be \string\cs, and the rest of your^^J
definition appears as input.}
\newhelp\defbhelp@{I've ignored your definition, because it might^^J
conflict with other uses that are important to me.}
\def\define{\FN@\define@}
\def\define@{\ifcat\noexpand\next\relax
 \expandafter\define@@\else\errhelp\defahelp@                               
 \err@{\string\define\space must be followed by a control
 sequence}\expandafter\def\expandafter\nextii@\fi}                          
\def\undefined@@@@@@@@@@{}
\def\preloaded@@@@@@@@@@{}
\def\next@@@@@@@@@@{}
\def\define@@#1{\ifx#1\relax\errhelp\defbhelp@                              
 \err@{\string#1\space is already defined}\DN@{\DNii@}\else
 \expandafter\ifx\csname\expandafter\eat@\string                            
 #1@@@@@@@@@@\endcsname\undefined@@@@@@@@@@\errhelp\defbhelp@
 \err@{\string#1\space can't be defined}\DN@{\DNii@}\else
 \expandafter\ifx\csname\expandafter\eat@\string#1\endcsname\relax          
 \global\let#1\undefined\DN@{\def#1}\else\errhelp\defbhelp@
 \err@{\string#1\space is already defined}\DN@{\DNii@}\fi
 \fi\fi\next@}
\let\redefine\def
\def\predefine#1#2{\let#1#2}
\def\undefine#1{\let#1\undefined}
\newdimen\captionwidth@
\captionwidth@\hsize
\advance\captionwidth@-1.5in
\def\pagewidth#1{\hsize#1\relax
 \captionwidth@\hsize\advance\captionwidth@-1.5in}
\def\pageheight#1{\vsize#1\relax}
\def\hcorrection#1{\advance\hoffset#1\relax}
\def\vcorrection#1{\advance\voffset#1\relax}

\let\graveaccent\`
\let\acuteaccent\'
\let\tildeaccent\~
\let\hataccent\^
\let\underscore\_
\let\B\=
\let\D\.
\let\ic@\/
\def\/{\unskip\ic@}
\def\textfonti{\the\textfont\@ne}
\def\t#1#2{{\edef\next@{\the\font}\textfonti\accent"7F \next@#1#2}}
\def~{\unskip\nobreak\ \ignorespaces}
\def\.{.\spacefactor\@m}
\atdef@;{\leavevmode\null;}
\atdef@:{\leavevmode\null:}
\atdef@?{\leavevmode\null?}
\def\@{\char64 }
\def\relaxnext@{\let\next\relax}
\atdef@-{\relaxnext@\leavevmode
 \DN@{\ifx\next-\DN@-{\FN@\nextii@}\else
  \DN@{\leavevmode\hbox{-}}\fi\next@}%
 \DNii@{\ifx\next-\DN@-{\leavevmode\hbox{---}}\else
  \DN@{\leavevmode\hbox{--}}\fi\next@}%
 \FN@\next@}
\def\srdr@{\kern.16667em}
\def\drsr@{\kern.02778em}
\def\sldl@{\kern.02778em}
\def\dlsl@{\kern.16667em}
\atdef@"{\unskip\relaxnext@
 \DN@{\ifx\next\space@\DN@. {\FN@\nextii@}\else
  \DN@.{\FN@\nextii@}\fi\next@.}%
 \DNii@{\ifx\next`\DN@`{\FN@\nextiii@}\else
  \ifx\next\lq\DN@\lq{\FN@\nextiii@}\else
  \DN@####1{\FN@\nextiv@}\fi\fi\next@}%
 \def\nextiii@{\ifx\next`\DN@`{\sldl@``}\else\ifx\next\lq
  \DN@\lq{\sldl@``}\else\DN@{\dlsl@`}\fi\fi\next@}%
 \def\nextiv@{\ifx\next'\DN@'{\srdr@''}\else
  \ifx\next\rq\DN@\rq{\srdr@''}\else\DN@{\drsr@'}\fi\fi\next@}%
 \FN@\next@}

\def\textfontii{\the\textfont\tw@}
\def\lbrace@{\delimiter"4266308 }
\def\rbrace@{\delimiter"5267309 }
\def\{{\RIfM@\lbrace@\else{\textfontii f}\spacefactor\@m\fi}
\def\}{\RIfM@\rbrace@\else
 \let\@sf\empty\ifhmode\edef\@sf{\spacefactor\the\spacefactor}\fi
 {\textfontii g}\@sf\relax\fi}
\let\lbrace\{
\let\rbrace\}
\def\AmSTeX{{\textfontii A}\kern-.1667em\lower.5ex\hbox
 {\textfontii M}\kern-.125em{\textfontii S}-\TeX}
\def\vmodeerr@#1{\Err@{\string#1\space not allowed between paragraphs}}
\def\mathmodeerr@#1{\Err@{\string#1\space not allowed in math mode}}
\def\linebreak{\RIfM@\mathmodeerr@\linebreak\else
 \ifhmode\unskip\unkern\break\else\vmodeerr@\linebreak\fi\fi}

\newskip\saveskip@
\def\allowlinebreak{\RIfM@\mathmodeerr@\allowlinebreak\else
 \ifhmode\saveskip@\lastskip\unskip
 \allowbreak\ifdim\saveskip@>\z@\hskip\saveskip@\fi
 \else\vmodeerr@\allowlinebreak\fi\fi}
\def\nolinebreak{\RIfM@\mathmodeerr@\nolinebreak\else
 \ifhmode\saveskip@\lastskip\unskip
 \nobreak\ifdim\saveskip@>\z@\hskip\saveskip@\fi
 \else\vmodeerr@\nolinebreak\fi\fi}
\def\newline{\relaxnext@
 \DN@{\RIfM@\expandafter\mathmodeerr@\expandafter\newline\else
  \ifhmode\ifx\next\par\else
  \expandafter\unskip\expandafter\null\expandafter\hfill\expandafter\break\fi
  \else
  \expandafter\vmodeerr@\expandafter\newline\fi\fi}%
 \FN@\next@}
\def\dmatherr@#1{\Err@{\string#1\space not allowed in display math mode}}
\def\nondmatherr@#1{\Err@{\string#1\space not allowed in non-display math
 mode}}
\def\onlydmatherr@#1{\Err@{\string#1\space allowed only in display math mode}}
\def\nonmatherr@#1{\Err@{\string#1\space allowed only in math mode}}
\def\mathbreak{\RIfMIfI@\break\else
 \dmatherr@\mathbreak\fi\else\nonmatherr@\mathbreak\fi}
\def\nomathbreak{\RIfMIfI@\nobreak\else
 \dmatherr@\nomathbreak\fi\else\nonmatherr@\nomathbreak\fi}
\def\allowmathbreak{\RIfMIfI@\allowbreak\else
 \dmatherr@\allowmathbreak\fi\else\nonmatherr@\allowmathbreak\fi}
\def\pagebreak{\RIfM@
 \ifinner\nondmatherr@\pagebreak\else\postdisplaypenalty-\@M\fi
 \else\ifvmode\removelastskip\break\else\vadjust{\break}\fi\fi}
\def\nopagebreak{\RIfM@
 \ifinner\nondmatherr@\nopagebreak\else\postdisplaypenalty\@M\fi
 \else\ifvmode\nobreak\else\vadjust{\nobreak}\fi\fi}
\def\nonvmodeerr@#1{\Err@{\string#1\space not allowed within a paragraph
 or in math}}
\def\vnonvmode@#1#2{\relaxnext@\DNii@{\ifx\next\par\DN@{#1}\else
 \DN@{#2}\fi\next@}%
 \ifvmode\DN@{#1}\else
 \DN@{\FN@\nextii@}\fi\next@}
\def\newpage{\vnonvmode@{\vfill\break}{\nonvmodeerr@\newpage}}
\def\smallpagebreak{\vnonvmode@\smallbreak{\nonvmodeerr@\smallpagebreak}}
\def\medpagebreak{\vnonvmode@\medbreak{\nonvmodeerr@\medpagebreak}}
\def\bigpagebreak{\vnonvmode@\bigbreak{\nonvmodeerr@\bigpagebreak}}
\def\NoBlackBoxes{\global\overfullrule\z@}
\def\BlackBoxes{\global\overfullrule5\p@}
\def\Invalid@#1{\def#1{\Err@{\Invalid@@\string#1}}}
\def\Invalid@@{Invalid use of }
\Invalid@\caption
\Invalid@\captionwidth
\newdimen\smallcaptionwidth@
\def\topspace{\mid@false\ins@}
\def\midspace{\mid@true\ins@}
\newif\ifmid@
\def\captionfont@{}
\def\ins@#1{\relaxnext@\allowbreak
 \smallcaptionwidth@\captionwidth@\gdef\thespace@{#1}%
 \DN@{\ifx\next\space@\DN@. {\FN@\nextii@}\else
  \DN@.{\FN@\nextii@}\fi\next@.}%
 \DNii@{\ifx\next\caption\DN@\caption{\FN@\nextiii@}%
  \else\let\next@\nextiv@\fi\next@}%
 \def\nextiv@{\vnonvmode@
  {\ifmid@\expandafter\midinsert\else\expandafter\topinsert\fi
   \vbox to\thespace@{}\endinsert}
  {\ifmid@\nonvmodeerr@\midspace\else\nonvmodeerr@\topspace\fi}}%
 \def\nextiii@{\ifx\next\captionwidth\expandafter\nextv@
  \else\expandafter\nextvi@\fi}%
 \def\nextv@\captionwidth##1##2{\smallcaptionwidth@##1\relax\nextvi@{##2}}%
 \def\nextvi@##1{\def\thecaption@{\captionfont@##1}%
  \DN@{\ifx\next\space@\DN@. {\FN@\nextvii@}\else
   \DN@.{\FN@\nextvii@}\fi\next@.}%
  \FN@\next@}%
 \def\nextvii@{\vnonvmode@
  {\ifmid@\expandafter\midinsert\else
  \expandafter\topinsert\fi\vbox to\thespace@{}\nobreak\smallskip
  \setboxz@h{\noindent\ignorespaces\thecaption@\unskip}%
  \ifdim\wdz@>\smallcaptionwidth@\centerline{\vbox{\hsize\smallcaptionwidth@
   \noindent\ignorespaces\thecaption@\unskip}}%
  \else\centerline{\boxz@}\fi\endinsert}
  {\ifmid@\nonvmodeerr@\midspace
  \else\nonvmodeerr@\topspace\fi}}%
 \FN@\next@}
\def\newcodes@{\catcode`\\=12 \catcode`\{=12 \catcode`\}=12 \catcode`\#=12
 \catcode`\%=12\relax}
\def\oldcodes@{\catcode`\\=0 \catcode`\{=1 \catcode`\}=2 \catcode`\#=6
 \catcode`\%=14\relax}
\def\comment{\newcodes@\endlinechar=10 \comment@}
{\lccode`\0=`\\
\lowercase{\gdef\comment@#1^^J{\comment@@#10endcomment\comment@@@}%
\gdef\comment@@#10endcomment{\FN@\comment@@@}%
\gdef\comment@@@#1\comment@@@{\ifx\next\comment@@@\let\next\comment@
 \else\def\next{\oldcodes@\endlinechar=`\^^M\relax}%
 \fi\next}}}
\def\pr@m@s{\ifx'\next\DN@##1{\prim@s}\else\let\next@\egroup\fi\next@}
\def\prime{{\null\prime@\null}}
\mathchardef\prime@="0230
\let\dsize\displaystyle

\let\ssize\scriptstyle

\def\,{\RIfM@\mskip\thinmuskip\relax\else\kern.16667em\fi}
\def\!{\RIfM@\mskip-\thinmuskip\relax\else\kern-.16667em\fi}
\let\thinspace\,
\let\negthinspace\!
\def\medspace{\RIfM@\mskip\medmuskip\relax\else\kern.222222em\fi}
\def\negmedspace{\RIfM@\mskip-\medmuskip\relax\else\kern-.222222em\fi}
\def\thickspace{\RIfM@\mskip\thickmuskip\relax\else\kern.27777em\fi}
\let\;\thickspace
\def\negthickspace{\RIfM@\mskip-\thickmuskip\relax\else
 \kern-.27777em\fi}
\atdef@,{\RIfM@\mskip.1\thinmuskip\else\leavevmode\null,\fi}
\atdef@!{\RIfM@\mskip-.1\thinmuskip\else\leavevmode\null!\fi}
\atdef@.{\RIfM@&&\else\leavevmode.\spacefactor3000 \fi}
\def\and{\DOTSB\;\mathchar"3026 \;}

\def\frac#1#2{{#1\over#2}}

\newdimen\ex@
\ex@.2326ex
\Invalid@\thickness
\def\thickfrac{\relaxnext@
 \DN@{\ifx\next\thickness\let\next@\nextii@\else
 \DN@{\nextii@\thickness1}\fi\next@}%
 \DNii@\thickness##1##2##3{{##2\above##1\ex@##3}}%
 \FN@\next@}

\def\thickfracwithdelims#1#2{\relaxnext@\def\ldelim@{#1}\def\rdelim@{#2}%
 \DN@{\ifx\next\thickness\let\next@\nextii@\else
 \DN@{\nextii@\thickness1}\fi\next@}%
 \DNii@\thickness##1##2##3{{##2\abovewithdelims
 \ldelim@\rdelim@##1\ex@##3}}%
 \FN@\next@}

\def\:{\nobreak\hskip.1111em\mathpunct{}\nonscript\mkern-\thinmuskip{:}\hskip
 .3333emplus.0555em\relax}
\def\snug{\unskip\kern-\mathsurround}
\def\topsmash{\top@true\bot@false\smash@}
\def\botsmash{\top@false\bot@true\smash@}
\newif\iftop@
\newif\ifbot@
\def\smash{\top@true\bot@true\smash@}
\def\smash@{\RIfM@\expandafter\mathpalette\expandafter\mathsm@sh\else
 \expandafter\makesm@sh\fi}
\def\finsm@sh{\iftop@\ht\z@\z@\fi\ifbot@\dp\z@\z@\fi\leavevmode\boxz@}
\def\LimitsOnSums{\global\let\slimits@\displaylimits}
\def\NoLimitsOnSums{\global\let\slimits@\nolimits}
\LimitsOnSums
\mathchardef\coprod@="1360       \def\coprod{\DOTSB\coprod@\slimits@}
\mathchardef\bigvee@="1357       \def\bigvee{\DOTSB\bigvee@\slimits@}
\mathchardef\bigwedge@="1356     \def\bigwedge{\DOTSB\bigwedge@\slimits@}
\mathchardef\biguplus@="1355     \def\biguplus{\DOTSB\biguplus@\slimits@}
\mathchardef\bigcap@="1354       \def\bigcap{\DOTSB\bigcap@\slimits@}
\mathchardef\bigcup@="1353       \def\bigcup{\DOTSB\bigcup@\slimits@}
\mathchardef\prod@="1351         \def\prod{\DOTSB\prod@\slimits@}
\mathchardef\sum@="1350          \def\sum{\DOTSB\sum@\slimits@}
\mathchardef\bigotimes@="134E    \def\bigotimes{\DOTSB\bigotimes@\slimits@}
\mathchardef\bigoplus@="134C     \def\bigoplus{\DOTSB\bigoplus@\slimits@}
\mathchardef\bigodot@="134A      \def\bigodot{\DOTSB\bigodot@\slimits@}
\mathchardef\bigsqcup@="1346     \def\bigsqcup{\DOTSB\bigsqcup@\slimits@}
\def\LimitsOnInts{\global\let\ilimits@\displaylimits}
\def\NoLimitsOnInts{\global\let\ilimits@\nolimits}
\NoLimitsOnInts
\def\int{\DOTSI\intop\ilimits@}
\def\oint{\DOTSI\ointop\ilimits@}
\def\intic@{\mathchoice{\hskip.5em}{\hskip.4em}{\hskip.4em}{\hskip.4em}}
\def\negintic@{\mathchoice
 {\hskip-.5em}{\hskip-.4em}{\hskip-.4em}{\hskip-.4em}}
\def\intkern@{\mathchoice{\!\!\!}{\!\!}{\!\!}{\!\!}}
\def\intdots@{\mathchoice{\plaincdots@}
 {{\cdotp}\mkern1.5mu{\cdotp}\mkern1.5mu{\cdotp}}
 {{\cdotp}\mkern1mu{\cdotp}\mkern1mu{\cdotp}}
 {{\cdotp}\mkern1mu{\cdotp}\mkern1mu{\cdotp}}}
\newcount\intno@
\def\iint{\DOTSI\intno@\tw@\FN@\ints@}
\def\iiint{\DOTSI\intno@\thr@@\FN@\ints@}
\def\iiiint{\DOTSI\intno@4 \FN@\ints@}
\def\idotsint{\DOTSI\intno@\z@\FN@\ints@}
\def\ints@{\findlimits@\ints@@}
\newif\iflimtoken@
\newif\iflimits@
\def\findlimits@{\limtoken@true\ifx\next\limits\limits@true
 \else\ifx\next\nolimits\limits@false\else
 \limtoken@false\ifx\ilimits@\nolimits\limits@false\else
 \ifinner\limits@false\else\limits@true\fi\fi\fi\fi}
\def\multint@{\int\ifnum\intno@=\z@\intdots@                                
 \else\intkern@\fi                                                          
 \ifnum\intno@>\tw@\int\intkern@\fi                                         
 \ifnum\intno@>\thr@@\int\intkern@\fi                                       
 \int}                                                                      
\def\multintlimits@{\intop\ifnum\intno@=\z@\intdots@\else\intkern@\fi
 \ifnum\intno@>\tw@\intop\intkern@\fi
 \ifnum\intno@>\thr@@\intop\intkern@\fi\intop}
\def\ints@@{\iflimtoken@                                                    
 \def\ints@@@{\iflimits@\negintic@\mathop{\intic@\multintlimits@}\limits    
  \else\multint@\nolimits\fi                                                
  \eat@}                                                                    
 \else                                                                      
 \def\ints@@@{\iflimits@\negintic@
  \mathop{\intic@\multintlimits@}\limits\else
  \multint@\nolimits\fi}\fi\ints@@@}
\def\LimitsOnNames{\global\let\nlimits@\displaylimits}
\def\NoLimitsOnNames{\global\let\nlimits@\nolimits@}
\LimitsOnNames
\def\nolimits@{\relaxnext@
 \DN@{\ifx\next\limits\DN@\limits{\nolimits}\else
  \let\next@\nolimits\fi\next@}%
 \FN@\next@}
\def\newmcodes@{\mathcode`\'="0027 \mathcode`\*="002A \mathcode`\.="613A
 \mathcode`\-="002D \mathcode`\/="002F \mathcode`\:="603A }
\def\operatorname#1{\mathop{\newmcodes@\kern\z@\fam\z@#1}\nolimits@}
\def\operatornamewithlimits#1{\mathop{\newmcodes@\kern\z@\fam\z@#1}\nlimits@}
\def\qopname@#1{\mathop{\fam\z@#1}\nolimits@}
\def\qopnamewl@#1{\mathop{\fam\z@#1}\nlimits@}
\def\arccos{\qopname@{arccos}}
\def\arcsin{\qopname@{arcsin}}
\def\arctan{\qopname@{arctan}}
\def\arg{\qopname@{arg}}
\def\cos{\qopname@{cos}}
\def\cosh{\qopname@{cosh}}
\def\cot{\qopname@{cot}}
\def\coth{\qopname@{coth}}
\def\csc{\qopname@{csc}}
\def\deg{\qopname@{deg}}
\def\det{\qopnamewl@{det}}
\def\dim{\qopname@{dim}}
\def\exp{\qopname@{exp}}
\def\gcd{\qopnamewl@{gcd}}
\def\hom{\qopname@{hom}}
\def\inf{\qopnamewl@{inf}}
\def\injlim{\qopnamewl@{inj\,lim}}
\def\ker{\qopname@{ker}}
\def\lg{\qopname@{lg}}
\def\lim{\qopnamewl@{lim}}
\def\liminf{\qopnamewl@{lim\,inf}}
\def\limsup{\qopnamewl@{lim\,sup}}
\def\ln{\qopname@{ln}}
\def\log{\qopname@{log}}
\def\max{\qopnamewl@{max}}
\def\min{\qopnamewl@{min}}
\def\Pr{\qopnamewl@{Pr}}
\def\projlim{\qopnamewl@{proj\,lim}}
\def\sec{\qopname@{sec}}
\def\sin{\qopname@{sin}}
\def\sinh{\qopname@{sinh}}
\def\sup{\qopnamewl@{sup}}
\def\tan{\qopname@{tan}}
\def\tanh{\qopname@{tanh}}
\def\varinjlim{\mathop{\vtop{\ialign{##\crcr
 \hfil\rm lim\hfil\crcr\noalign{\nointerlineskip}\rightarrowfill\crcr
 \noalign{\nointerlineskip\kern-\ex@}\crcr}}}}
\def\varprojlim{\mathop{\vtop{\ialign{##\crcr
 \hfil\rm lim\hfil\crcr\noalign{\nointerlineskip}\leftarrowfill\crcr
 \noalign{\nointerlineskip\kern-\ex@}\crcr}}}}
\def\varliminf{\mathop{\underline{\vrule height\z@ depth.2exwidth\z@
 \hbox{\rm lim}}}}

\newdimen\buffer@
\buffer@\fontdimen13 \tenex
\newdimen\buffer
\buffer\buffer@

\def\ResetBuffer{\fontdimen13 \tenex\buffer@\global\buffer\buffer@}
\def\shave#1{\mathop{\hbox{$\m@th\fontdimen13 \tenex\z@                     
 \displaystyle{#1}$}}\fontdimen13 \tenex\buffer}

\Invalid@\\
\def\Let@{\relax\iffalse{\fi\let\\=\cr\iffalse}\fi}
\Invalid@\vspace
\def\vspace@{\def\vspace##1{\crcr\noalign{\vskip##1\relax}}}
\def\multilimits@{\bgroup\vspace@\Let@
 \baselineskip\fontdimen10 \scriptfont\tw@
 \advance\baselineskip\fontdimen12 \scriptfont\tw@
 \lineskip\thr@@\fontdimen8 \scriptfont\thr@@
 \lineskiplimit\lineskip
 \vbox\bgroup\ialign\bgroup\hfil$\m@th\scriptstyle{##}$\hfil\crcr}
\def\Sb{_\multilimits@}
\def\endSb{\crcr\egroup\egroup\egroup}
\def\Sp{^\multilimits@}

\def\spreadlines#1{\RIfMIfI@\onlydmatherr@\spreadlines\else
 \openup#1\relax\fi\else\onlydmatherr@\spreadlines\fi}
\def\Mathstrut@{\copy\Mathstrutbox@}
\newbox\Mathstrutbox@
\setbox\Mathstrutbox@\null
\setboxz@h{$\m@th($}
\ht\Mathstrutbox@\ht\z@
\dp\Mathstrutbox@\dp\z@
\newdimen\spreadmlines@
\def\spreadmatrixlines#1{\RIfMIfI@
 \onlydmatherr@\spreadmatrixlines\else
 \spreadmlines@#1\relax\fi\else\onlydmatherr@\spreadmatrixlines\fi}
\def\matrix{\null\,\vcenter\bgroup\Let@\vspace@
 \normalbaselines\openup\spreadmlines@\ialign
 \bgroup\hfil$\m@th##$\hfil&&\quad\hfil$\m@th##$\hfil\crcr
 \Mathstrut@\crcr\noalign{\kern-\baselineskip}}
\def\endmatrix{\crcr\Mathstrut@\crcr\noalign{\kern-\baselineskip}\egroup
 \egroup\,}
\def\format{\crcr\egroup\iffalse{\fi\ifnum`}=0 \fi\format@}
\newtoks\hashtoks@
\hashtoks@{#}
\def\format@#1\\{\def\preamble@{#1}%
 \def\l{$\m@th\the\hashtoks@$\hfil}%
 \def\c{\hfil$\m@th\the\hashtoks@$\hfil}%
 \def\r{\hfil$\m@th\the\hashtoks@$}%
 \edef\Preamble@{\preamble@}\ifnum`{=0 \fi\iffalse}\fi
 \ialign\bgroup\span\Preamble@\crcr}
\def\smallmatrix{\null\,\vcenter\bgroup\vspace@\Let@
 \baselineskip9\ex@\lineskip\ex@
 \ialign\bgroup\hfil$\m@th\scriptstyle{##}$\hfil&&\thickspace\hfil
 $\m@th\scriptstyle{##}$\hfil\crcr}
\def\endsmallmatrix{\crcr\egroup\egroup\,}

\newmuskip\dotsspace@
\dotsspace@1.5mu
\def\strip@#1 {#1}
\def\spacehdots#1\for#2{\multispan{#2}\xleaders
 \hbox{$\m@th\mkern\strip@#1 \dotsspace@.\mkern\strip@#1 \dotsspace@$}\hfill}
\def\hdotsfor#1{\spacehdots\@ne\for{#1}}
\def\multispan@#1{\omit\mscount#1\unskip\loop\ifnum\mscount>\@ne\sp@n\repeat}
\def\spaceinnerhdots#1\for#2\after#3{\multispan@{\strip@#2 }#3\xleaders
 \hbox{$\m@th\mkern\strip@#1 \dotsspace@.\mkern\strip@#1 \dotsspace@$}\hfill}
\def\innerhdotsfor#1\after#2{\spaceinnerhdots\@ne\for#1\after{#2}}
\def\cases{\bgroup\spreadmlines@\jot\left\{\,\matrix\format\l&\quad\l\\}
\def\endcases{\endmatrix\right.\egroup}
\newif\ifinany@
\newif\ifinalign@
\newif\ifingather@
\def\strut@{\copy\strutbox@}
\newbox\strutbox@
\setbox\strutbox@\hbox{\vrule height8\p@ depth3\p@ width\z@}
\def\topaligned{\null\,\vtop\aligned@}
\def\botaligned{\null\,\vbox\aligned@}
\def\aligned{\null\,\vcenter\aligned@}
\def\aligned@{\bgroup\vspace@\Let@
 \ifinany@\else\openup\jot\fi\ialign
 \bgroup\hfil\strut@$\m@th\displaystyle{##}$&
 $\m@th\displaystyle{{}##}$\hfil\crcr}
\def\endaligned{\crcr\egroup\egroup}

\def\alignedat#1{\null\,\vcenter\bgroup\doat@{#1}\vspace@\Let@
 \ifinany@\else\openup\jot\fi\ialign\bgroup\span\preamble@@\crcr}
\newcount\atcount@
\def\doat@#1{\toks@{\hfil\strut@$\m@th
 \displaystyle{\the\hashtoks@}$&$\m@th\displaystyle
 {{}\the\hashtoks@}$\hfil}
 \atcount@#1\relax\advance\atcount@\m@ne                                    
 \loop\ifnum\atcount@>\z@\toks@=\expandafter{\the\toks@&\hfil$\m@th
 \displaystyle{\the\hashtoks@}$&$\m@th
 \displaystyle{{}\the\hashtoks@}$\hfil}\advance
  \atcount@\m@ne\repeat                                                     
 \xdef\preamble@{\the\toks@}\xdef\preamble@@{\preamble@}}

\def\gathered{\null\,\vcenter\bgroup\vspace@\Let@
 \ifinany@\else\openup\jot\fi\ialign
 \bgroup\hfil\strut@$\m@th\displaystyle{##}$\hfil\crcr}
\def\endgathered{\crcr\egroup\egroup}
\newif\iftagsleft@
\def\TagsOnLeft{\global\tagsleft@true}
\def\TagsOnRight{\global\tagsleft@false}
\TagsOnLeft
\newif\ifmathtags@
\def\TagsAsMath{\global\mathtags@true}
\def\TagsAsText{\global\mathtags@false}
\TagsAsText
\def\tagform@#1{\hbox{\rm(\ignorespaces#1\unskip)}}
\def\thetag{\leavevmode\tagform@}
\def\tag#1$${\iftagsleft@\leqno\else\eqno\fi                                
 \maketag@#1\maketag@                                                       
 $$}                                                                        
\def\maketag@{\FN@\maketag@@}
\def\maketag@@{\ifx\next"\expandafter\maketag@@@\else\expandafter\maketag@@@@
 \fi}
\def\maketag@@@"#1"#2\maketag@{\hbox{\rm#1}}                                
\def\maketag@@@@#1\maketag@{\ifmathtags@\tagform@{$\m@th#1$}\else
 \tagform@{#1}\fi}
\interdisplaylinepenalty\@M
\def\allowdisplaybreaks{\RIfMIfI@
 \onlydmatherr@\allowdisplaybreaks\else
 \interdisplaylinepenalty\z@\fi\else\onlydmatherr@\allowdisplaybreaks\fi}
\Invalid@\allowdisplaybreak
\Invalid@\displaybreak
\Invalid@\intertext
\def\allowdisplaybreak@{\def\allowdisplaybreak{\crcr\noalign{\allowbreak}}}
\def\displaybreak@{\def\displaybreak{\crcr\noalign{\break}}}
\def\intertext@{\def\intertext##1{\crcr\noalign{\vskip\belowdisplayskip
 \vbox{\normalbaselines\noindent##1}\vskip\abovedisplayskip}}}
\newskip\centering@
\centering@\z@ plus\@m\p@
\def\align{\relax\ifingather@\DN@{\csname align (in
  \string\gather)\endcsname}\else
 \ifmmode\ifinner\DN@{\onlydmatherr@\align}\else
  \let\next@\align@\fi
 \else\DN@{\onlydmatherr@\align}\fi\fi\next@}
\newhelp\andhelp@
{An extra & here is so disastrous that you should probably exit^^J
and fix things up.}
\newif\iftag@
\newcount\and@
\def\align@{\inalign@true\inany@true
 \vspace@\allowdisplaybreak@\displaybreak@\intertext@
 \def\tag{\global\tag@true\ifnum\and@=\z@\DN@{&&}\else
          \DN@{&}\fi\next@}%
 \iftagsleft@\DN@{\csname align \endcsname}\else
  \DN@{\csname align \space\endcsname}\fi\next@}
\def\Tag@{\iftag@\else\errhelp\andhelp@\err@{Extra & on this line}\fi}
\newdimen\lwidth@
\newdimen\rwidth@
\newdimen\maxlwidth@
\newdimen\maxrwidth@
\newdimen\totwidth@
\def\measure@#1\endalign{\lwidth@\z@\rwidth@\z@\maxlwidth@\z@\maxrwidth@\z@
 \global\and@\z@                                                            
 \setbox@ne\vbox                                                            
  {\everycr{\noalign{\global\tag@false\global\and@\z@}}\Let@                
  \halign{\setboxz@h{$\m@th\displaystyle{\@lign##}$}
   \global\lwidth@\wdz@                                                     
   \ifdim\lwidth@>\maxlwidth@\global\maxlwidth@\lwidth@\fi                  
   \global\advance\and@\@ne                                                 
   &\setboxz@h{$\m@th\displaystyle{{}\@lign##}$}\global\rwidth@\wdz@        
   \ifdim\rwidth@>\maxrwidth@\global\maxrwidth@\rwidth@\fi                  
   \global\advance\and@\@ne                                                
   &\Tag@
   \eat@{##}\crcr#1\crcr}}
 \totwidth@\maxlwidth@\advance\totwidth@\maxrwidth@}                       
\def\displ@y@{\global\dt@ptrue\openup\jot
 \everycr{\noalign{\global\tag@false\global\and@\z@\ifdt@p\global\dt@pfalse
 \vskip-\lineskiplimit\vskip\normallineskiplimit\else
 \penalty\interdisplaylinepenalty\fi}}}
\def\black@#1{\noalign{\ifdim#1>\displaywidth
 \dimen@\prevdepth\nointerlineskip                                          
 \vskip-\ht\strutbox@\vskip-\dp\strutbox@                                   
 \vbox{\noindent\hbox to#1{\strut@\hfill}}
 \prevdepth\dimen@                                                          
 \fi}}
\expandafter\def\csname align \space\endcsname#1\endalign
 {\measure@#1\endalign\global\and@\z@                                       
 \ifingather@\everycr{\noalign{\global\and@\z@}}\else\displ@y@\fi           
 \Let@\tabskip\centering@                                                   
 \halign to\displaywidth
  {\hfil\strut@\setboxz@h{$\m@th\displaystyle{\@lign##}$}
  \global\lwidth@\wdz@\boxz@\global\advance\and@\@ne                        
  \tabskip\z@skip                                                           
  &\setboxz@h{$\m@th\displaystyle{{}\@lign##}$}
  \global\rwidth@\wdz@\boxz@\hfill\global\advance\and@\@ne                  
  \tabskip\centering@                                                       
  &\setboxz@h{\@lign\strut@\maketag@##\maketag@}
  \dimen@\displaywidth\advance\dimen@-\totwidth@
  \divide\dimen@\tw@\advance\dimen@\maxrwidth@\advance\dimen@-\rwidth@     
  \ifdim\dimen@<\tw@\wdz@\llap{\vtop{\normalbaselines\null\boxz@}}
  \else\llap{\boxz@}\fi                                                    
  \tabskip\z@skip                                                          
  \crcr#1\crcr                                                             
  \black@\totwidth@}}                                                      
\newdimen\lineht@
\expandafter\def\csname align \endcsname#1\endalign{\measure@#1\endalign
 \global\and@\z@
 \ifdim\totwidth@>\displaywidth\let\displaywidth@\totwidth@\else
  \let\displaywidth@\displaywidth\fi                                        
 \ifingather@\everycr{\noalign{\global\and@\z@}}\else\displ@y@\fi
 \Let@\tabskip\centering@\halign to\displaywidth
  {\hfil\strut@\setboxz@h{$\m@th\displaystyle{\@lign##}$}%
  \global\lwidth@\wdz@\global\lineht@\ht\z@                                 
  \boxz@\global\advance\and@\@ne
  \tabskip\z@skip&\setboxz@h{$\m@th\displaystyle{{}\@lign##}$}%
  \global\rwidth@\wdz@\ifdim\ht\z@>\lineht@\global\lineht@\ht\z@\fi         
  \boxz@\hfil\global\advance\and@\@ne
  \tabskip\centering@&\kern-\displaywidth@                                  
  \setboxz@h{\@lign\strut@\maketag@##\maketag@}%
  \dimen@\displaywidth\advance\dimen@-\totwidth@
  \divide\dimen@\tw@\advance\dimen@\maxlwidth@\advance\dimen@-\lwidth@
  \ifdim\dimen@<\tw@\wdz@
   \rlap{\vbox{\normalbaselines\boxz@\vbox to\lineht@{}}}\else
   \rlap{\boxz@}\fi
  \tabskip\displaywidth@\crcr#1\crcr\black@\totwidth@}}
\expandafter\def\csname align (in \string\gather)\endcsname
  #1\endalign{\vcenter{\align@#1\endalign}}
\Invalid@\endalign
\newif\ifxat@
\def\alignat{\RIfMIfI@\DN@{\onlydmatherr@\alignat}\else
 \DN@{\csname alignat \endcsname}\fi\else
 \DN@{\onlydmatherr@\alignat}\fi\next@}
\newif\ifmeasuring@
\newbox\savealignat@
\expandafter\def\csname alignat \endcsname#1#2\endalignat                   
 {\inany@true\xat@false
 \def\tag{\global\tag@true\count@#1\relax\multiply\count@\tw@
  \xdef\tag@{}\loop\ifnum\count@>\and@\xdef\tag@{&\tag@}\advance\count@\m@ne
  \repeat\tag@}%
 \vspace@\allowdisplaybreak@\displaybreak@\intertext@
 \displ@y@\measuring@true                                                   
 \setbox\savealignat@\hbox{$\m@th\displaystyle\Let@
  \attag@{#1}
  \vbox{\halign{\span\preamble@@\crcr#2\crcr}}$}%
 \measuring@false                                                           
 \Let@\attag@{#1}
 \tabskip\centering@\halign to\displaywidth
  {\span\preamble@@\crcr#2\crcr                                             
  \black@{\wd\savealignat@}}}                                               
\Invalid@\endalignat
\def\xalignat{\RIfMIfI@
 \DN@{\onlydmatherr@\xalignat}\else
 \DN@{\csname xalignat \endcsname}\fi\else
 \DN@{\onlydmatherr@\xalignat}\fi\next@}
\expandafter\def\csname xalignat \endcsname#1#2\endxalignat
 {\inany@true\xat@true
 \def\tag{\global\tag@true\def\tag@{}\count@#1\relax\multiply\count@\tw@
  \loop\ifnum\count@>\and@\xdef\tag@{&\tag@}\advance\count@\m@ne\repeat\tag@}%
 \vspace@\allowdisplaybreak@\displaybreak@\intertext@
 \displ@y@\measuring@true\setbox\savealignat@\hbox{$\m@th\displaystyle\Let@
 \attag@{#1}\vbox{\halign{\span\preamble@@\crcr#2\crcr}}$}%
 \measuring@false\Let@
 \attag@{#1}\tabskip\centering@\halign to\displaywidth
 {\span\preamble@@\crcr#2\crcr\black@{\wd\savealignat@}}}
\def\attag@#1{\let\Maketag@\maketag@\let\TAG@\Tag@                          
 \let\Tag@=0\let\maketag@=0
 \ifmeasuring@\def\llap@##1{\setboxz@h{##1}\hbox to\tw@\wdz@{}}%
  \def\rlap@##1{\setboxz@h{##1}\hbox to\tw@\wdz@{}}\else
  \let\llap@\llap\let\rlap@\rlap\fi                                         
 \toks@{\hfil\strut@$\m@th\displaystyle{\@lign\the\hashtoks@}$\tabskip\z@skip
  \global\advance\and@\@ne&$\m@th\displaystyle{{}\@lign\the\hashtoks@}$\hfil
  \ifxat@\tabskip\centering@\fi\global\advance\and@\@ne}
 \iftagsleft@
  \toks@@{\tabskip\centering@&\Tag@\kern-\displaywidth
   \rlap@{\@lign\maketag@\the\hashtoks@\maketag@}%
   \global\advance\and@\@ne\tabskip\displaywidth}\else
  \toks@@{\tabskip\centering@&\Tag@\llap@{\@lign\maketag@
   \the\hashtoks@\maketag@}\global\advance\and@\@ne\tabskip\z@skip}\fi      
 \atcount@#1\relax\advance\atcount@\m@ne
 \loop\ifnum\atcount@>\z@
 \toks@=\expandafter{\the\toks@&\hfil$\m@th\displaystyle{\@lign
  \the\hashtoks@}$\global\advance\and@\@ne
  \tabskip\z@skip&$\m@th\displaystyle{{}\@lign\the\hashtoks@}$\hfil\ifxat@
  \tabskip\centering@\fi\global\advance\and@\@ne}\advance\atcount@\m@ne
 \repeat                                                                    
 \xdef\preamble@{\the\toks@\the\toks@@}
 \xdef\preamble@@{\preamble@}
 \let\maketag@\Maketag@\let\Tag@\TAG@}                                      
\Invalid@\endxalignat
\def\xxalignat{\RIfMIfI@
 \DN@{\onlydmatherr@\xxalignat}\else\DN@{\csname xxalignat
  \endcsname}\fi\else
 \DN@{\onlydmatherr@\xxalignat}\fi\next@}
\expandafter\def\csname xxalignat \endcsname#1#2\endxxalignat{\inany@true
 \vspace@\allowdisplaybreak@\displaybreak@\intertext@
 \displ@y\setbox\savealignat@\hbox{$\m@th\displaystyle\Let@
 \xxattag@{#1}\vbox{\halign{\span\preamble@@\crcr#2\crcr}}$}%
 \Let@\xxattag@{#1}\tabskip\z@skip\halign to\displaywidth
 {\span\preamble@@\crcr#2\crcr\black@{\wd\savealignat@}}}
\def\xxattag@#1{\toks@{\tabskip\z@skip\hfil\strut@
 $\m@th\displaystyle{\the\hashtoks@}$&%
 $\m@th\displaystyle{{}\the\hashtoks@}$\hfil\tabskip\centering@&}%
 \atcount@#1\relax\advance\atcount@\m@ne\loop\ifnum\atcount@>\z@
 \toks@=\expandafter{\the\toks@&\hfil$\m@th\displaystyle{\the\hashtoks@}$%
  \tabskip\z@skip&$\m@th\displaystyle{{}\the\hashtoks@}$\hfil
  \tabskip\centering@}\advance\atcount@\m@ne\repeat
 \xdef\preamble@{\the\toks@\tabskip\z@skip}\xdef\preamble@@{\preamble@}}
\Invalid@\endxxalignat
\newdimen\gwidth@
\newdimen\gmaxwidth@
\def\gmeasure@#1\endgather{\gwidth@\z@\gmaxwidth@\z@\setbox@ne\vbox{\Let@
 \halign{\setboxz@h{$\m@th\displaystyle{##}$}\global\gwidth@\wdz@
 \ifdim\gwidth@>\gmaxwidth@\global\gmaxwidth@\gwidth@\fi
 &\eat@{##}\crcr#1\crcr}}}
\def\gather{\RIfMIfI@\DN@{\onlydmatherr@\gather}\else
 \ingather@true\inany@true\def\tag{&}%
 \vspace@\allowdisplaybreak@\displaybreak@\intertext@
 \displ@y\Let@
 \iftagsleft@\DN@{\csname gather \endcsname}\else
  \DN@{\csname gather \space\endcsname}\fi\fi
 \else\DN@{\onlydmatherr@\gather}\fi\next@}
\expandafter\def\csname gather \space\endcsname#1\endgather
 {\gmeasure@#1\endgather\tabskip\centering@
 \halign to\displaywidth{\hfil\strut@\setboxz@h{$\m@th\displaystyle{##}$}%
 \global\gwidth@\wdz@\boxz@\hfil&
 \setboxz@h{\strut@{\maketag@##\maketag@}}%
 \dimen@\displaywidth\advance\dimen@-\gwidth@
 \ifdim\dimen@>\tw@\wdz@\llap{\boxz@}\else
 \llap{\vtop{\normalbaselines\null\boxz@}}\fi
 \tabskip\z@skip\crcr#1\crcr\black@\gmaxwidth@}}
\newdimen\glineht@
\expandafter\def\csname gather \endcsname#1\endgather{\gmeasure@#1\endgather
 \ifdim\gmaxwidth@>\displaywidth\let\gdisplaywidth@\gmaxwidth@\else
 \let\gdisplaywidth@\displaywidth\fi\tabskip\centering@\halign to\displaywidth
 {\hfil\strut@\setboxz@h{$\m@th\displaystyle{##}$}%
 \global\gwidth@\wdz@\global\glineht@\ht\z@\boxz@\hfil&\kern-\gdisplaywidth@
 \setboxz@h{\strut@{\maketag@##\maketag@}}%
 \dimen@\displaywidth\advance\dimen@-\gwidth@
 \ifdim\dimen@>\tw@\wdz@\rlap{\boxz@}\else
 \rlap{\vbox{\normalbaselines\boxz@\vbox to\glineht@{}}}\fi
 \tabskip\gdisplaywidth@\crcr#1\crcr\black@\gmaxwidth@}}
\newif\ifctagsplit@
\def\CenteredTagsOnSplits{\global\ctagsplit@true}
\def\TopOrBottomTagsOnSplits{\global\ctagsplit@false}
\TopOrBottomTagsOnSplits
\def\split{\relax\ifinany@\let\next@\insplit@\else
 \ifmmode\ifinner\def\next@{\onlydmatherr@\split}\else
 \let\next@\outsplit@\fi\else
 \def\next@{\onlydmatherr@\split}\fi\fi\next@}
\def\insplit@{\global\setbox\z@\vbox\bgroup\vspace@\Let@\ialign\bgroup
 \hfil\strut@$\m@th\displaystyle{##}$&$\m@th\displaystyle{{}##}$\hfill\crcr}
\def\endsplit{\crcr\egroup\egroup\iftagsleft@\expandafter\lendsplit@\else
 \expandafter\rendsplit@\fi}
\def\rendsplit@{\global\setbox9 \vbox
 {\unvcopy\z@\global\setbox8 \lastbox\unskip}
 \setbox@ne\hbox{\unhcopy8 \unskip\global\setbox\tw@\lastbox
 \unskip\global\setbox\thr@@\lastbox}
 \global\setbox7 \hbox{\unhbox\tw@\unskip}
 \ifinalign@\ifctagsplit@                                                   
  \gdef\split@{\hbox to\wd\thr@@{}&
   \vcenter{\vbox{\moveleft\wd\thr@@\boxz@}}}
 \else\gdef\split@{&\vbox{\moveleft\wd\thr@@\box9}\crcr
  \box\thr@@&\box7}\fi                                                      
 \else                                                                      
  \ifctagsplit@\gdef\split@{\vcenter{\boxz@}}\else
  \gdef\split@{\box9\crcr\hbox{\box\thr@@\box7}}\fi
 \fi
 \split@}                                                                   
\def\lendsplit@{\global\setbox9\vtop{\unvcopy\z@}
 \setbox@ne\vbox{\unvcopy\z@\global\setbox8\lastbox}
 \setbox@ne\hbox{\unhcopy8\unskip\setbox\tw@\lastbox
  \unskip\global\setbox\thr@@\lastbox}
 \ifinalign@\ifctagsplit@                                                   
  \gdef\split@{\hbox to\wd\thr@@{}&
  \vcenter{\vbox{\moveleft\wd\thr@@\box9}}}
  \else                                                                     
  \gdef\split@{\hbox to\wd\thr@@{}&\vbox{\moveleft\wd\thr@@\box9}}\fi
 \else
  \ifctagsplit@\gdef\split@{\vcenter{\box9}}\else
  \gdef\split@{\box9}\fi
 \fi\split@}
\def\outsplit@#1$${\align\insplit@#1\endalign$$}
\newdimen\multlinegap@
\multlinegap@1em
\newdimen\multlinetaggap@
\multlinetaggap@1em
\def\MultlineGap#1{\global\multlinegap@#1\relax}
\def\multlinegap#1{\RIfMIfI@\onlydmatherr@\multlinegap\else
 \multlinegap@#1\relax\fi\else\onlydmatherr@\multlinegap\fi}
\def\nomultlinegap{\multlinegap{\z@}}
\def\multline{\RIfMIfI@
 \DN@{\onlydmatherr@\multline}\else
 \DN@{\multline@}\fi\else
 \DN@{\onlydmatherr@\multline}\fi\next@}
\newif\iftagin@
\def\tagin@#1{\tagin@false\in@\tag{#1}\ifin@\tagin@true\fi}
\def\multline@#1$${\inany@true\vspace@\allowdisplaybreak@\displaybreak@
 \tagin@{#1}\iftagsleft@\DN@{\multline@l#1$$}\else
 \DN@{\multline@r#1$$}\fi\next@}
\newdimen\mwidth@
\def\rmmeasure@#1\endmultline{%
 \def\shoveleft##1{##1}\def\shoveright##1{##1}
 \setbox@ne\vbox{\Let@\halign{\setboxz@h
  {$\m@th\@lign\displaystyle{}##$}\global\mwidth@\wdz@
  \crcr#1\crcr}}}
\newdimen\mlineht@
\newif\ifzerocr@
\newif\ifonecr@
\def\lmmeasure@#1\endmultline{\global\zerocr@true\global\onecr@false
 \everycr{\noalign{\ifonecr@\global\onecr@false\fi
  \ifzerocr@\global\zerocr@false\global\onecr@true\fi}}
  \def\shoveleft##1{##1}\def\shoveright##1{##1}%
 \setbox@ne\vbox{\Let@\halign{\setboxz@h
  {$\m@th\@lign\displaystyle{}##$}\ifonecr@\global\mwidth@\wdz@
  \global\mlineht@\ht\z@\fi\crcr#1\crcr}}}
\newbox\mtagbox@
\newdimen\ltwidth@
\newdimen\rtwidth@
\def\multline@l#1$${\iftagin@\DN@{\lmultline@@#1$$}\else
 \DN@{\setbox\mtagbox@\null\ltwidth@\z@\rtwidth@\z@
  \lmultline@@@#1$$}\fi\next@}
\def\lmultline@@#1\endmultline\tag#2$${%
 \setbox\mtagbox@\hbox{\maketag@#2\maketag@}
 \lmmeasure@#1\endmultline\dimen@\mwidth@\advance\dimen@\wd\mtagbox@
 \advance\dimen@\multlinetaggap@                                            
 \ifdim\dimen@>\displaywidth\ltwidth@\z@\else\ltwidth@\wd\mtagbox@\fi       
 \lmultline@@@#1\endmultline$$}
\def\lmultline@@@{\displ@y
 \def\shoveright##1{##1\hfilneg\hskip\multlinegap@}%
 \def\shoveleft##1{\setboxz@h{$\m@th\displaystyle{}##1$}%
  \setbox@ne\hbox{$\m@th\displaystyle##1$}%
  \hfilneg
  \iftagin@
   \ifdim\ltwidth@>\z@\hskip\ltwidth@\hskip\multlinetaggap@\fi
  \else\hskip\multlinegap@\fi\hskip.5\wd@ne\hskip-.5\wdz@##1}
  \halign\bgroup\Let@\hbox to\displaywidth
   {\strut@$\m@th\displaystyle\hfil{}##\hfil$}\crcr
   \hfilneg                                                                 
   \iftagin@                                                                
    \ifdim\ltwidth@>\z@                                                     
     \box\mtagbox@\hskip\multlinetaggap@                                    
    \else
     \rlap{\vbox{\normalbaselines\hbox{\strut@\box\mtagbox@}%
     \vbox to\mlineht@{}}}\fi                                               
   \else\hskip\multlinegap@\fi}                                             
\def\multline@r#1$${\iftagin@\DN@{\rmultline@@#1$$}\else
 \DN@{\setbox\mtagbox@\null\ltwidth@\z@\rtwidth@\z@
  \rmultline@@@#1$$}\fi\next@}
\def\rmultline@@#1\endmultline\tag#2$${\ltwidth@\z@
 \setbox\mtagbox@\hbox{\maketag@#2\maketag@}%
 \rmmeasure@#1\endmultline\dimen@\mwidth@\advance\dimen@\wd\mtagbox@
 \advance\dimen@\multlinetaggap@
 \ifdim\dimen@>\displaywidth\rtwidth@\z@\else\rtwidth@\wd\mtagbox@\fi
 \rmultline@@@#1\endmultline$$}
\def\rmultline@@@{\displ@y
 \def\shoveright##1{##1\hfilneg\iftagin@\ifdim\rtwidth@>\z@
  \hskip\rtwidth@\hskip\multlinetaggap@\fi\else\hskip\multlinegap@\fi}%
 \def\shoveleft##1{\setboxz@h{$\m@th\displaystyle{}##1$}%
  \setbox@ne\hbox{$\m@th\displaystyle##1$}%
  \hfilneg\hskip\multlinegap@\hskip.5\wd@ne\hskip-.5\wdz@##1}%
 \halign\bgroup\Let@\hbox to\displaywidth
  {\strut@$\m@th\displaystyle\hfil{}##\hfil$}\crcr
 \hfilneg\hskip\multlinegap@}
\def\endmultline{\iftagsleft@\expandafter\lendmultline@\else
 \expandafter\rendmultline@\fi}
\def\lendmultline@{\hfilneg\hskip\multlinegap@\crcr\egroup}
\def\rendmultline@{\iftagin@                                                
 \ifdim\rtwidth@>\z@                                                        
  \hskip\multlinetaggap@\box\mtagbox@                                       
 \else\llap{\vtop{\normalbaselines\null\hbox{\strut@\box\mtagbox@}}}\fi     
 \else\hskip\multlinegap@\fi                                                
 \hfilneg\crcr\egroup}
\def\bmod{\mskip-\medmuskip\mkern5mu\mathbin{\fam\z@ mod}\penalty900
 \mkern5mu\mskip-\medmuskip}
\def\pmod#1{\allowbreak\ifinner\mkern8mu\else\mkern18mu\fi
 ({\fam\z@ mod}\,\,#1)}
\def\pod#1{\allowbreak\ifinner\mkern8mu\else\mkern18mu\fi(#1)}
\def\mod#1{\allowbreak\ifinner\mkern12mu\else\mkern18mu\fi{\fam\z@ mod}\,\,#1}
\newcount\cfraccount@
\def\cfrac{\bgroup\bgroup\advance\cfraccount@\@ne\strut
 \iffalse{\fi\def\\{\over\displaystyle}\iffalse}\fi}
\def\lcfrac{\bgroup\bgroup\advance\cfraccount@\@ne\strut
 \iffalse{\fi\def\\{\hfill\over\displaystyle}\iffalse}\fi}
\def\rcfrac{\bgroup\bgroup\advance\cfraccount@\@ne\strut\hfill
 \iffalse{\fi\def\\{\over\displaystyle}\iffalse}\fi}
\def\gloop@#1\repeat{\gdef\body{#1}\iterate}
\def\endcfrac{\gloop@\ifnum\cfraccount@>\z@\global\advance\cfraccount@\m@ne
 \egroup\hskip-\nulldelimiterspace\egroup\repeat}
\def\binrel@#1{\setboxz@h{\thinmuskip0mu
  \medmuskip\m@ne mu\thickmuskip\@ne mu$#1\m@th$}%
 \setbox@ne\hbox{\thinmuskip0mu\medmuskip\m@ne mu\thickmuskip
  \@ne mu${}#1{}\m@th$}%
 \setbox\tw@\hbox{\hskip\wd@ne\hskip-\wdz@}}
\def\overset#1\to#2{\binrel@{#2}\ifdim\wd\tw@<\z@
 \mathbin{\mathop{\kern\z@#2}\limits^{#1}}\else\ifdim\wd\tw@>\z@
 \mathrel{\mathop{\kern\z@#2}\limits^{#1}}\else
 {\mathop{\kern\z@#2}\limits^{#1}}{}\fi\fi}
\def\underset#1\to#2{\binrel@{#2}\ifdim\wd\tw@<\z@
 \mathbin{\mathop{\kern\z@#2}\limits_{#1}}\else\ifdim\wd\tw@>\z@
 \mathrel{\mathop{\kern\z@#2}\limits_{#1}}\else
 {\mathop{\kern\z@#2}\limits_{#1}}{}\fi\fi}
\def\oversetbrace#1\to#2{\overbrace{#2}^{#1}}
\def\undersetbrace#1\to#2{\underbrace{#2}_{#1}}
\def\sideset#1\and#2\to#3{%
 \setbox@ne\hbox{$\dsize{\vphantom{#3}}#1{#3}\m@th$}%
 \setbox\tw@\hbox{$\dsize{#3}#2\m@th$}%
 \hskip\wd@ne\hskip-\wd\tw@\mathop{\hskip\wd\tw@\hskip-\wd@ne
  {\vphantom{#3}}#1{#3}#2}}
\def\rightarrowfill@#1{$#1\m@th\mathord-\mkern-6mu\cleaders
 \hbox{$#1\mkern-2mu\mathord-\mkern-2mu$}\hfill
 \mkern-6mu\mathord\rightarrow$}
\def\leftarrowfill@#1{$#1\m@th\mathord\leftarrow\mkern-6mu\cleaders
 \hbox{$#1\mkern-2mu\mathord-\mkern-2mu$}\hfill\mkern-6mu\mathord-$}
\def\leftrightarrowfill@#1{$#1\m@th\mathord\leftarrow\mkern-6mu\cleaders
 \hbox{$#1\mkern-2mu\mathord-\mkern-2mu$}\hfill
 \mkern-6mu\mathord\rightarrow$}
\def\overrightarrow{\mathpalette\overrightarrow@}
\def\overrightarrow@#1#2{\vbox{\ialign{##\crcr\rightarrowfill@#1\crcr
 \noalign{\kern-\ex@\nointerlineskip}$\m@th\hfil#1#2\hfil$\crcr}}}

\def\overleftarrow{\mathpalette\overleftarrow@}
\def\overleftarrow@#1#2{\vbox{\ialign{##\crcr\leftarrowfill@#1\crcr
 \noalign{\kern-\ex@\nointerlineskip}$\m@th\hfil#1#2\hfil$\crcr}}}
\def\overleftrightarrow{\mathpalette\overleftrightarrow@}
\def\overleftrightarrow@#1#2{\vbox{\ialign{##\crcr\leftrightarrowfill@#1\crcr
 \noalign{\kern-\ex@\nointerlineskip}$\m@th\hfil#1#2\hfil$\crcr}}}
\def\underrightarrow{\mathpalette\underrightarrow@}
\def\underrightarrow@#1#2{\vtop{\ialign{##\crcr$\m@th\hfil#1#2\hfil$\crcr
 \noalign{\nointerlineskip}\rightarrowfill@#1\crcr}}}

\def\underleftarrow{\mathpalette\underleftarrow@}
\def\underleftarrow@#1#2{\vtop{\ialign{##\crcr$\m@th\hfil#1#2\hfil$\crcr
 \noalign{\nointerlineskip}\leftarrowfill@#1\crcr}}}
\def\underleftrightarrow{\mathpalette\underleftrightarrow@}
\def\underleftrightarrow@#1#2{\vtop{\ialign{##\crcr$\m@th\hfil#1#2\hfil$\crcr
 \noalign{\nointerlineskip}\leftrightarrowfill@#1\crcr}}}
\let\DOTSI\relax
\let\DOTSB\relax

\newif\ifmath@
{\uccode`7=`\\ \uccode`8=`m \uccode`9=`a \uccode`0=`t \uccode`!=`h
 \uppercase{\gdef\math@#1#2#3#4#5#6\math@{\global\math@false\ifx 7#1\ifx 8#2%
 \ifx 9#3\ifx 0#4\ifx !#5\xdef\meaning@{#6}\global\math@true\fi\fi\fi\fi\fi}}}
\newif\ifmathch@
{\uccode`7=`c \uccode`8=`h \uccode`9=`\"
 \uppercase{\gdef\mathch@#1#2#3#4#5#6\mathch@{\global\mathch@false
  \ifx 7#1\ifx 8#2\ifx 9#5\global\mathch@true\xdef\meaning@{9#6}\fi\fi\fi}}}
\newcount\classnum@
\def\getmathch@#1.#2\getmathch@{\classnum@#1 \divide\classnum@4096
 \ifcase\number\classnum@\or\or\gdef\thedots@{\dotsb@}\or
 \gdef\thedots@{\dotsb@}\fi}
\newif\ifmathbin@
{\uccode`4=`b \uccode`5=`i \uccode`6=`n
 \uppercase{\gdef\mathbin@#1#2#3{\relaxnext@
  \DNii@##1\mathbin@{\ifx\space@\next\global\mathbin@true\fi}%
 \global\mathbin@false\DN@##1\mathbin@{}%
 \ifx 4#1\ifx 5#2\ifx 6#3\DN@{\FN@\nextii@}\fi\fi\fi\next@}}}
\newif\ifmathrel@
{\uccode`4=`r \uccode`5=`e \uccode`6=`l
 \uppercase{\gdef\mathrel@#1#2#3{\relaxnext@
  \DNii@##1\mathrel@{\ifx\space@\next\global\mathrel@true\fi}%
 \global\mathrel@false\DN@##1\mathrel@{}%
 \ifx 4#1\ifx 5#2\ifx 6#3\DN@{\FN@\nextii@}\fi\fi\fi\next@}}}
\newif\ifmacro@
{\uccode`5=`m \uccode`6=`a \uccode`7=`c
 \uppercase{\gdef\macro@#1#2#3#4\macro@{\global\macro@false
  \ifx 5#1\ifx 6#2\ifx 7#3\global\macro@true
  \xdef\meaning@{\macro@@#4\macro@@}\fi\fi\fi}}}
\def\macro@@#1->#2\macro@@{#2}
\newif\ifDOTS@
\newcount\DOTSCASE@
{\uccode`6=`\\ \uccode`7=`D \uccode`8=`O \uccode`9=`T \uccode`0=`S
 \uppercase{\gdef\DOTS@#1#2#3#4#5{\global\DOTS@false\DN@##1\DOTS@{}%
  \ifx 6#1\ifx 7#2\ifx 8#3\ifx 9#4\ifx 0#5\let\next@\DOTS@@\fi\fi\fi\fi\fi
  \next@}}}
{\uccode`3=`B \uccode`4=`I \uccode`5=`X
 \uppercase{\gdef\DOTS@@#1{\relaxnext@
  \DNii@##1\DOTS@{\ifx\space@\next\global\DOTS@true\fi}%
  \DN@{\FN@\nextii@}%
  \ifx 3#1\global\DOTSCASE@\z@\else
  \ifx 4#1\global\DOTSCASE@\@ne\else
  \ifx 5#1\global\DOTSCASE@\tw@\else\DN@##1\DOTS@{}%
  \fi\fi\fi\next@}}}
\newif\ifnot@
{\uccode`5=`\\ \uccode`6=`n \uccode`7=`o \uccode`8=`t
 \uppercase{\gdef\not@#1#2#3#4{\relaxnext@
  \DNii@##1\not@{\ifx\space@\next\global\not@true\fi}%
 \global\not@false\DN@##1\not@{}%
 \ifx 5#1\ifx 6#2\ifx 7#3\ifx 8#4\DN@{\FN@\nextii@}\fi\fi\fi
 \fi\next@}}}
\newif\ifkeybin@
\def\keybin@{\keybin@true
 \ifx\next+\else\ifx\next=\else\ifx\next<\else\ifx\next>\else\ifx\next-\else
 \ifx\next*\else\ifx\next:\else\keybin@false\fi\fi\fi\fi\fi\fi\fi}
\def\dots{\RIfM@\expandafter\mdots@\else\expandafter\tdots@\fi}
\def\tdots@{\unskip\relaxnext@
 \DN@{$\m@th\mathinner{\ldotp\ldotp\ldotp}\,
   \ifx\next,\,$\else\ifx\next.\,$\else\ifx\next;\,$\else\ifx\next:\,$\else
   \ifx\next?\,$\else\ifx\next!\,$\else$ \fi\fi\fi\fi\fi\fi}%
 \ \FN@\next@}
\def\mdots@{\FN@\mdots@@}
\def\mdots@@{\gdef\thedots@{\dotso@}
 \ifx\next\boldkey\gdef\thedots@\boldkey{\boldkeydots@}\else                
 \ifx\next\boldsymbol\gdef\thedots@\boldsymbol{\boldsymboldots@}\else       
 \ifx,\next\gdef\thedots@{\dotsc}
 \else\ifx\not\next\gdef\thedots@{\dotsb@}
 \else\keybin@
 \ifkeybin@\gdef\thedots@{\dotsb@}
 \else\xdef\meaning@{\meaning\next..........}\xdef\meaning@@{\meaning@}
  \expandafter\math@\meaning@\math@
  \ifmath@
   \expandafter\mathch@\meaning@\mathch@
   \ifmathch@\expandafter\getmathch@\meaning@\getmathch@\fi                 
  \else\expandafter\macro@\meaning@@\macro@                                 
  \ifmacro@                                                                
   \expandafter\not@\meaning@\not@\ifnot@\gdef\thedots@{\dotsb@}
  \else\expandafter\DOTS@\meaning@\DOTS@
  \ifDOTS@
   \ifcase\number\DOTSCASE@\gdef\thedots@{\dotsb@}%
    \or\gdef\thedots@{\dotsi}\else\fi                                      
  \else\expandafter\math@\meaning@\math@                                   
  \ifmath@\expandafter\mathbin@\meaning@\mathbin@
  \ifmathbin@\gdef\thedots@{\dotsb@}
  \else\expandafter\mathrel@\meaning@\mathrel@
  \ifmathrel@\gdef\thedots@{\dotsb@}
  \fi\fi\fi\fi\fi\fi\fi\fi\fi\fi\fi\fi
 \thedots@}
\def\plainldots@{\mathinner{\ldotp\ldotp\ldotp}}
\def\plaincdots@{\mathinner{\cdotp\cdotp\cdotp}}
\def\dotsi{\!\plaincdots@}
\let\dotsb@\plaincdots@
\newif\ifextra@
\newif\ifrightdelim@
\def\rightdelim@{\global\rightdelim@true                                    
 \ifx\next)\else                                                            
 \ifx\next]\else
 \ifx\next\rbrack\else
 \ifx\next\}\else
 \ifx\next\rbrace\else
 \ifx\next\rangle\else
 \ifx\next\rceil\else
 \ifx\next\rfloor\else
 \ifx\next\rgroup\else
 \ifx\next\rmoustache\else
 \ifx\next\right\else
 \ifx\next\bigr\else
 \ifx\next\biggr\else
 \ifx\next\Bigr\else                                                        
 \ifx\next\Biggr\else\global\rightdelim@false
 \fi\fi\fi\fi\fi\fi\fi\fi\fi\fi\fi\fi\fi\fi\fi}
\def\extra@{%
 \global\extra@false\rightdelim@\ifrightdelim@\global\extra@true            
 \else\ifx\next$\global\extra@true                                          
 \else\xdef\meaning@{\meaning\next..........}
 \expandafter\macro@\meaning@\macro@\ifmacro@                               
 \expandafter\DOTS@\meaning@\DOTS@
 \ifDOTS@
 \ifnum\DOTSCASE@=\tw@\global\extra@true                                    
 \fi\fi\fi\fi\fi}
\newif\ifbold@
\def\dotso@{\relaxnext@
 \ifbold@
  \let\next\delayed@
  \DNii@{\extra@\plainldots@\ifextra@\,\fi}%
 \else
  \DNii@{\DN@{\extra@\plainldots@\ifextra@\,\fi}\FN@\next@}%
 \fi
 \nextii@}
\def\extrap@#1{%
 \ifx\next,\DN@{#1\,}\else
 \ifx\next;\DN@{#1\,}\else
 \ifx\next.\DN@{#1\,}\else\extra@
 \ifextra@\DN@{#1\,}\else
 \let\next@#1\fi\fi\fi\fi\next@}
\def\ldots{\DN@{\extrap@\plainldots@}%
 \FN@\next@}
\def\cdots{\DN@{\extrap@\plaincdots@}%
 \FN@\next@}

\def\dotsc{\relaxnext@
 \DN@{\ifx\next;\plainldots@\,\else
  \ifx\next.\plainldots@\,\else\extra@\plainldots@
  \ifextra@\,\fi\fi\fi}%
 \FN@\next@}
\def\cdot{\mathchar"2201 }
\def\longrightarrow{\DOTSB\relbar\joinrel\rightarrow}

\def\mapsto{\DOTSB\mapstochar\rightarrow}

\def\dddot#1{{\mathop{#1}\limits^{\vbox to-1.4\ex@{\kern-\tw@\ex@
 \hbox{\rm...}\vss}}}}
\def\ddddot#1{{\mathop{#1}\limits^{\vbox to-1.4\ex@{\kern-\tw@\ex@
 \hbox{\rm....}\vss}}}}
\def\sphat{^{\mathchoice{}{}%
 {\,\,\botsmash{\hbox{\lower4\ex@\hbox{$\m@th\widehat{\null}$}}}}%
 {\,\botsmash{\hbox{\lower3\ex@\hbox{$\m@th\hat{\null}$}}}}}}

\def\spacute{^{\!\botsmash{\hbox{\lower\@ne ex\hbox{\'{}}}}}}
\def\spgrave{^{\mathchoice{}{}{}{\!}%
 \botsmash{\hbox{\lower\@ne ex\hbox{\`{}}}}}}
\def\spdot{^{\hbox{\raise\ex@\hbox{\rm.}}}}
\def\spddot{^{\hbox{\raise\ex@\hbox{\rm..}}}}
\def\spdddot{^{\hbox{\raise\ex@\hbox{\rm...}}}}
\def\spddddot{^{\hbox{\raise\ex@\hbox{\rm....}}}}
\def\spbreve{^{\!\botsmash{\hbox{\lower4\ex@\hbox{\u{}}}}}}

\def\textonlyfont@#1#2{\def#1{\RIfM@
 \Err@{Use \string#1\space only in text}\else#2\fi}}
\textonlyfont@\rm\tenrm
\textonlyfont@\it\tenit
\textonlyfont@\sl\tensl
\textonlyfont@\bf\tenbf
\def\oldnos#1{\RIfM@{\mathcode`\,="013B \fam\@ne#1}\else
 \leavevmode\hbox{$\m@th\mathcode`\,="013B \fam\@ne#1$}\fi}
\def\text{\RIfM@\expandafter\text@\else\expandafter\text@@\fi}
\def\text@@#1{\leavevmode\hbox{#1}}
\def\mathhexbox@#1#2#3{\text{$\m@th\mathchar"#1#2#3$}}
\def\dag{{\mathhexbox@279}}
\def\ddag{{\mathhexbox@27A}}
\def\S{{\mathhexbox@278}}
\def\P{{\mathhexbox@27B}}
\newif\iffirstchoice@
\firstchoice@true
\def\text@#1{\mathchoice
 {\hbox{\everymath{\displaystyle}\def\textfonti{\the\textfont\@ne}%
  \def\textfontii{\the\textfont\tw@}\textdef@@ T#1}}
 {\hbox{\firstchoice@false
  \everymath{\textstyle}\def\textfonti{\the\textfont\@ne}%
  \def\textfontii{\the\textfont\tw@}\textdef@@ T#1}}
 {\hbox{\firstchoice@false
  \everymath{\scriptstyle}\def\textfonti{\the\scriptfont\@ne}%
  \def\textfontii{\the\scriptfont\tw@}\textdef@@ S\rm#1}}
 {\hbox{\firstchoice@false
  \everymath{\scriptscriptstyle}\def\textfonti
  {\the\scriptscriptfont\@ne}%
  \def\textfontii{\the\scriptscriptfont\tw@}\textdef@@ s\rm#1}}}
\def\textdef@@#1{\textdef@#1\rm\textdef@#1\bf\textdef@#1\sl\textdef@#1\it}
\def\rmfam{0}
\def\textdef@#1#2{%
 \DN@{\csname\expandafter\eat@\string#2fam\endcsname}%
 \if S#1\edef#2{\the\scriptfont\next@\relax}%
 \else\if s#1\edef#2{\the\scriptscriptfont\next@\relax}%
 \else\edef#2{\the\textfont\next@\relax}\fi\fi}
\scriptfont\itfam\tenit \scriptscriptfont\itfam\tenit
\scriptfont\slfam\tensl \scriptscriptfont\slfam\tensl
\newif\iftopfolded@
\newif\ifbotfolded@
\def\topfoldedtext{\topfolded@true\botfolded@false\foldedtext@}
\def\botfoldedtext{\botfolded@true\topfolded@false\foldedtext@}
\def\foldedtext{\topfolded@false\botfolded@false\foldedtext@}
\Invalid@\foldedwidth
\def\foldedtext@{\relaxnext@
 \DN@{\ifx\next\foldedwidth\let\next@\nextii@\else
  \DN@{\nextii@\foldedwidth{.3\hsize}}\fi\next@}%
 \DNii@\foldedwidth##1##2{\setbox\z@\vbox
  {\normalbaselines\hsize##1\relax
  \tolerance1600 \noindent\ignorespaces##2}\ifbotfolded@\boxz@\else
  \iftopfolded@\vtop{\unvbox\z@}\else\vcenter{\boxz@}\fi\fi}%
 \FN@\next@}
\def\bold{\RIfM@\expandafter\bold@\else
 \expandafter\nonmatherr@\expandafter\bold\fi}
\def\bold@#1{{\bold@@{#1}}}
\def\bold@@#1{\fam\bffam\relax#1}
\def\slanted{\RIfM@\expandafter\slanted@\else
 \expandafter\nonmatherr@\expandafter\slanted\fi}
\def\slanted@#1{{\slanted@@{#1}}}
\def\slanted@@#1{\fam\slfam\relax#1}
\def\roman{\RIfM@\expandafter\roman@\else
 \expandafter\nonmatherr@\expandafter\roman\fi}
\def\roman@#1{{\roman@@{#1}}}
\def\roman@@#1{\fam\rmfam\relax#1}
\def\italic{\RIfM@\expandafter\italic@\else
 \expandafter\nonmatherr@\expandafter\italic\fi}
\def\italic@#1{{\italic@@{#1}}}
\def\italic@@#1{\fam\itfam\relax#1}
\def\Cal{\RIfM@\expandafter\Cal@\else
 \expandafter\nonmatherr@\expandafter\Cal\fi}
\def\Cal@#1{{\Cal@@{#1}}}
\def\Cal@@#1{\noaccents@\fam\tw@#1}
\mathchardef\Gamma="0000
\mathchardef\Delta="0001
\mathchardef\Theta="0002
\mathchardef\Lambda="0003
\mathchardef\Xi="0004
\mathchardef\Pi="0005
\mathchardef\Sigma="0006
\mathchardef\Upsilon="0007
\mathchardef\Phi="0008
\mathchardef\Psi="0009
\mathchardef\Omega="000A
\mathchardef\varGamma="0100
\mathchardef\varDelta="0101
\mathchardef\varTheta="0102
\mathchardef\varLambda="0103
\mathchardef\varXi="0104
\mathchardef\varPi="0105
\mathchardef\varSigma="0106
\mathchardef\varUpsilon="0107
\mathchardef\varPhi="0108
\mathchardef\varPsi="0109
\mathchardef\varOmega="010A
\newif\ifmsamloaded@
\newif\ifmsbmloaded@
\newif\ifeufmloaded@
\let\alloc@@\alloc@
\def\hexnumber@#1{\ifcase#1 0\or 1\or 2\or 3\or 4\or 5\or 6\or 7\or 8\or
 9\or A\or B\or C\or D\or E\or F\fi}
\edef\bffam@{\hexnumber@\bffam}
\def\loadmsam{\msamloaded@true
 \font@\tenmsa=msam10
 \font@\sevenmsa=msam7
 \font@\fivemsa=msam5
 \alloc@@8\fam\chardef\sixt@@n\msafam
 \textfont\msafam=\tenmsa
 \scriptfont\msafam=\sevenmsa
 \scriptscriptfont\msafam=\fivemsa
 \edef\msafam@{\hexnumber@\msafam}%
 \mathchardef\dabar@"0\msafam@39
 \def\dashrightarrow{\mathrel{\dabar@\dabar@\mathchar"0\msafam@4B}}%
 \def\dashleftarrow{\mathrel{\mathchar"0\msafam@4C\dabar@\dabar@}}%
 \let\dasharrow\dashrightarrow
 \def\ulcorner{\delimiter"4\msafam@70\msafam@70 }
 \def\urcorner{\delimiter"5\msafam@71\msafam@71 }
 \def\llcorner{\delimiter"4\msafam@78\msafam@78 }
 \def\lrcorner{\delimiter"5\msafam@79\msafam@79 }
 \def\yen{{\mathhexbox@\msafam@55 }}
 \def\checkmark{{\mathhexbox@\msafam@58 }}
 \def\circledR{{\mathhexbox@\msafam@72 }}
 \def\maltese{{\mathhexbox@\msafam@7A }}}
\def\loadmsbm{\msbmloaded@true
 \font@\tenmsb=msbm10
 \font@\sevenmsb=msbm7
 \font@\fivemsb=msbm5
 \alloc@@8\fam\chardef\sixt@@n\msbfam
 \textfont\msbfam=\tenmsb
 \scriptfont\msbfam=\sevenmsb
 \scriptscriptfont\msbfam=\fivemsb
 \edef\msbfam@{\hexnumber@\msbfam}%
 }
\def\widehat#1{\ifmsbmloaded@
  \setboxz@h{$\m@th#1$}\ifdim\wdz@>\tw@ em\mathaccent"0\msbfam@5B{#1}\else
  \mathaccent"0362{#1}\fi
 \else\mathaccent"0362{#1}\fi}
\def\widetilde#1{\ifmsbmloaded@
  \setboxz@h{$\m@th#1$}\ifdim\wdz@>\tw@ em\mathaccent"0\msbfam@5D{#1}\else
  \mathaccent"0365{#1}\fi
 \else\mathaccent"0365{#1}\fi}
\def\newsymbol#1#2#3#4#5{\define#1{}\let\next@\relax
 \ifnum#2=\@ne\ifmsamloaded@\let\next@\msafam@\fi\else
 \ifnum#2=\tw@\ifmsbmloaded@\let\next@\msbfam@\fi\fi\fi
 \ifx\next@\relax
  \ifnum#2>\tw@\Err@{\Invalid@@\string\newsymbol}\else
  \ifnum#2=\@ne\Err@{You must first \string\loadmsam}\else
   \Err@{You must first \string\loadmsbm}\fi\fi
 \else
  \mathchardef#1="#3\next@#4#5
 \fi}

\def\Bbb{\RIfM@\expandafter\Bbb@\else
 \expandafter\nonmatherr@\expandafter\Bbb\fi}
\def\Bbb@#1{{\Bbb@@{#1}}}
\def\Bbb@@#1{\noaccents@\fam\msbfam\relax#1}
\def\loadeufm{\eufmloaded@true
 \font@\teneufm=eufm10
 \font@\seveneufm=eufm7
 \font@\fiveeufm=eufm5
 \alloc@@8\fam\chardef\sixt@@n\eufmfam
 \textfont\eufmfam=\teneufm
 \scriptfont\eufmfam=\seveneufm
 \scriptscriptfont\eufmfam=\fiveeufm}
\def\frak{\RIfM@\expandafter\frak@\else
 \expandafter\nonmatherr@\expandafter\frak\fi}
\def\frak@#1{{\frak@@{#1}}}
\def\frak@@#1{\fam\eufmfam\relax#1}

\newif\ifcmmibloaded@
\newif\ifcmbsyloaded@
\def\loadbold{\cmmibloaded@true\cmbsyloaded@true
 \font@\tencmmib=cmmib10 \font@\sevencmmib=cmmib7 \font@\fivecmmib=cmmib5
 \skewchar\tencmmib='177 \skewchar\sevencmmib='177 \skewchar\fivecmmib='177
 \alloc@@8\fam\chardef\sixt@@n\cmmibfam
 \textfont\cmmibfam=\tencmmib
 \scriptfont\cmmibfam=\sevencmmib
 \scriptscriptfont\cmmibfam=\fivecmmib
 \edef\cmmibfam@{\hexnumber@\cmmibfam}%
 \font@\tencmbsy=cmbsy10 \font@\sevencmbsy=cmbsy7 \font@\fivecmbsy=cmbsy5
 \skewchar\tencmbsy='60 \skewchar\sevencmbsy='60 \skewchar\fivecmbsy='60
 \alloc@@8\fam\chardef\sixt@@n\cmbsyfam
 \textfont\cmbsyfam=\tencmbsy
 \scriptfont\cmbsyfam=\sevencmbsy
 \scriptscriptfont\cmbsyfam=\fivecmbsy
 \edef\cmbsyfam@{\hexnumber@\cmbsyfam}}
\def\mathchari@#1#2#3{\ifcmmibloaded@\mathchar"#1\cmmibfam@#2#3 \else
 \Err@{First bold symbol font not loaded}\fi}
\def\mathcharii@#1#2#3{\ifcmbsyloaded@\mathchar"#1\cmbsyfam@#2#3 \else
 \Err@{Second bold symbol font not loaded}\fi}
\def\boldkey#1{\ifcat\noexpand#1A%
  \ifcmmibloaded@{\fam\cmmibfam#1}\else
   \Err@{First bold symbol font not loaded}\fi
 \else
 \ifx#1!\mathchar"5\bffam@21 \else
 \ifx#1(\mathchar"4\bffam@28 \else\ifx#1)\mathchar"5\bffam@29 \else
 \ifx#1+\mathchar"2\bffam@2B \else\ifx#1:\mathchar"3\bffam@3A \else
 \ifx#1;\mathchar"6\bffam@3B \else\ifx#1=\mathchar"3\bffam@3D \else
 \ifx#1?\mathchar"5\bffam@3F \else\ifx#1[\mathchar"4\bffam@5B \else
 \ifx#1]\mathchar"5\bffam@5D \else
 \ifx#1,\mathchari@63B \else
 \ifx#1-\mathcharii@200 \else
 \ifx#1.\mathchari@03A \else
 \ifx#1/\mathchari@03D \else
 \ifx#1<\mathchari@33C \else
 \ifx#1>\mathchari@33E \else
 \ifx#1*\mathcharii@203 \else
 \ifx#1|\mathcharii@06A \else
 \ifx#10\bold0\else\ifx#11\bold1\else\ifx#12\bold2\else\ifx#13\bold3\else
 \ifx#14\bold4\else\ifx#15\bold5\else\ifx#16\bold6\else\ifx#17\bold7\else
 \ifx#18\bold8\else\ifx#19\bold9\else
  \Err@{\string\boldkey\space can't be used with #1}%
 \fi\fi\fi\fi\fi\fi\fi\fi\fi\fi\fi\fi\fi\fi\fi
 \fi\fi\fi\fi\fi\fi\fi\fi\fi\fi\fi\fi\fi\fi}
\def\boldsymbol#1{%
 \DN@{\Err@{You can't use \string\boldsymbol\space with \string#1}#1}%
 \ifcat\noexpand#1A%
   \let\next@\relax
  \ifcmmibloaded@{\fam\cmmibfam#1}\else\Err@{First bold symbol
   font not loaded}\fi
 \else
  \xdef\meaning@{\meaning#1.........}%
  \expandafter\math@\meaning@\math@
  \ifmath@
   \expandafter\mathch@\meaning@\mathch@
   \ifmathch@
    \expandafter\boldsymbol@@\meaning@\boldsymbol@@
   \fi
  \else
   \expandafter\macro@\meaning@\macro@
   \expandafter\delim@\meaning@\delim@
   \ifdelim@
    \expandafter\delim@@\meaning@\delim@@
   \else
    \boldsymbol@{#1}%
   \fi
  \fi
 \fi
 \next@}
\def\mathhexboxii@#1#2{\ifcmbsyloaded@\mathhexbox@{\cmbsyfam@}{#1}{#2}\else
  \Err@{Second bold symbol font not loaded}\fi}
\def\boldsymbol@#1{\let\next@\relax\let\next#1%
 \ifx\next\cdot\mathcharii@201 \else
 \ifx\next\prime{{\null\mathcharii@030 \null}}\else
 \ifx\next\lbrack\mathchar"4\bffam@5B \else
 \ifx\next\rbrack\mathchar"5\bffam@5D \else
 \ifx\next\{\mathcharii@466 \else
 \ifx\next\lbrace\mathcharii@466 \else
 \ifx\next\}\mathcharii@567 \else
 \ifx\next\rbrace\mathcharii@567 \else
 \ifx\next\surd{{\mathcharii@170}}\else
 \ifx\next\S{{\mathhexboxii@78}}\else
 \ifx\next\P{{\mathhexboxii@7B}}\else
 \ifx\next\dag{{\mathhexboxii@79}}\else
 \ifx\next\ddag{{\mathhexboxii@7A}}\else
 \DN@{\Err@{You can't use \string\boldsymbol\space with \string#1}#1}%
 \fi\fi\fi\fi\fi\fi\fi\fi\fi\fi\fi\fi\fi}
\def\boldsymbol@@#1.#2\boldsymbol@@{\classnum@#1 \count@@@\classnum@        
 \divide\classnum@4096 \count@\classnum@                                    
 \multiply\count@4096 \advance\count@@@-\count@ \count@@\count@@@           
 \divide\count@@@\@cclvi \count@\count@@                                    
 \multiply\count@@@\@cclvi \advance\count@@-\count@@@                       
 \divide\count@@@\@cclvi                                                    
 \multiply\classnum@4096 \advance\classnum@\count@@                         
 \ifnum\count@@@=\z@                                                        
  \count@"\bffam@ \multiply\count@\@cclvi
  \advance\classnum@\count@
  \DN@{\mathchar\number\classnum@}%
 \else
  \ifnum\count@@@=\@ne                                                      
   \ifcmmibloaded@
   \count@"\cmmibfam@ \multiply\count@\@cclvi
   \advance\classnum@\count@
   \DN@{\mathchar\number\classnum@}%
   \else\DN@{\Err@{First bold symbol font not loaded}}\fi
  \else
   \ifnum\count@@@=\tw@                                                    
  \ifcmbsyloaded@
    \count@"\cmbsyfam@ \multiply\count@\@cclvi
    \advance\classnum@\count@
    \DN@{\mathchar\number\classnum@}%
  \else\DN@{\Err@{Second bold symbol font not loaded}}\fi
  \fi
 \fi
\fi}
\newif\ifdelim@
\newcount\delimcount@
{\uccode`6=`\\ \uccode`7=`d \uccode`8=`e \uccode`9=`l
 \uppercase{\gdef\delim@#1#2#3#4#5\delim@
  {\delim@false\ifx 6#1\ifx 7#2\ifx 8#3\ifx 9#4\delim@true
   \xdef\meaning@{#5}\fi\fi\fi\fi}}}
\def\delim@@#1"#2#3#4#5#6\delim@@{\if#32%
\let\next@\relax
 \ifcmbsyloaded@
 \mathcharii@#2#4#5 \else\Err@{Second bold family not loaded}\fi\fi}
\def\vert{\delimiter"026A30C }
\def\Vert{\delimiter"026B30D }
\let\|\Vert

\def\boldkeydots@#1{\bold@true\let\next=#1\let\delayed@=#1\mdots@@
 \boldkey#1\bold@false}  
\def\boldsymboldots@#1{\bold@true\let\next#1\let\delayed@#1\mdots@@
 \boldsymbol#1\bold@false}
\newif\ifeufbloaded@
\def\loadeufb{\eufbloaded@true
 \font@\teneufb=eufb10
 \font@\seveneufb=eufb7
 \font@\fiveeufb=eufb5
 \alloc@@8\fam\chardef\sixt@@n\eufbfam
 \textfont\eufbfam=\teneufb
 \scriptfont\eufbfam=\seveneufb
 \scriptscriptfont\eufbfam=\fiveeufb
 \edef\eufbfam@{\hexnumber@\eufbfam}}
\newif\ifeusmloaded@
\def\loadeusm{\eusmloaded@true
 \font@\teneusm=eusm10
 \font@\seveneusm=eusm7
 \font@\fiveeusm=eusm5
 \alloc@@8\fam\chardef\sixt@@n\eusmfam
 \textfont\eusmfam=\teneusm
 \scriptfont\eusmfam=\seveneusm
 \scriptscriptfont\eusmfam=\fiveeusm
 \edef\eusmfam@{\hexnumber@\eusmfam}}
\newif\ifeusbloaded@
\def\loadeusb{\eusbloaded@true
 \font@\teneusb=eusb10
 \font@\seveneusb=eusb7
 \font@\fiveeusb=eusb5
 \alloc@@8\fam\chardef\sixt@@n\eusbfam
 \textfont\eusbfam=\teneusb
 \scriptfont\eusbfam=\seveneusb
 \scriptscriptfont\eusbfam=\fiveeusb
 \edef\eusbfam@{\hexnumber@\eusbfam}}
\newif\ifeurmloaded@
\def\loadeurm{\eurmloaded@true
 \font@\teneurm=eurm10
 \font@\seveneurm=eurm7
 \font@\fiveeurm=eurm5
 \alloc@@8\fam\chardef\sixt@@n\eurmfam
 \textfont\eurmfam=\teneurm
 \scriptfont\eurmfam=\seveneurm
 \scriptscriptfont\eurmfam=\fiveeurm
 \edef\eurmfam@{\hexnumber@\eurmfam}}
\newif\ifeurbloaded@
\def\loadeurb{\eurbloaded@true
 \font@\teneurb=eurb10
 \font@\seveneurb=eurb7
 \font@\fiveeurb=eurb5
 \alloc@@8\fam\chardef\sixt@@n\eurbfam
 \textfont\eurbfam=\teneurb
 \scriptfont\eurbfam=\seveneurb
 \scriptscriptfont\eurbfam=\fiveeurb
 \edef\eurbfam@{\hexnumber@\eurbfam}}
\def\accentclass@{7}
\def\noaccents@{\def\accentclass@{0}}
\def\makeacc@#1#2{\def#1{\mathaccent"\accentclass@#2 }}
\makeacc@\hat{05E}
\makeacc@\check{014}
\makeacc@\tilde{07E}
\makeacc@\acute{013}
\makeacc@\grave{012}
\makeacc@\dot{05F}
\makeacc@\ddot{07F}
\makeacc@\breve{015}
\makeacc@\bar{016}

\newcount\skewcharcount@
\newcount\familycount@
\def\theskewchar@{\familycount@\@ne
 \global\skewcharcount@\the\skewchar\textfont\@ne                           
 \ifnum\fam>\m@ne\ifnum\fam<16
  \global\familycount@\the\fam\relax
  \global\skewcharcount@\the\skewchar\textfont\the\fam\relax\fi\fi          
 \ifnum\skewcharcount@>\m@ne
  \ifnum\skewcharcount@<128
  \multiply\familycount@256
  \global\advance\skewcharcount@\familycount@
  \global\advance\skewcharcount@28672
  \mathchar\skewcharcount@\else
  \global\skewcharcount@\m@ne\fi\else
 \global\skewcharcount@\m@ne\fi}                                            
\newcount\pointcount@
\def\getpoints@#1.#2\getpoints@{\pointcount@#1 }
\newdimen\accentdimen@
\newcount\accentmu@
\def\dimentomu@{\multiply\accentdimen@ 100
 \expandafter\getpoints@\the\accentdimen@\getpoints@
 \multiply\pointcount@18
 \divide\pointcount@\@m
 \global\accentmu@\pointcount@}
\def\Makeacc@#1#2{\def#1{\RIfM@\DN@{\mathaccent@
 {"\accentclass@#2 }}\else\DN@{\nonmatherr@{#1}}\fi\next@}}
\def\unbracefonts@{\let\Cal@\Cal@@\let\roman@\roman@@\let\bold@\bold@@
 \let\slanted@\slanted@@}
\def\mathaccent@#1#2{\ifnum\fam=\m@ne\xdef\thefam@{1}\else
 \xdef\thefam@{\the\fam}\fi                                                 
 \accentdimen@\z@                                                           
 \setboxz@h{\unbracefonts@$\m@th\fam\thefam@\relax#2$}
 \ifdim\accentdimen@=\z@\DN@{\mathaccent#1{#2}}
  \setbox@ne\hbox{\unbracefonts@$\m@th\fam\thefam@\relax#2\theskewchar@$}
  \setbox\tw@\hbox{$\m@th\ifnum\skewcharcount@=\m@ne\else
   \mathchar\skewcharcount@\fi$}
  \global\accentdimen@\wd@ne\global\advance\accentdimen@-\wdz@
  \global\advance\accentdimen@-\wd\tw@                                     
  \global\multiply\accentdimen@\tw@
  \dimentomu@\global\advance\accentmu@\@ne                                 
 \else\DN@{{\mathaccent#1{#2\mkern\accentmu@ mu}%
    \mkern-\accentmu@ mu}{}}\fi                                             
 \next@}\Makeacc@\Hat{05E}
\Makeacc@\Check{014}
\Makeacc@\Tilde{07E}
\Makeacc@\Acute{013}
\Makeacc@\Grave{012}
\Makeacc@\Dot{05F}
\Makeacc@\Ddot{07F}
\Makeacc@\Breve{015}
\Makeacc@\Bar{016}
\def\Vec{\RIfM@\DN@{\mathaccent@{"017E }}\else
 \DN@{\nonmatherr@\Vec}\fi\next@}
\def\newbox@{\alloc@4\box\chardef\insc@unt}
\def\accentedsymbol#1#2{\expandafter\newbox@\csname\expandafter
  \eat@\string#1@box\endcsname
 \expandafter\setbox\csname\expandafter\eat@
  \string#1@box\endcsname\hbox{$\m@th#2$}\define
  #1{\expandafter\copy\csname\expandafter\eat@\string#1@box\endcsname{}}}
\def\sqrt#1{\radical"270370 {#1}}
\let\underline@\underline
\let\overline@\overline
\def\underline#1{\underline@{#1}}
\def\overline#1{\overline@{#1}}
\Invalid@\leftroot
\Invalid@\uproot
\newcount\uproot@
\newcount\leftroot@
\def\root{\relaxnext@
  \DN@{\ifx\next\uproot\let\next@\nextii@\else
   \ifx\next\leftroot\let\next@\nextiii@\else
   \let\next@\plainroot@\fi\fi\next@}%
  \DNii@\uproot##1{\uproot@##1\relax\FN@\nextiv@}%
  \def\nextiv@{\ifx\next\space@\DN@. {\FN@\nextv@}\else
   \DN@.{\FN@\nextv@}\fi\next@.}%
  \def\nextv@{\ifx\next\leftroot\let\next@\nextvi@\else
   \let\next@\plainroot@\fi\next@}%
  \def\nextvi@\leftroot##1{\leftroot@##1\relax\plainroot@}%
   \def\nextiii@\leftroot##1{\leftroot@##1\relax\FN@\nextvii@}%
  \def\nextvii@{\ifx\next\space@
   \DN@. {\FN@\nextviii@}\else
   \DN@.{\FN@\nextviii@}\fi\next@.}%
  \def\nextviii@{\ifx\next\uproot\let\next@\nextix@\else
   \let\next@\plainroot@\fi\next@}%
  \def\nextix@\uproot##1{\uproot@##1\relax\plainroot@}%
  \bgroup\uproot@\z@\leftroot@\z@\FN@\next@}
\def\plainroot@#1\of#2{\setbox\rootbox\hbox{$\m@th\scriptscriptstyle{#1}$}%
 \mathchoice{\r@@t\displaystyle{#2}}{\r@@t\textstyle{#2}}
 {\r@@t\scriptstyle{#2}}{\r@@t\scriptscriptstyle{#2}}\egroup}
\def\r@@t#1#2{\setboxz@h{$\m@th#1\sqrt{#2}$}%
 \dimen@\ht\z@\advance\dimen@-\dp\z@
 \setbox@ne\hbox{$\m@th#1\mskip\uproot@ mu$}\advance\dimen@ by1.667\wd@ne
 \mkern-\leftroot@ mu\mkern5mu\raise.6\dimen@\copy\rootbox
 \mkern-10mu\mkern\leftroot@ mu\boxz@}
\def\boxed#1{\setboxz@h{$\m@th\displaystyle{#1}$}\dimen@.4\ex@
 \advance\dimen@3\ex@\advance\dimen@\dp\z@
 \hbox{\lower\dimen@\hbox{%
 \vbox{\hrule height.4\ex@
 \hbox{\vrule width.4\ex@\hskip3\ex@\vbox{\vskip3\ex@\boxz@\vskip3\ex@}%
 \hskip3\ex@\vrule width.4\ex@}\hrule height.4\ex@}%
 }}}
\let\ampersand@\relax
\newdimen\minaw@
\minaw@11.11128\ex@
\newdimen\minCDaw@
\minCDaw@2.5pc
\def\minCDarrowwidth#1{\RIfMIfI@\onlydmatherr@\minCDarrowwidth
 \else\minCDaw@#1\relax\fi\else\onlydmatherr@\minCDarrowwidth\fi}
\newif\ifCD@
\def\CD{\bgroup\vspace@\relax\let\ampersand@&\iffalse}\fi
 \CD@true\vcenter\bgroup\Let@\tabskip\z@skip\baselineskip20\ex@
 \lineskip3\ex@\lineskiplimit3\ex@\halign\bgroup
 &\hfill$\m@th##$\hfill\crcr}
\def\endCD{\crcr\egroup\egroup\egroup}
\newdimen\bigaw@
\atdef@>#1>#2>{\ampersand@                                                  
 \setboxz@h{$\m@th\ssize\;{#1}\;\;$}
 \setbox@ne\hbox{$\m@th\ssize\;{#2}\;\;$}
 \setbox\tw@\hbox{$\m@th#2$}
 \ifCD@\global\bigaw@\minCDaw@\else\global\bigaw@\minaw@\fi                 
 \ifdim\wdz@>\bigaw@\global\bigaw@\wdz@\fi
 \ifdim\wd@ne>\bigaw@\global\bigaw@\wd@ne\fi                                
 \ifCD@\hskip.5em\fi                                                        
 \ifdim\wd\tw@>\z@
  \mathrel{\mathop{\hbox to\bigaw@{\rightarrowfill}}\limits^{#1}_{#2}}
 \else\mathrel{\mathop{\hbox to\bigaw@{\rightarrowfill}}\limits^{#1}}\fi    
 \ifCD@\hskip.5em\fi                                                       
 \ampersand@}                                                              
\atdef@<#1<#2<{\ampersand@\setboxz@h{$\m@th\ssize\;\;{#1}\;$}%
 \setbox@ne\hbox{$\m@th\ssize\;\;{#2}\;$}\setbox\tw@\hbox{$\m@th#2$}%
 \ifCD@\global\bigaw@\minCDaw@\else\global\bigaw@\minaw@\fi
 \ifdim\wdz@>\bigaw@\global\bigaw@\wdz@\fi
 \ifdim\wd@ne>\bigaw@\global\bigaw@\wd@ne\fi
 \ifCD@\hskip.5em\fi
 \ifdim\wd\tw@>\z@
  \mathrel{\mathop{\hbox to\bigaw@{\leftarrowfill}}\limits^{#1}_{#2}}\else
  \mathrel{\mathop{\hbox to\bigaw@{\leftarrowfill}}\limits^{#1}}\fi
 \ifCD@\hskip.5em\fi\ampersand@}
\atdef@)#1)#2){\ampersand@
 \setboxz@h{$\m@th\ssize\;{#1}\;\;$}%
 \setbox@ne\hbox{$\m@th\ssize\;{#2}\;\;$}%
 \setbox\tw@\hbox{$\m@th#2$}%
 \ifCD@
 \global\bigaw@\minCDaw@\else\global\bigaw@\minaw@\fi
 \ifdim\wdz@>\bigaw@\global\bigaw@\wdz@\fi
 \ifdim\wd@ne>\bigaw@\global\bigaw@\wd@ne\fi
 \ifCD@\hskip.5em\fi
 \ifdim\wd\tw@>\z@
  \mathrel{\mathop{\hbox to\bigaw@{\rightarrowfill}}\limits^{#1}_{#2}}%
 \else\mathrel{\mathop{\hbox to\bigaw@{\rightarrowfill}}\limits^{#1}}\fi
 \ifCD@\hskip.5em\fi
 \ampersand@}
\atdef@(#1(#2({\ampersand@\setboxz@h{$\m@th\ssize\;\;{#1}\;$}%
 \setbox@ne\hbox{$\m@th\ssize\;\;{#2}\;$}\setbox\tw@\hbox{$\m@th#2$}%
 \ifCD@\global\bigaw@\minCDaw@\else\global\bigaw@\minaw@\fi
 \ifdim\wdz@>\bigaw@\global\bigaw@\wdz@\fi
 \ifdim\wd@ne>\bigaw@\global\bigaw@\wd@ne\fi
 \ifCD@\hskip.5em\fi
 \ifdim\wd\tw@>\z@
  \mathrel{\mathop{\hbox to\bigaw@{\leftarrowfill}}\limits^{#1}_{#2}}\else
  \mathrel{\mathop{\hbox to\bigaw@{\leftarrowfill}}\limits^{#1}}\fi
 \ifCD@\hskip.5em\fi\ampersand@}
\atdef@ A#1A#2A{\llap{$\m@th\vcenter{\hbox
 {$\ssize#1$}}$}\Big\uparrow\rlap{$\m@th\vcenter{\hbox{$\ssize#2$}}$}&&}
\atdef@ V#1V#2V{\llap{$\m@th\vcenter{\hbox
 {$\ssize#1$}}$}\Big\downarrow\rlap{$\m@th\vcenter{\hbox{$\ssize#2$}}$}&&}
\atdef@={&\hskip.5em\mathrel
 {\vbox{\hrule width\minCDaw@\vskip3\ex@\hrule width
 \minCDaw@}}\hskip.5em&}
\atdef@|{\Big\Vert&&}
\atdef@@\vert{\Big\Vert&&}
\def\pretend#1\haswidth#2{\setboxz@h{$\m@th\scriptstyle{#2}$}\hbox
 to\wdz@{\hfill$\m@th\scriptstyle{#1}$\hfill}}
\def\pmb{\RIfM@\expandafter\mathpalette\expandafter\pmb@\else
 \expandafter\pmb@@\fi}
\def\pmb@@#1{\leavevmode\setboxz@h{#1}\kern-.025em\copy\z@\kern-\wdz@
 \kern-.05em\copy\z@\kern-\wdz@\kern-.025em\raise.0433em\boxz@}
\def\binrel@@#1{\ifdim\wd2<\z@\mathbin{#1}\else\ifdim\wd\tw@>\z@
 \mathrel{#1}\else{#1}\fi\fi}
\newdimen\pmbraise@
\def\pmb@#1#2{\setbox\thr@@\hbox{$\m@th#1{#2}$}%
 \setbox4 \hbox{$\m@th#1\mkern.7794mu$}\pmbraise@\wd4
 \binrel@{#2}\binrel@@{\mkern-.45mu\copy\thr@@\kern-\wd\thr@@
 \mkern-.9mu\copy\thr@@\kern-\wd\thr@@\mkern-.45mu\raise\pmbraise@\box\thr@@}}
\def\documentstyle#1{\input #1.sty\relax}
\font\dummyft@=dummy
\fontdimen1 \dummyft@=\z@
\fontdimen2 \dummyft@=\z@
\fontdimen3 \dummyft@=\z@
\fontdimen4 \dummyft@=\z@
\fontdimen5 \dummyft@=\z@
\fontdimen6 \dummyft@=\z@
\fontdimen7 \dummyft@=\z@
\fontdimen8 \dummyft@=\z@
\fontdimen9 \dummyft@=\z@
\fontdimen10 \dummyft@=\z@
\fontdimen11 \dummyft@=\z@
\fontdimen12 \dummyft@=\z@
\fontdimen13 \dummyft@=\z@
\fontdimen14 \dummyft@=\z@
\fontdimen15 \dummyft@=\z@
\fontdimen16 \dummyft@=\z@
\fontdimen17 \dummyft@=\z@
\fontdimen18 \dummyft@=\z@
\fontdimen19 \dummyft@=\z@
\fontdimen20 \dummyft@=\z@
\fontdimen21 \dummyft@=\z@
\fontdimen22 \dummyft@=\z@
\def\fontlist@{\\{\tenrm}\\{\sevenrm}\\{\fiverm}\\{\teni}\\{\seveni}%
 \\{\fivei}\\{\tensy}\\{\sevensy}\\{\fivesy}\\{\tenex}\\{\tenbf}\\{\sevenbf}%
 \\{\fivebf}\\{\tensl}\\{\tenit}}
\def\font@#1=#2 {\rightappend@#1\to\fontlist@\font#1=#2 }
\def\dodummy@{{\def\\##1{\global\let##1\dummyft@}\fontlist@}}
\def\nopages@{\output={\setbox\z@\box255 \deadcycles\z@}%
 \alloc@5\toks\toksdef\@cclvi\output}
\let\galleys\nopages@
\newif\ifsyntax@
\newcount\countxviii@
\def\syntax{\syntax@true\dodummy@\countxviii@\count18
 \loop\ifnum\countxviii@>\m@ne\textfont\countxviii@=\dummyft@
 \scriptfont\countxviii@=\dummyft@\scriptscriptfont\countxviii@=\dummyft@
 \advance\countxviii@\m@ne\repeat                                           
 \dummyft@\tracinglostchars\z@\nopages@\frenchspacing\hbadness\@M}
\def\S@{S } \def\G@{G } \def\P@{P }
\newif\ifbadans@
\def\printoptions{\W@{Do you want S(yntax check),
  G(alleys) or P(ages)?^^JType S, G or P, follow by <return>: }\loop
 \read\m@ne to\ans@
 \xdef\next@{\def\noexpand\Ans@{\ans@}}\uppercase\expandafter{\next@}
 \ifx\Ans@\S@\badans@false\syntax\else
 \ifx\Ans@\G@\badans@false\galleys\else
 \ifx\Ans@\P@\badans@false\else
 \badans@true\fi\fi\fi
 \ifbadans@\W@{Type S, G or P, follow by <return>: }%
 \repeat}
\def\alloc@#1#2#3#4#5{\global\advance\count1#1by\@ne
 \ch@ck#1#4#2\allocationnumber=\count1#1
 \global#3#5=\allocationnumber
 \ifalloc@\wlog{\string#5=\string#2\the\allocationnumber}\fi}
\def\document{\def\alloclist@{}\def\fontlist@{}}
\let\enddocument\bye

\let\proclaim\undefined
\let\footnote\undefined
\let\=\undefined
\let\>\undefined

\catcode`\@=\active

%
%
\def\next{AMSPPT}\ifx\styname\next \endinput\fi
\catcode`\@=11
\def\styname{AMSPPT}
\def\styversion{2.0}
{\W@{\styname.STY - Version \styversion}\W@{}}
\hyphenation{acad-e-my acad-e-mies af-ter-thought anom-aly anom-alies
an-ti-deriv-a-tive an-tin-o-my an-tin-o-mies apoth-e-o-ses apoth-e-o-sis
ap-pen-dix ar-che-typ-al as-sign-a-ble as-sist-ant-ship as-ymp-tot-ic
asyn-chro-nous at-trib-uted at-trib-ut-able bank-rupt bank-rupt-cy
bi-dif-fer-en-tial blue-print busier busiest cat-a-stroph-ic
cat-a-stroph-i-cally con-gress cross-hatched data-base de-fin-i-tive
de-riv-a-tive dis-trib-ute dri-ver dri-vers eco-nom-ics econ-o-mist
elit-ist equi-vari-ant ex-quis-ite ex-tra-or-di-nary flow-chart
for-mi-da-ble forth-right friv-o-lous ge-o-des-ic ge-o-det-ic geo-met-ric
griev-ance griev-ous griev-ous-ly hexa-dec-i-mal ho-lo-no-my ho-mo-thetic
ideals idio-syn-crasy in-fin-ite-ly in-fin-i-tes-i-mal ir-rev-o-ca-ble
key-stroke lam-en-ta-ble light-weight mal-a-prop-ism man-u-script
mar-gin-al meta-bol-ic me-tab-o-lism meta-lan-guage me-trop-o-lis
met-ro-pol-i-tan mi-nut-est mol-e-cule mono-chrome mono-pole mo-nop-oly
mono-spline mo-not-o-nous mul-ti-fac-eted mul-ti-plic-able non-euclid-ean
non-iso-mor-phic non-smooth par-a-digm par-a-bol-ic pa-rab-o-loid
pa-ram-e-trize para-mount pen-ta-gon phe-nom-e-non post-script pre-am-ble
pro-ce-dur-al pro-hib-i-tive pro-hib-i-tive-ly pseu-do-dif-fer-en-tial
pseu-do-fi-nite pseu-do-nym qua-drat-ics quad-ra-ture qua-si-smooth
qua-si-sta-tion-ary qua-si-tri-an-gu-lar quin-tes-sence quin-tes-sen-tial
re-arrange-ment rec-tan-gle ret-ri-bu-tion retro-fit retro-fit-ted
right-eous right-eous-ness ro-bot ro-bot-ics sched-ul-ing se-mes-ter
semi-def-i-nite semi-ho-mo-thet-ic set-up se-vere-ly side-step sov-er-eign
spe-cious spher-oid spher-oid-al star-tling star-tling-ly
sta-tis-tics sto-chas-tic straight-est strange-ness strat-a-gem strong-hold
sum-ma-ble symp-to-matic syn-chro-nous topo-graph-i-cal tra-vers-a-ble
tra-ver-sal tra-ver-sals treach-ery turn-around un-at-tached un-err-ing-ly
white-space wide-spread wing-spread wretch-ed wretch-ed-ly Brown-ian
Eng-lish Euler-ian Feb-ru-ary Gauss-ian Grothen-dieck Hamil-ton-ian
Her-mit-ian Jan-u-ary Japan-ese Kor-te-weg Le-gendre Lip-schitz
Lip-schitz-ian Mar-kov-ian Noe-ther-ian No-vem-ber Rie-mann-ian
Schwarz-schild Sep-tem-ber}
\Invalid@\nofrills
\Invalid@\usualspace
\newif\ifnofrills@
\def\nofrills@#1#2{\relaxnext@
  \DN@{\ifx\next\nofrills
    \nofrills@true\let#2\relax\DN@\nofrills{\nextii@}%
  \else
    \nofrills@false\def#2{#1}\let\next@\nextii@\fi
\next@}}
\def\usualspace@#1{\ifnofrills@\def\usualspace{#1}\fi}
\def\addto#1#2{\csname \expandafter\eat@\string#1@\endcsname
  \expandafter{\the\csname \expandafter\eat@\string#1@\endcsname#2}}
\newdimen\bigsize@
\def\big@#1#2{{\hbox{$\left#2\vcenter to#1\bigsize@{}%
  \right.\nulldelimiterspace\z@\m@th$}}}
\def\big{\big@\@ne}
\def\Big{\big@{1.5}}
\def\bigg{\big@\tw@}
\def\Bigg{\big@{2.5}}
\def\raggedcenter@{\leftskip\z@ plus.4\hsize \rightskip\leftskip
 \parfillskip\z@ \parindent\z@ \spaceskip.3333em \xspaceskip.5em
 \pretolerance9999\tolerance9999 \exhyphenpenalty\@M
 \hyphenpenalty\@M \let\\\linebreak}
\def\upperspecialchars{\def\ss{SS}\let\i=I\let\j=J\let\ae\AE\let\oe\OE
  \let\o\O\let\aa\AA\let\l\L}
\def\uppercasetext@#1{%
  {\spaceskip1.2\fontdimen2\the\font plus1.2\fontdimen3\the\font
   \upperspecialchars\uctext@#1$\m@th\aftergroup\eat@$}}
\def\uctext@#1$#2${\endash@#1-\endash@$#2$\uctext@}
\def\endash@#1-#2\endash@{\uppercase{#1}\if\notempty{#2}--\endash@#2\endash@\fi}
\def\runaway@#1{\DN@{#1}\ifx\envir@\next@
  \Err@{You seem to have a missing or misspelled \string\end#1 ...}%
  \let\envir@\empty\fi}
\newif\iftemp@
\def\notempty#1{TT\fi\def\test@{#1}\ifx\test@\empty\temp@false
  \else\temp@true\fi \iftemp@}
\font@\tensmc=cmcsc10
\font@\sevenex=cmex7
\font@\sevenit=cmti7
\font@\eightrm=cmr8 
\font@\sixrm=cmr6 
\font@\eighti=cmmi8     \skewchar\eighti='177 
\font@\sixi=cmmi6       \skewchar\sixi='177   
\font@\eightsy=cmsy8    \skewchar\eightsy='60 
\font@\sixsy=cmsy6      \skewchar\sixsy='60   
\font@\eightex=cmex8
\font@\eightbf=cmbx8 
\font@\sixbf=cmbx6   
\font@\eightit=cmti8 
\font@\eightsl=cmsl8 
\font@\eightsmc=cmcsc8
\font@\eighttt=cmtt8 
\loadmsam
\loadmsbm
\loadeufm
\input amssym.tex\relax
\newtoks\tenpoint@
\def\tenpoint{\normalbaselineskip12\p@
 \abovedisplayskip12\p@ plus3\p@ minus9\p@
 \belowdisplayskip\abovedisplayskip
 \abovedisplayshortskip\z@ plus3\p@
 \belowdisplayshortskip7\p@ plus3\p@ minus4\p@
 \textonlyfont@\rm\tenrm \textonlyfont@\it\tenit
 \textonlyfont@\sl\tensl \textonlyfont@\bf\tenbf
 \textonlyfont@\smc\tensmc \textonlyfont@\tt\tentt
 \ifsyntax@ \def\big##1{{\hbox{$\left##1\right.$}}}%
  \let\Big\big \let\bigg\big \let\Bigg\big
 \else
  \textfont\z@=\tenrm  \scriptfont\z@=\sevenrm  \scriptscriptfont\z@=\fiverm
  \textfont\@ne=\teni  \scriptfont\@ne=\seveni  \scriptscriptfont\@ne=\fivei
  \textfont\tw@=\tensy \scriptfont\tw@=\sevensy \scriptscriptfont\tw@=\fivesy
  \textfont\thr@@=\tenex \scriptfont\thr@@=\sevenex
        \scriptscriptfont\thr@@=\sevenex
  \textfont\itfam=\tenit \scriptfont\itfam=\sevenit
        \scriptscriptfont\itfam=\sevenit
  \textfont\bffam=\tenbf \scriptfont\bffam=\sevenbf
        \scriptscriptfont\bffam=\fivebf
  \setbox\strutbox\hbox{\vrule height8.5\p@ depth3.5\p@ width\z@}%
  \setbox\strutbox@\hbox{\lower.5\normallineskiplimit\vbox{%
        \kern-\normallineskiplimit\copy\strutbox}}%
 \setbox\z@\vbox{\hbox{$($}\kern\z@}\bigsize@=1.2\ht\z@
 \fi
 \normalbaselines\rm\ex@.2326ex\jot3\ex@\the\tenpoint@}
\newtoks\eightpoint@
\def\eightpoint{\normalbaselineskip10\p@
 \abovedisplayskip10\p@ plus2.4\p@ minus7.2\p@
 \belowdisplayskip\abovedisplayskip
 \abovedisplayshortskip\z@ plus2.4\p@
 \belowdisplayshortskip5.6\p@ plus2.4\p@ minus3.2\p@
 \textonlyfont@\rm\eightrm \textonlyfont@\it\eightit
 \textonlyfont@\sl\eightsl \textonlyfont@\bf\eightbf
 \textonlyfont@\smc\eightsmc \textonlyfont@\tt\eighttt
 \ifsyntax@\def\big##1{{\hbox{$\left##1\right.$}}}%
  \let\Big\big \let\bigg\big \let\Bigg\big
 \else
  \textfont\z@=\eightrm \scriptfont\z@=\sixrm \scriptscriptfont\z@=\fiverm
  \textfont\@ne=\eighti \scriptfont\@ne=\sixi \scriptscriptfont\@ne=\fivei
  \textfont\tw@=\eightsy \scriptfont\tw@=\sixsy \scriptscriptfont\tw@=\fivesy
  \textfont\thr@@=\eightex \scriptfont\thr@@=\sevenex
   \scriptscriptfont\thr@@=\sevenex
  \textfont\itfam=\eightit \scriptfont\itfam=\sevenit
   \scriptscriptfont\itfam=\sevenit
  \textfont\bffam=\eightbf \scriptfont\bffam=\sixbf
   \scriptscriptfont\bffam=\fivebf
 \setbox\strutbox\hbox{\vrule height7\p@ depth3\p@ width\z@}%
 \setbox\strutbox@\hbox{\raise.5\normallineskiplimit\vbox{%
   \kern-\normallineskiplimit\copy\strutbox}}%
 \setbox\z@\vbox{\hbox{$($}\kern\z@}\bigsize@=1.2\ht\z@
 \fi
 \normalbaselines\eightrm\ex@.2326ex\jot3\ex@\the\eightpoint@}
\parindent1pc
\normallineskiplimit\p@
\newdimen\indenti \indenti=2pc
\def\pageheight#1{\vsize#1}
\def\pagewidth#1{\hsize#1%
   \captionwidth@\hsize \advance\captionwidth@-2\indenti}
\pagewidth{30pc} \pageheight{47pc}
\def\topmatter{%
 \ifx\undefined\msafam
 \else\font@\eightmsa=msam8 \font@\sixmsa=msam6
   \ifsyntax@\else \addto\tenpoint{\textfont\msafam=\tenmsa
              \scriptfont\msafam=\sevenmsa \scriptscriptfont\msafam=\fivemsa}%
     \addto\eightpoint{\textfont\msafam=\eightmsa \scriptfont\msafam=\sixmsa
              \scriptscriptfont\msafam=\fivemsa}%
   \fi
 \fi
 \ifx\undefined\msbfam
 \else\font@\eightmsb=msbm8 \font@\sixmsb=msbm6
   \ifsyntax@\else \addto\tenpoint{\textfont\msbfam=\tenmsb
         \scriptfont\msbfam=\sevenmsb \scriptscriptfont\msbfam=\fivemsb}%
     \addto\eightpoint{\textfont\msbfam=\eightmsb \scriptfont\msbfam=\sixmsb
         \scriptscriptfont\msbfam=\fivemsb}%
   \fi
 \fi
 \ifx\undefined\eufmfam
 \else \font@\eighteufm=eufm8 \font@\sixeufm=eufm6
   \ifsyntax@\else \addto\tenpoint{\textfont\eufmfam=\teneufm
       \scriptfont\eufmfam=\seveneufm \scriptscriptfont\eufmfam=\fiveeufm}%
     \addto\eightpoint{\textfont\eufmfam=\eighteufm
       \scriptfont\eufmfam=\sixeufm \scriptscriptfont\eufmfam=\fiveeufm}%
   \fi
 \fi
 \ifx\undefined\eufbfam
 \else \font@\eighteufb=eufb8 \font@\sixeufb=eufb6
   \ifsyntax@\else \addto\tenpoint{\textfont\eufbfam=\teneufb
      \scriptfont\eufbfam=\seveneufb \scriptscriptfont\eufbfam=\fiveeufb}%
    \addto\eightpoint{\textfont\eufbfam=\eighteufb
      \scriptfont\eufbfam=\sixeufb \scriptscriptfont\eufbfam=\fiveeufb}%
   \fi
 \fi
 \ifx\undefined\eusmfam
 \else \font@\eighteusm=eusm8 \font@\sixeusm=eusm6
   \ifsyntax@\else \addto\tenpoint{\textfont\eusmfam=\teneusm
       \scriptfont\eusmfam=\seveneusm \scriptscriptfont\eusmfam=\fiveeusm}%
     \addto\eightpoint{\textfont\eusmfam=\eighteusm
       \scriptfont\eusmfam=\sixeusm \scriptscriptfont\eusmfam=\fiveeusm}%
   \fi
 \fi
 \ifx\undefined\eusbfam
 \else \font@\eighteusb=eusb8 \font@\sixeusb=eusb6
   \ifsyntax@\else \addto\tenpoint{\textfont\eusbfam=\teneusb
       \scriptfont\eusbfam=\seveneusb \scriptscriptfont\eusbfam=\fiveeusb}%
     \addto\eightpoint{\textfont\eusbfam=\eighteusb
       \scriptfont\eusbfam=\sixeusb \scriptscriptfont\eusbfam=\fiveeusb}%
   \fi
 \fi
 \ifx\undefined\eurmfam
 \else \font@\eighteurm=eurm8 \font@\sixeurm=eurm6
   \ifsyntax@\else \addto\tenpoint{\textfont\eurmfam=\teneurm
       \scriptfont\eurmfam=\seveneurm \scriptscriptfont\eurmfam=\fiveeurm}%
     \addto\eightpoint{\textfont\eurmfam=\eighteurm
       \scriptfont\eurmfam=\sixeurm \scriptscriptfont\eurmfam=\fiveeurm}%
   \fi
 \fi
 \ifx\undefined\eurbfam
 \else \font@\eighteurb=eurb8 \font@\sixeurb=eurb6
   \ifsyntax@\else \addto\tenpoint{\textfont\eurbfam=\teneurb
       \scriptfont\eurbfam=\seveneurb \scriptscriptfont\eurbfam=\fiveeurb}%
    \addto\eightpoint{\textfont\eurbfam=\eighteurb
       \scriptfont\eurbfam=\sixeurb \scriptscriptfont\eurbfam=\fiveeurb}%
   \fi
 \fi
 \ifx\undefined\cmmibfam
 \else \font@\eightcmmib=cmmib8 \font@\sixcmmib=cmmib6
   \ifsyntax@\else \addto\tenpoint{\textfont\cmmibfam=\tencmmib
       \scriptfont\cmmibfam=\sevencmmib \scriptscriptfont\cmmibfam=\fivecmmib}%
    \addto\eightpoint{\textfont\cmmibfam=\eightcmmib
       \scriptfont\cmmibfam=\sixcmmib \scriptscriptfont\cmmibfam=\fivecmmib}%
   \fi
 \fi
 \ifx\undefined\cmbsyfam
 \else \font@\eightcmbsy=cmbsy8 \font@\sixcmbsy=cmbsy6
   \ifsyntax@\else \addto\tenpoint{\textfont\cmbsyfam=\tencmbsy
      \scriptfont\cmbsyfam=\sevencmbsy \scriptscriptfont\cmbsyfam=\fivecmbsy}%
    \addto\eightpoint{\textfont\cmbsyfam=\eightcmbsy
      \scriptfont\cmbsyfam=\sixcmbsy \scriptscriptfont\cmbsyfam=\fivecmbsy}%
   \fi
 \fi
 \let\topmatter\relax}
\def\chapterno@{\uppercase\expandafter{\romannumeral\chaptercount@}}
\newcount\chaptercount@
\def\chapter{\nofrills@{\afterassignment\chapterno@
                        CHAPTER \global\chaptercount@=}\chapter@
 \DNii@##1{\leavevmode\hskip-\leftskip
   \rlap{\vbox to\z@{\vss\centerline{\eightpoint
   \chapter@##1\unskip}\baselineskip2pc\null}}\hskip\leftskip
   \nofrills@false}%
 \FN@\next@}
\newbox\titlebox@
\def\title{\nofrills@{\uppercasetext@}\title@%
 \DNii@##1\endtitle{\global\setbox\titlebox@\vtop{\tenpoint\bf
 \raggedcenter@\ignorespaces
 \baselineskip1.3\baselineskip\title@{##1}\endgraf}%
 \ifmonograph@ \edef\next{\the\leftheadtoks}\ifx\next\empty
    \leftheadtext{##1}\fi
 \fi
 \edef\next{\the\rightheadtoks}\ifx\next\empty \rightheadtext{##1}\fi
 }\FN@\next@}
\newbox\authorbox@
\def\author#1\endauthor{\global\setbox\authorbox@
 \vbox{\tenpoint\smc\raggedcenter@\ignorespaces
 #1\endgraf}\relaxnext@ \edef\next{\the\leftheadtoks}%
 \ifx\next\empty\leftheadtext{#1}\fi}
\newbox\affilbox@
\def\affil#1\endaffil{\global\setbox\affilbox@
 \vbox{\tenpoint\raggedcenter@\ignorespaces#1\endgraf}}
\newcount\addresscount@
\addresscount@\z@
\def\address#1\endaddress{\global\advance\addresscount@\@ne
  \expandafter\gdef\csname address\number\addresscount@\endcsname
  {\vskip12\p@ minus6\p@\noindent\eightpoint\smc\ignorespaces#1\par}}
\def\email{\nofrills@{\eightpoint{\it E-mail\/}:\enspace}\email@
  \DNii@##1\endemail{%
  \expandafter\gdef\csname email\number\addresscount@\endcsname
  {\def\usualspace{{\it\enspace}}\smallskip\noindent\eightpoint\email@
  \ignorespaces##1\par}}%
 \FN@\next@}
\def\thedate@{}
\def\date#1\enddate{\gdef\thedate@{\tenpoint\ignorespaces#1\unskip}}
\def\thethanks@{}
\def\thanks#1\endthanks{\gdef\thethanks@{\eightpoint\ignorespaces#1.\unskip}}
\def\thekeywords@{}
\def\keywords{\nofrills@{{\it Key words and phrases.\enspace}}\keywords@
 \DNii@##1\endkeywords{\def\thekeywords@{\def\usualspace{{\it\enspace}}%
 \eightpoint\keywords@\ignorespaces##1\unskip.}}%
 \FN@\next@}
\def\thesubjclass@{}
\def\subjclass{\nofrills@{{\rm1980 {\it Mathematics Subject
   Classification\/} (1985 {\it Revision\/}).\enspace}}\subjclass@
 \DNii@##1\endsubjclass{\def\thesubjclass@{\def\usualspace
  {{\rm\enspace}}\eightpoint\subjclass@\ignorespaces##1\unskip.}}%
 \FN@\next@}
\newbox\abstractbox@
\def\abstract{\nofrills@{{\smc Abstract.\enspace}}\abstract@
 \DNii@{\setbox\abstractbox@\vbox\bgroup\noindent$$\vbox\bgroup
  \def\envir@{abstract}\advance\hsize-2\indenti
  \usualspace@{{\enspace}}\eightpoint \noindent\abstract@\ignorespaces}%
 \FN@\next@}
\def\endabstract{\par\unskip\egroup$$\egroup}
\def\widestnumber#1#2{\begingroup\let\head\null\let\subhead\empty
   \let\subsubhead\subhead
   \ifx#1\head\global\setbox\tocheadbox@\hbox{#2.\enspace}%
   \else\ifx#1\subhead\global\setbox\tocsubheadbox@\hbox{#2.\enspace}%
   \else\ifx#1\key\bgroup\let\endrefitem@\egroup
        \key#2\endrefitem@\global\refindentwd\wd\keybox@
   \else\ifx#1\no\bgroup\let\endrefitem@\egroup
        \no#2\endrefitem@\global\refindentwd\wd\nobox@
   \else\ifx#1\page\global\setbox\pagesbox@\hbox{\quad\bf#2}%
   \else\ifx#1\item\setboxz@h{#2}\global\rosteritemwd\wdz@
        \global\advance\rosteritemwd by.5\parindent
   \else\message{\string\widestnumber is not defined for this option
   (\string#1)}%
\fi\fi\fi\fi\fi\fi\endgroup}
\newif\ifmonograph@
\def\Monograph{\monograph@true \let\headmark\rightheadtext
  \let\varindent@\indent \def\headfont@{\bf}\def\proclaimfont@{\smc}%
  \def\demofont@{\smc}}
\let\varindent@\noindent
\newbox\tocheadbox@    \newbox\tocsubheadbox@
\newbox\tocbox@
\def\toc{\toc@{Contents}}
\def\newtocdefs{%
   \def \title##1\endtitle
       {\penaltyandskip@\z@\smallskipamount
        \hangindent\wd\tocheadbox@\noindent{\bf##1}}%
   \def \chapter##1{%
        Chapter \uppercase\expandafter{\romannumeral##1.\unskip}\enspace}%
   \def \specialhead##1\endspecialhead
       {\par\hangindent\wd\tocheadbox@ \noindent##1\par}%
   \def \head##1 ##2\endhead
       {\par\hangindent\wd\tocheadbox@ \noindent
        \if\notempty{##1}\hbox to\wd\tocheadbox@{\hfil##1\enspace}\fi
        ##2\par}%
   \def \subhead##1 ##2\endsubhead
       {\par\vskip-\parskip {\normalbaselines
        \advance\leftskip\wd\tocheadbox@
        \hangindent\wd\tocsubheadbox@ \noindent
        \if\notempty{##1}\hbox to\wd\tocsubheadbox@{##1\unskip\hfil}\fi
         ##2\par}}%
   \def \subsubhead##1 ##2\endsubsubhead
       {\par\vskip-\parskip {\normalbaselines
        \advance\leftskip\wd\tocheadbox@
        \hangindent\wd\tocsubheadbox@ \noindent
        \if\notempty{##1}\hbox to\wd\tocsubheadbox@{##1\unskip\hfil}\fi
        ##2\par}}}
\def\toc@#1{\relaxnext@
   \def\page##1%
       {\unskip\penalty0\null\hfil
        \rlap{\hbox to\wd\pagesbox@{\quad\hfil##1}}\hfilneg\penalty\@M}%
 \DN@{\ifx\next\nofrills\DN@\nofrills{\nextii@}%
      \else\DN@{\nextii@{{#1}}}\fi
      \next@}%
 \DNii@##1{%
\ifmonograph@\bgroup\else\setbox\tocbox@\vbox\bgroup
   \centerline{\headfont@\ignorespaces##1\unskip}\nobreak
   \vskip\belowheadskip \fi
   \setbox\tocheadbox@\hbox{0.\enspace}%
   \setbox\tocsubheadbox@\hbox{0.0.\enspace}%
   \leftskip\indenti \rightskip\leftskip
   \setbox\pagesbox@\hbox{\bf\quad000}\advance\rightskip\wd\pagesbox@
   \newtocdefs
 }%
 \FN@\next@}
\def\endtoc{\par\egroup}
\let\pretitle\relax
\let\preauthor\relax
\let\preaffil\relax
\let\predate\relax
\let\preabstract\relax
\let\prepaper\relax
\def\dedicatory #1\enddedicatory{\def\preabstract{{\medskip
  \eightpoint\it \raggedcenter@#1\endgraf}}}
\def\thetranslator@{}
\def\translator#1\endtranslator{\def\thetranslator@{\nobreak\medskip
 \line{\eightpoint\hfil Translated by \uppercase{#1}\qquad\qquad}\nobreak}}
\outer\def\endtopmatter{\runaway@{abstract}%
 \edef\next{\the\leftheadtoks}\ifx\next\empty
  \expandafter\leftheadtext\expandafter{\the\rightheadtoks}\fi
 \ifmonograph@\else
   \ifx\thesubjclass@\empty\else \makefootnote@{}{\thesubjclass@}\fi
   \ifx\thekeywords@\empty\else \makefootnote@{}{\thekeywords@}\fi
   \ifx\thethanks@\empty\else \makefootnote@{}{\thethanks@}\fi
 \fi
  \pretitle
  \ifmonograph@ \topskip7pc \else \topskip4pc \fi
  \box\titlebox@
  \topskip10pt
  \preauthor
  \ifvoid\authorbox@\else \vskip2.5pc plus1pc \unvbox\authorbox@\fi
  \preaffil
  \ifvoid\affilbox@\else \vskip1pc plus.5pc \unvbox\affilbox@\fi
  \predate
  \ifx\thedate@\empty\else \vskip1pc plus.5pc \line{\hfil\thedate@\hfil}\fi
  \preabstract
  \ifvoid\abstractbox@\else \vskip1.5pc plus.5pc \unvbox\abstractbox@ \fi
  \ifvoid\tocbox@\else\vskip1.5pc plus.5pc \unvbox\tocbox@\fi
  \prepaper
  \vskip2pc plus1pc
}
\def\document{\let\fontlist@\relax\let\alloclist@\relax
  \tenpoint}
\newskip\aboveheadskip       \aboveheadskip\bigskipamount
\newdimen\belowheadskip      \belowheadskip6\p@
\def\headfont@{\smc}
\def\penaltyandskip@#1#2{\relax\ifdim\lastskip<#2\relax\removelastskip
      \ifnum#1=\z@\else\penalty@#1\relax\fi\vskip#2%
  \else\ifnum#1=\z@\else\penalty@#1\relax\fi\fi}
\def\nobreak{\penalty\@M
  \ifvmode\def\penalty@{\let\penalty@\penalty\count@@@}%
  \everypar{\let\penalty@\penalty\everypar{}}\fi}
\let\penalty@\penalty
\def\heading#1\endheading{\head#1\endhead}
\def\subheading#1{\subhead#1\endsubhead}
\def\specialheadfont@{\bf}
\outer\def\specialhead{\par\penaltyandskip@{-200}\aboveheadskip
  \begingroup\interlinepenalty\@M\rightskip\z@ plus\hsize \let\\\linebreak
  \specialheadfont@\noindent\ignorespaces}
\def\endspecialhead{\par\endgroup\nobreak\vskip\belowheadskip}
\outer\def\head#1\endhead{\par\penaltyandskip@{-200}\aboveheadskip
 {\headfont@\raggedcenter@\interlinepenalty\@M
 \ignorespaces#1\endgraf}\nobreak
 \vskip\belowheadskip
 \headmark{#1}}
\let\headmark\eat@
\newskip\subheadskip       \subheadskip\medskipamount
\def\subheadfont@{\bf}
\outer\def\subhead{\nofrills@{.\enspace}\subhead@
 \DNii@##1\endsubhead{\par\penaltyandskip@{-100}\subheadskip
  \varindent@{\usualspace@{{\subheadfont@\enspace}}%
 \subheadfont@\ignorespaces##1\unskip\subhead@}\ignorespaces}%
 \FN@\next@}
\outer\def\subsubhead{\nofrills@{.\enspace}\subsubhead@
 \DNii@##1\endsubsubhead{\par\penaltyandskip@{-50}\medskipamount
      {\usualspace@{{\it\enspace}}%
  \it\ignorespaces##1\unskip\subsubhead@}\ignorespaces}%
 \FN@\next@}
\def\proclaimheadfont@{\bf}
\outer\def\proclaim{\runaway@{proclaim}\def\envir@{proclaim}%
  \nofrills@{.\enspace}\proclaim@
 \DNii@##1{\penaltyandskip@{-100}\medskipamount\varindent@
   \usualspace@{{\proclaimheadfont@\enspace}}\proclaimheadfont@
   \ignorespaces##1\unskip\proclaim@
  \sl\ignorespaces}%
 \FN@\next@}
\outer\def\endproclaim{\let\envir@\relax\par\rm
  \penaltyandskip@{55}\medskipamount}
\def\demoheadfont@{\it}
\def\demo{\runaway@{proclaim}\nofrills@{.\enspace}\demo@
     \DNii@##1{\par\penaltyandskip@\z@\medskipamount
  {\usualspace@{{\demoheadfont@\enspace}}%
  \varindent@\demoheadfont@\ignorespaces##1\unskip\demo@}\rm
  \ignorespaces}\FN@\next@}

\def\qed{\ifhmode\unskip\nobreak\fi\quad\ifmmode\square\else$\m@th\square$\fi}

\def\definition{\runaway@{proclaim}%
  \nofrills@{.\proclaimheadfont@\enspace}\definition@
        \DNii@##1{\penaltyandskip@{-100}\medskipamount
        {\usualspace@{{\proclaimheadfont@\enspace}}%
        \varindent@\proclaimheadfont@\ignorespaces##1\unskip\definition@}%
        \rm \ignorespaces}\FN@\next@}

\newdimen\rosteritemwd
\newcount\rostercount@
\newif\iffirstitem@
\let\plainitem@\item
\newtoks\everypartoks@
\def\par@{\everypartoks@\expandafter{\the\everypar}\everypar{}}
\def\roster{\edef\leftskip@{\leftskip\the\leftskip}%
 \relaxnext@
 \rostercount@\z@  
 \def\item{\FN@\rosteritem@}%
 \DN@{\ifx\next\runinitem\let\next@\nextii@\else
  \let\next@\nextiii@\fi\next@}%
 \DNii@\runinitem  
  {\unskip  
   \DN@{\ifx\next[\let\next@\nextii@\else
    \ifx\next"\let\next@\nextiii@\else\let\next@\nextiv@\fi\fi\next@}%
   \DNii@[####1]{\rostercount@####1\relax
    \enspace{\rm(\number\rostercount@)}~\ignorespaces}%
   \def\nextiii@"####1"{\enspace{\rm####1}~\ignorespaces}%
   \def\nextiv@{\enspace{\rm(1)}\rostercount@\@ne~}%
   \par@\firstitem@false  
   \FN@\next@}%
 \def\nextiii@{\par\par@  
  \penalty\@m\smallskip\vskip-\parskip
  \firstitem@true}%
 \FN@\next@}
\def\rosteritem@{\iffirstitem@\firstitem@false\else\par\vskip-\parskip\fi
 \leftskip3\parindent\noindent  
 \DNii@[##1]{\rostercount@##1\relax
  \llap{\hbox to2.5\parindent{\hss\rm(\number\rostercount@)}%
   \hskip.5\parindent}\ignorespaces}%
 \def\nextiii@"##1"{%
  \llap{\hbox to2.5\parindent{\hss\rm##1}\hskip.5\parindent}\ignorespaces}%
 \def\nextiv@{\advance\rostercount@\@ne
  \llap{\hbox to2.5\parindent{\hss\rm(\number\rostercount@)}%
   \hskip.5\parindent}}%
 \ifx\next[\let\next@\nextii@\else\ifx\next"\let\next@\nextiii@\else
  \let\next@\nextiv@\fi\fi\next@}

\newif\ifnextRunin@
\def\endroster{\relaxnext@
 \par\leftskip@  
 \penalty-50 \vskip-\parskip\smallskip  
 \DN@{\ifx\next\Runinitem\let\next@\relax
  \else\nextRunin@false\let\item\plainitem@  
   \ifx\next\par 
    \DN@\par{\everypar\expandafter{\the\everypartoks@}}%
   \else  
    \DN@{\noindent\everypar\expandafter{\the\everypartoks@}}%
  \fi\fi\next@}%
 \FN@\next@}
\newcount\rosterhangafter@
\def\Runinitem#1\roster\runinitem{\relaxnext@
 \rostercount@\z@ 
 \def\item{\FN@\rosteritem@}%
 \def\runinitem@{#1}%
 \DN@{\ifx\next[\let\next\nextii@\else\ifx\next"\let\next\nextiii@
  \else\let\next\nextiv@\fi\fi\next}%
 \DNii@[##1]{\rostercount@##1\relax
  \def\item@{{\rm(\number\rostercount@)}}\nextv@}%
 \def\nextiii@"##1"{\def\item@{{\rm##1}}\nextv@}%
 \def\nextiv@{\advance\rostercount@\@ne
  \def\item@{{\rm(\number\rostercount@)}}\nextv@}%
 \def\nextv@{\setbox\z@\vbox  
  {\ifnextRunin@\noindent\fi  
  \runinitem@\unskip\enspace\item@~\par  
  \global\rosterhangafter@\prevgraf}%
  \firstitem@false  
  \ifnextRunin@\else\par\fi  
  \hangafter\rosterhangafter@\hangindent3\parindent
  \ifnextRunin@\noindent\fi  
  \runinitem@\unskip\enspace 
  \item@~\ifnextRunin@\else\par@\fi  
  \nextRunin@true\ignorespaces}%
 \FN@\next@}
\def\footmarkform@#1{$\m@th^{#1}$}
\let\thefootnotemark\footmarkform@
\def\makefootnote@#1#2{\insert\footins
 {\interlinepenalty\interfootnotelinepenalty
 \eightpoint\splittopskip\ht\strutbox\splitmaxdepth\dp\strutbox
 \floatingpenalty\@MM\leftskip\z@\rightskip\z@\spaceskip\z@\xspaceskip\z@
 \leavevmode{#1}\footstrut\ignorespaces#2\unskip\lower\dp\strutbox
 \vbox to\dp\strutbox{}}}
\newcount\footmarkcount@
\footmarkcount@\z@
\def\footnotemark{\let\@sf\empty\relaxnext@
 \ifhmode\edef\@sf{\spacefactor\the\spacefactor}\/\fi
 \DN@{\ifx[\next\let\next@\nextii@\else
  \ifx"\next\let\next@\nextiii@\else
  \let\next@\nextiv@\fi\fi\next@}%
 \DNii@[##1]{\footmarkform@{##1}\@sf}%
 \def\nextiii@"##1"{{##1}\@sf}%
 \def\nextiv@{\iffirstchoice@\global\advance\footmarkcount@\@ne\fi
  \footmarkform@{\number\footmarkcount@}\@sf}%
 \FN@\next@}
\def\footnotetext{\relaxnext@
 \DN@{\ifx[\next\let\next@\nextii@\else
  \ifx"\next\let\next@\nextiii@\else
  \let\next@\nextiv@\fi\fi\next@}%
 \DNii@[##1]##2{\makefootnote@{\footmarkform@{##1}}{##2}}%
 \def\nextiii@"##1"##2{\makefootnote@{##1}{##2}}%
 \def\nextiv@##1{\makefootnote@{\footmarkform@{\number\footmarkcount@}}{##1}}%
 \FN@\next@}
\def\footnote{\let\@sf\empty\relaxnext@
 \ifhmode\edef\@sf{\spacefactor\the\spacefactor}\/\fi
 \DN@{\ifx[\next\let\next@\nextii@\else
  \ifx"\next\let\next@\nextiii@\else
  \let\next@\nextiv@\fi\fi\next@}%
 \DNii@[##1]##2{\footnotemark[##1]\footnotetext[##1]{##2}}%
 \def\nextiii@"##1"##2{\footnotemark"##1"\footnotetext"##1"{##2}}%
 \def\nextiv@##1{\footnotemark\footnotetext{##1}}%
 \FN@\next@}
\def\adjustfootnotemark#1{\advance\footmarkcount@#1\relax}
\def\footnoterule{\kern-3\p@
  \hrule width 5pc\kern 2.6\p@} 
\def\captionfont@{\smc}
\def\topcaption#1#2\endcaption{%
  {\dimen@\hsize \advance\dimen@-\captionwidth@
   \rm\raggedcenter@ \advance\leftskip.5\dimen@ \rightskip\leftskip
  {\captionfont@#1}%
  \if\notempty{#2}.\enspace\ignorespaces#2\fi
  \endgraf}\nobreak\bigskip}
\def\botcaption#1#2\endcaption{%
  \nobreak\bigskip
  \setboxz@h{\captionfont@#1\if\notempty{#2}.\enspace\rm#2\fi}%
  {\dimen@\hsize \advance\dimen@-\captionwidth@
   \leftskip.5\dimen@ \rightskip\leftskip
   \noindent \ifdim\wdz@>\captionwidth@ 
   \else\hfil\fi 
  {\captionfont@#1}\if\notempty{#2}.\enspace\rm#2\fi\endgraf}}
\def\@ins{\par\begingroup\def\vspace##1{\vskip##1\relax}%
  \def\captionwidth##1{\captionwidth@##1\relax}%
  \setbox\z@\vbox\bgroup} 
\def\block{\RIfMIfI@\nondmatherr@\block\fi
       \else\ifvmode\vskip\abovedisplayskip\noindent\fi
        $$\def\endblock{\par\egroup$$}\fi
  \vbox\bgroup\advance\hsize-2\indenti\noindent}
\def\endblock{\par\egroup}
\def\cite#1{{\rm[{\citefont@\m@th#1}]}}
\def\citefont@{\rm}
\def\refsfont@{\eightpoint}
\outer\def\Refs{\runaway@{proclaim}%
 \relaxnext@ \DN@{\ifx\next\nofrills\DN@\nofrills{\nextii@}\else
  \DN@{\nextii@{References}}\fi\next@}%
 \DNii@##1{\penaltyandskip@{-200}\aboveheadskip
  \line{\hfil\headfont@\ignorespaces##1\unskip\hfil}\nobreak
  \vskip\belowheadskip
  \begingroup\refsfont@\sfcode`.=\@m}%
 \FN@\next@}

\newbox\nobox@            \newbox\keybox@           \newbox\bybox@
\newbox\paperbox@         \newbox\paperinfobox@     \newbox\jourbox@
\newbox\volbox@           \newbox\issuebox@         \newbox\yrbox@
\newbox\pagesbox@         \newbox\bookbox@          \newbox\bookinfobox@
\newbox\publbox@          \newbox\publaddrbox@      \newbox\finalinfobox@
\newbox\edsbox@           \newbox\langbox@
\newif\iffirstref@        \newif\iflastref@
\newif\ifprevjour@        \newif\ifbook@            \newif\ifprevinbook@
\newif\ifquotes@          \newif\ifbookquotes@      \newif\ifpaperquotes@
\newdimen\bysamerulewd@
\setboxz@h{\refsfont@\kern3em}
\bysamerulewd@\wdz@
\newdimen\refindentwd
\setboxz@h{\refsfont@ 00. }
\refindentwd\wdz@
\outer\def\ref{\begingroup \noindent\hangindent\refindentwd
 \firstref@true \def\nofrills{\def\refkern@{\kern3sp}}%
 \ref@}
\def\ref@{\book@false \bgroup\let\endrefitem@\egroup \ignorespaces}
\def\moreref{\endrefitem@\endref@\firstref@false\ref@}%
\def\transl{\endrefitem@\endref@\firstref@false
  \book@false
  \prepunct@
  \setboxz@h\bgroup \aftergroup\unhbox\aftergroup\z@
    \def\endrefitem@{\unskip\refkern@\egroup}\ignorespaces}%
\def\emptyifempty@{\dimen@\wd\currbox@
  \advance\dimen@-\wd\z@ \advance\dimen@-.1\p@
  \ifdim\dimen@<\z@ \setbox\currbox@\copy\voidb@x \fi}
\let\refkern@\relax
\def\endrefitem@{\unskip\refkern@\egroup
  \setboxz@h{\refkern@}\emptyifempty@}\ignorespaces
\def\refdef@#1#2#3{\edef\next@{\noexpand\endrefitem@
  \let\noexpand\currbox@\csname\expandafter\eat@\string#1box@\endcsname
    \noexpand\setbox\noexpand\currbox@\hbox\bgroup}%
  \toks@\expandafter{\next@}%
  \if\notempty{#2#3}\toks@\expandafter{\the\toks@
  \def\endrefitem@{\unskip#3\refkern@\egroup
  \setboxz@h{#2#3\refkern@}\emptyifempty@}#2}\fi
  \toks@\expandafter{\the\toks@\ignorespaces}%
  \edef#1{\the\toks@}}
\refdef@\no{}{. }
\refdef@\key{[\m@th}{] }
\refdef@\by{}{}
\def\bysame{\by\hbox to\bysamerulewd@{\hrulefill}\thinspace
   \kern0sp}
\def\manyby{\message{\string\manyby is no longer necessary; \string\by
  can be used instead, starting with version 2.0 of \styname.STY}\by}
\refdef@\paper{\ifpaperquotes@``\fi\it}{}
\refdef@\paperinfo{}{}
\def\jour{\endrefitem@\let\currbox@\jourbox@
  \setbox\currbox@\hbox\bgroup
  \def\endrefitem@{\unskip\refkern@\egroup
    \setboxz@h{\refkern@}\emptyifempty@
    \ifvoid\jourbox@\else\prevjour@true\fi}%
\ignorespaces}
\refdef@\vol{\ifbook@\else\bf\fi}{}
\refdef@\issue{no. }{}
\refdef@\yr{}{}
\refdef@\pages{}{}
\def\page{\endrefitem@\def\pp@{\def\pp@{pp.~}p.~}\let\currbox@\pagesbox@
  \setbox\currbox@\hbox\bgroup\ignorespaces}
\def\pp@{pp.~}
\def\book{\endrefitem@ \let\currbox@\bookbox@
 \setbox\currbox@\hbox\bgroup\def\endrefitem@{\unskip\refkern@\egroup
  \setboxz@h{\ifbookquotes@``\fi}\emptyifempty@
  \ifvoid\bookbox@\else\book@true\fi}%
  \ifbookquotes@``\fi\it\ignorespaces}
\def\inbook{\endrefitem@
  \let\currbox@\bookbox@\setbox\currbox@\hbox\bgroup
  \def\endrefitem@{\unskip\refkern@\egroup
  \setboxz@h{\ifbookquotes@``\fi}\emptyifempty@
  \ifvoid\bookbox@\else\book@true\previnbook@true\fi}%
  \ifbookquotes@``\fi\ignorespaces}
\refdef@\eds{(}{, eds.)}
\def\ed{\endrefitem@\let\currbox@\edsbox@
 \setbox\currbox@\hbox\bgroup
 \def\endrefitem@{\unskip, ed.)\refkern@\egroup
  \setboxz@h{(, ed.)}\emptyifempty@}(\ignorespaces}
\refdef@\bookinfo{}{}
\refdef@\publ{}{}
\refdef@\publaddr{}{}
\refdef@\finalinfo{}{}
\refdef@\lang{(}{)}

\let\refdef@\relax 
\def\ppunbox@#1{\ifvoid#1\else\prepunct@\unhbox#1\fi}
\def\nocomma@#1{\ifvoid#1\else\changepunct@3\prepunct@\unhbox#1\fi}
\def\changepunct@#1{\ifnum\lastkern<3 \unkern\kern#1sp\fi}
\def\prepunct@{\count@\lastkern\unkern
  \ifnum\lastpenalty=0
    \let\penalty@\relax
  \else
    \edef\penalty@{\penalty\the\lastpenalty\relax}%
  \fi
  \unpenalty
  \let\refspace@\ \ifcase\count@,
\or;\or.\or 
  \or\let\refspace@\relax
  \else,\fi
  \ifquotes@''\quotes@false\fi \penalty@ \refspace@
}
\def\transferpenalty@#1{\dimen@\lastkern\unkern
  \ifnum\lastpenalty=0\unpenalty\let\penalty@\relax
  \else\edef\penalty@{\penalty\the\lastpenalty\relax}\unpenalty\fi
  #1\penalty@\kern\dimen@}
\def\endref{\endrefitem@\lastref@true\endref@
  \par\endgroup \prevjour@false \previnbook@false }
\def\endref@{%
\iffirstref@
  \ifvoid\nobox@\ifvoid\keybox@\indent\fi
  \else\hbox to\refindentwd{\hss\unhbox\nobox@}\fi
  \ifvoid\keybox@
  \else\ifdim\wd\keybox@>\refindentwd
         \box\keybox@
       \else\hbox to\refindentwd{\unhbox\keybox@\hfil}\fi\fi
  \kern4sp\ppunbox@\bybox@
\fi 
  \ifvoid\paperbox@
  \else\prepunct@\unhbox\paperbox@
    \ifpaperquotes@\quotes@true\fi\fi
  \ppunbox@\paperinfobox@
  \ifvoid\jourbox@
    \ifprevjour@ \nocomma@\volbox@
      \nocomma@\issuebox@
      \ifvoid\yrbox@\else\changepunct@3\prepunct@(\unhbox\yrbox@
        \transferpenalty@)\fi
      \ppunbox@\pagesbox@
    \fi 
  \else \prepunct@\unhbox\jourbox@
    \nocomma@\volbox@
    \nocomma@\issuebox@
    \ifvoid\yrbox@\else\changepunct@3\prepunct@(\unhbox\yrbox@
      \transferpenalty@)\fi
    \ppunbox@\pagesbox@
  \fi 
  \ifbook@\prepunct@\unhbox\bookbox@ \ifbookquotes@\quotes@true\fi \fi
  \nocomma@\edsbox@
  \ppunbox@\bookinfobox@
  \ifbook@\ifvoid\volbox@\else\prepunct@ vol.~\unhbox\volbox@
  \fi\fi
  \ppunbox@\publbox@ \ppunbox@\publaddrbox@
  \ifbook@ \ppunbox@\yrbox@
    \ifvoid\pagesbox@
    \else\prepunct@\pp@\unhbox\pagesbox@\fi
  \else
    \ifprevinbook@ \ppunbox@\yrbox@
      \ifvoid\pagesbox@\else\prepunct@\pp@\unhbox\pagesbox@\fi
    \fi \fi
  \ppunbox@\finalinfobox@
  \iflastref@
    \ifvoid\langbox@.\ifquotes@''\fi
    \else\changepunct@2\prepunct@\unhbox\langbox@\fi
  \else
    \ifvoid\langbox@\changepunct@1%
    \else\changepunct@3\prepunct@\unhbox\langbox@
      \changepunct@1\fi
  \fi
}
\outer\def\enddocument{%
 \runaway@{proclaim}%
\ifmonograph@ 
\else
 \nobreak
 \thetranslator@
 \count@\z@ \loop\ifnum\count@<\addresscount@\advance\count@\@ne
 \csname address\number\count@\endcsname
 \csname email\number\count@\endcsname
 \repeat
\fi
 \vfill\supereject\end}
\def\folio{{\foliofont@\ifnum\pageno<\z@ \romannumeral-\pageno
 \else\number\pageno \fi}}
\def\foliofont@{\eightrm}
\def\headlinefont@{\eightpoint}
\def\leftheadline{\rlap{\folio}\hfill \iftrue\topmark\fi \hfill}
\def\rightheadline{\hfill \expandafter
  \hfill \llap{\folio}}
\newtoks\leftheadtoks
\newtoks\rightheadtoks
\def\leftheadtext{\nofrills@{\uppercasetext@}\lht@
  \DNii@##1{\leftheadtoks\expandafter{\lht@{##1}}%
    \mark{\the\leftheadtoks\noexpand\else\the\rightheadtoks}
    \ifsyntax@\setboxz@h{\def\\{\unskip\space\ignorespaces}%
        \headlinefont@##1}\fi}%
  \FN@\next@}
\def\rightheadtext{\nofrills@{\uppercasetext@}\rht@
  \DNii@##1{\rightheadtoks\expandafter{\rht@{##1}}%
    \mark{\the\leftheadtoks\noexpand\else\the\rightheadtoks}%
    \ifsyntax@\setboxz@h{\def\\{\unskip\space\ignorespaces}%
        \headlinefont@##1}\fi}%
  \FN@\next@}
\headline={\def\chapter#1{\chapterno@. }%
  \def\\{\unskip\space\ignorespaces}\headlinefont@
  \ifodd\pageno \rightheadline \else \leftheadline\fi}
\def\NoRunningHeads{\global\runheads@false\global\let\headmark\eat@}

\def\logo@{\baselineskip2pc \hbox to\hsize{\hfil\eightpoint Typeset by
 \AmSTeX}}
\newif\iffirstpage@     \firstpage@true
\newif\ifrunheads@      \runheads@true
\output={\output@}
\def\output@{\shipout\vbox{%
 \iffirstpage@ \global\firstpage@false
  \pagebody \logo@ \makefootline%
 \else \ifrunheads@ \makeheadline \pagebody
       \else \pagebody \makefootline \fi
 \fi}%
 \advancepageno \ifnum\outputpenalty>-\@MM\else\dosupereject\fi}
\tenpoint
\catcode`\@=\active

\font\footfont=cmr10 at 8  pt
   
 \font\titfont=cmr10 at 18pt 
 \font\aufont= cmr10 at 14pt

 \magnification=1200
 \NoBlackBoxes





 \define\zz{ {{\bold{Z}_2}}}


 \define\calc{\Cal C}

 \define\calq{\Cal Q}
 
 \define\calz{\Cal Z}



 \define\cycd#1#2{{\calz}_{#1}(#2)}
 \define\cyc#1#2{{\calz}^{#1}(#2)}
 \define\cych#1{{{\calz}^{#1}}}
 \define\cycp#1#2{{\calz}^{#1}(\bbp(#2))}
 \define\cyf#1#2{{\cyc{#1}{#2}}^{fix}}
 \define\cyfd#1#2{{\cycd{#1}{#2}}^{fix}}

 \define\crl#1#2{{\calz}_{\bbr}^{#1}{(#2)}}

 \define\crd#1#2{\widetilde{\calz}_{\bbr}^{#1}{(#2)}}

 \define\crld#1#2{{\calz}_{#1,\bbr}{(#2)}}
 \define\crdd#1#2{\widetilde{\calz}_{#1,{\bbr}}{(#2)}}

 \define\crlh#1{{\calz}_{\bbr}^{#1}}
 \define\crdh#1{{\widetilde{\calz}_{\bbr}^{#1}}}
 \define\cyav#1#2{{{\cyc{#1}{#2}}^{av}}}
 \define\cyavd#1#2{{\cycd{#1}{#2}}^{av}}

 \define\cyaa#1#2{{\cyc{#1}{#2}}^{-}}
 \define\cyaad#1#2{{\cycd{#1}{#2}}^{-}}

 \define\cyq#1#2{{\calq}^{#1}(#2)}
 \define\cyqd#1#2{{\calq}_{#1}(#2)}

 \define\cqt#1#2{{\calz}_{\bbh}^{#1}{(#2)}}
 \define\cqtav#1#2{{\calz}^{#1}{(#2)}^{av}}
 \define\cqtrd#1#2{\widetilde{\calz}_{\bbh}^{#1}{(#2)}}

 \define\cyct#1#2{{\calz}^{#1}(#2)_\zz}
 \define\cyft#1#2{{\cyc{#1}{#2}}^{fix}_\zz}
\define\cxg#1#2{G^{#1}_{\bbc}(#2)}
 \define\reg#1#2{G^{#1}_{\bbr}(#2)}

 \define\cyaat#1#2{{\cyc{#1}{#2}}^{-}_\zz}



 \define\fflag#1#2{{#1}={#1}_{#2} \supset {#1}_{{#2}-1} \supset
 \ldots \supset {#1}_{0} }
 
 \define\vect#1{ {\Cal{V}ect}_{#1}}


 \define\chv#1#2#3{{\calc}^{#1}_{#2}(#3)}
 \define\chvd#1#2#3{{\calc}_{#1,#2}(#3)}
 \define\chm#1#2{{\calc}_{#1}(#2)}



 \define\Claim#1{\subheading{Claim #1}}

 \define\xrightarrow#1{\overset{#1}\to{\rightarrow}}


\font \fr = eufm10

\hfuzz1pc 


\define\gl{\text{\fr gl}}

\define\bbz{\Bbb Z}

\define\bbr{\Bbb R}
\define\bbc{\Bbb C}
\define\bbh{\Bbb H}
\define\bbp{\Bbb P}

\define\cu{\Cal U}

\define\ca{\Cal A}
\define\cc{\Cal C}
\define\cd{\Cal D}
\define\ce{\Cal E}
\define\ch{\Cal H}
\define\cb{\Cal B}
\define\cj{\Cal J}

\define\cg{\Cal G}

\define\cf{\Cal F}

\define\cv{\Cal V}

\define\cs{\Cal S}
\define\cz{\Cal Z}
\define\co#1{\Cal O_{#1}}
\define\ct{\Cal T}
\define\ci{\Cal I}
\define\cR{\Cal R}


\define\a{\alpha}
\redefine\b{\beta}

\redefine\d{\delta}
\define\r{\rho}
\define\s{\sigma}
\define\z{\zeta}

\redefine\D{\Delta}

\define\p#1{{\bbp}^{#1}}

\define\equdef{\overset\text{def}\to=}

\define\blbx{\hfill  $\square$}
\redefine\qed{\blbx}

\define\pf{\subheading{Proof}}
\define\Lemma#1{\subheading{Lemma #1}}
\define\Theorem#1{\subheading{Theorem #1}}
\define\Prop#1{\subheading{Proposition #1}}
\define\Cor#1{\subheading{Corollary #1}}
\define\Note#1{\subheading{Note #1}}
\define\Def#1{\subheading{Definition #1}}
\define\Remark#1{\subheading{Remark #1}}
\define\Ex#1{\subheading{Example #1}}
\define\arr{\longrightarrow}

\define\Map{\operatorname{Map}}

\redefine\Xi{X_{\infty}}

\define\jac#1#2{\left(\!\!\!\left(
\frac{\partial #1}{\partial #2}
\right)\!\!\!\right)}
\define\restrict#1{\left. #1 \right|_{t_{p+1} = \dots = t_n = 0}}

\define\SP#1#2{{\roman SP}^{#1}(#2)}

\define\coc#1#2#3{\cc^{#1}(#2;\, #3)}
\define\zoc#1#2#3{\cz^{#1}(#2;\, #3)}
\define\zyc#1#2#3{\cz^{#1}(#2 \times #3)}

\define\ar#1{\overset{#1}\to{\longrightarrow}}

\define\th#1#2{{\Bbb H}^{#1}(#2)}
\define\hth#1#2{\widehat{\Bbb H}^{#1}(#2)}


\define\bad#1#2{\cf_{#2}(#1)}

\define\pch#1{\bbp_{\bbc}(\bbh^{#1})}

\def\l{\ell}

\def\I#1#2{I_{#1, #2}}

\def\A{A}

\def\D{D}

\def\L{L}

\def\z2t{\text{$\bbz_2\ct$}}

\def\hH{\widehat{I\!\!H}}

\def\<{\left<}
\def\>{\right>}
\def\[{\left[}
\def\]{\right]}

\def\wt{\widetilde}
\def\vf{\varphi}

\def\AA{4}

\redefine\and{\qquad\text{and}\qquad}

\def\dbar{\overline{\partial}}

\def\th{^{\text{th}}} 
\def\for{\qquad\text{ for }}

\def\wh{\widehat}\
\def\Image{\operatorname{Image}}

\define\Dp#1{{\cd'}^{#1}(X,p)}
\define\Ep#1{{\ce}^{#1}(X,p)}
\define\Zp#1{{\cz}^{#1}_{\bbz}(X,p)}
\define\Hop#1#2{{H}^{#1}_{#2}(X,p)}
\define\Hp#1{{\hH}^{#1}(X,p)}
\define\ddp{\overline{d}}
\define\Dcx#1{\underline{\bbz}_{\cd}(#1)}
\define\HD#1#2{H^{#1}_{\cd}(X,\bbz(#2))}

       
\redefine\hH{\widehat{\bold H}}
\def\vff{e}
\def\ov{\overline}
\def\In{\ci}

\centerline{\bf \titfont }
\vskip .3in

\centerline{\titfont D-bar Sparks}
\vskip .1in
\centerline{\aufont by}

\vskip .1in
\centerline{\aufont F. Reese Harvey and H. Blaine Lawson, Jr.
\footnote{\footfont Research
 partially supported by the NSF}}

\vskip .5in
 
\centerline{\bf Abstract} \medskip
  \font\abstractfont=cmr10 at 10 pt

{{\parindent=.9in\narrower\abstractfont \noindent 
A $\dbar$-analogue of differential characters for complex manifolds
is introduced and studied using a new theory of
homological spark complexes. Many essentially different spark complexes
are shown to have isomorphic  groups of spark classes. This has many
consequences:  It leads to an analytic representation of 
$\co{}^\times$-gerbes with connection,  it yields  a soft
resolution of the  sheaf $\co{}^\times$ by currents on the manifold, 
and more generally it gives a Dolbeault-Federer representation of Deligne
cohomology as the cohomology of certain complexes of currents. 

}}
\medskip   
{{\parindent=.9in\narrower\abstractfont \noindent 
It is shown that the $\dbar$-spark classes $\hH^*(X)$ carry a functorial
ring structure.  Holomorphic bundles have Chern classes in this theory
which refine the integral classes and satisfy  Whitney duality. A version
of Bott vanishing for holomorphic foliations is proved  in this context.

}}

\vfill\eject\ \qquad
\vskip .5in            
               
\centerline{\bf Table of Contents}
\medskip   
       
\hskip 1in  0.        Introduction
 
\hskip 1in  1.      Homological Spark Complexes        
  
\hskip 1in  2.      Quasi-Isomorphisms of Spark Complexes  

\hskip 1in  3.      \v Cech-Dolbeault Sparks

\hskip 1in  4.     $\co{}^{\times}$-Gerbes with Connection

\hskip 1in  5.      The Case of  Degree 0

\hskip 1in  6.      The Case of  Degree 1

\hskip 1in  7.      $\dbar$-Sparks
 
\hskip 1in  8.     \v Cech-Dolbeault Hypersparks
  
\hskip 1in 9.       A Current Resolution of the Sheaf $\co{}^\times$
  
\hskip 1in  \!\! 10.  Ring Structure  

\hskip 1in  \!\! 11.  Functoriality
 
\hskip 1in  \!\! 12.  Refined Chern Classes for Holomorphic Bundles
 
\hskip 1in  \!\! 13.   Bott Vanishing for Holomorphic Foliations 
 
\hskip 1in  \!\! 14.  Generalizations and the Relation to Deligne Cohomology
 
\hskip 1in  \!\! 15.  Theorems in the General Setting
 
\hskip 1in  \!\! Appendix A.  Hypercohomology and Cone Complexes
  
\hskip 1in  \!\! Appendix B.  Non-smoothness of Integrally Flat Currents

\vfill\eject

\centerline{ 0. INTRODUCTION}\medskip
In 1972 Jeff Cheeger and Jim Simons developed a theory of differential characters
with applications to secondary characteristic invariants, conformal geometry, the theory
of foliations, and much more [4], [23],  [3], [5]. Over the intervening years the importance  of this theory has steadily grown -- it is now relevant to many areas of mathematics and mathematical physics.
The purpose of this paper is to establish a ``$\dbar$-analogue'' of this theory with new
applications to geometry and analysis. In fact for each integer $p\geq 0$
we establish a ring functor $ \Hp {*} $ on the category of complex manifolds and holomorphic maps,
with two natural transformations (graded ring homomorphisms)
$$
\d_1: \Hp {*-1} \ \arr\ \Zp {*} \and \d_2: \Hp {*-1}\ \arr\ H^{*}(X;\bbz)
$$
 where $\Zp {*}$ denote the closed differential forms, with certain
 integrality properties, in the truncated Dolbeault complex $\bigoplus_{j<p}\ce^{j,*-j}(X)$.
The kernel of $\d_1$ is Deligne cohomology ($\cong H^{*-1}(X,\co{}^{\times})$ when $p=1$).
The kernel of $\d_2$ corresponds to classes with smooth representatives.
These homomorphisms sit in an exact, functorial $3\times 3$ grid.

For the sake of exposition we shall restrict attention primarily to the case $p=1$,
which is quite natural and already of considerable interest. The results and  proofs
for $p>1$, which are strictly analogous, will be presented in \S14.

We first  introduce a number of distinct $\dbar$-spark complexes and examine them in
 detail in low degrees. One  of these is related to $\dbar$-gerbes with connection,
 another  is more analytical and concerns $\dbar$-equations relating
 smooth to rectifiable currents.  All these complexes are shown to have isomorphic
groups of spark classes ${\hH}^{*}(X,1)$. 
From this, one  establishes the  ring structure and functoriality. In certain presentations
of the theory, there are quite explicit formulas for the product.

This theory is then applied to define refined Chern classes 
$\hat{d}_k(E) \in \hH^{2k-1}(X,\co{}^\times)$ for holomorphic vector bundles  $E$.
These classes are constructed by using a hermitian metric and its associated
canonical connection, but as in ordinary Chern-Weil theory, the result is independent
of the choice.  The associated {\sl total class}  $ {\hat d}(E) = 1+{\hat d}_1(E)+{\hat d}_2(E)+\dots$
satisfies the duality
$$
 {\hat d}(E\oplus F)\  =\  {\hat d}(E) *  {\hat d}(F).  
$$
Furthermore, $\d_2\circ  {\hat d}(E) = c(E) \in H^*(X;\bbz)$ is  the ordinary total integral
Chern class of $E$.  These results, and their extension to  $p>1$, give a 
representation of Chern-Weil type for the Deligne characteristic classes of $E$.

In a subsequent section we consider holomorphic foliations and establish the following 
analogue of the Bott Vanishing Theorem.

\Theorem{}  {\sl Let $N$ be a holomorphic bundle of rank $q$ on a
complex manifold $X$.  If $N$ is (isomorphic to) the normal bundle of a
holomorphic foliation of $X$, then for every  polynomial $P$ of pure
cohomology degree $k>2q$, the associated refined Chern class satisfies
$$
P({\wh d}_1(N),...,{\wh d}_q(N)) \ \in\ H^{2k-1}(X;\,\bbc^{\times}) \subset
H^{2k-1}(X;\,{\co{}}^{\times})
$$
 }

A nice aspect of this theory is its natural  presentation of 
Deligne cohomology in terms of forms and currents.  In fact, we shall 
construct several different spark complexes, each yielding the same ring of 
spark classes containing Deligne cohomology. Hence, we also obtain several other
geometric constructions of Deligne cohomology.

We note that the standard short exact sequence
for the Deligne group  $\HD {k} p $ sits as the left column in our $3\times 3$-grid. 
When $X$ is compact Kaehler
and $k=2p$, this is the  sequence: 
$
0\to  \cj_p(X) \to \HD {2p} p  \to  \text{Hdg}^{p,p}(X)\to0
$
where $\cj_p(X)$ is Griffiths' intermediate Jacobian.
(See   the diagrams following Proposition 14.4.) 
In our context it is trivial to see that every holomorphic cycle of codimension-$p$ determines 
a class in $ \HD {2p} p$.  Furthermore, one   sees in a similar way that {\sl every maximally
complex cycle $M$ of codimension $2p+1$ determines a class   $[M]\in \HD{2p+1} p$, and if 
$M$ bounds a holomorphic chain of (complex) codimension $p$, then $[M]=0$} (cf. Prop. 14.6).

Recall that a powerful property of  cohomology theory
is its broad range of distinct formulations.  Some presentations of the theory
make it easy to compute, while others, such as de Rham theory, current
theory or harmonic theory, lead to non-trivial assertions in analysis.
Differential characters are similar in nature.  There are many distinct formulations:
some relatively simple and computable, and others rather more complicated, involving
\v Cech-deRham complexes [19] or complexes of currents [20], [10], [14].   These latter approaches  relate differential characters for example to refined characteristic classes
for singular connections [16], [17], [22], to Morse Theory [18] and to harmonic theory.

 In [19] the  homological apparatus of {\sl spark complexes} was introduced
to  establish the equivalence of the many approaches to differential characters. 
In \S 1 below, that apparatus is generalized to treat a  wide range
of  situations.  In particular, on a complex manifold one can replace the 
deRham component of differential character theory by the $\dbar$-complex of 
$(0,q)$-forms, or more generally, the deRham complex truncated at level $p$. 
The machine developed in \S 1  is of   independent interest and applies
to a broad range of interesting situations.

Using this machine we  show that   a variety of $\dbar$-spark complexes 
are equivalent and therefore lead to isomorphic groups of spark classes. Thus we are able to
relate the classes of $\dbar$-gerbes with connection, which are defined within the \v Cech-Dolbeault double complex,  to classes of currents satisfying a 
$\dbar$-spark equation of the form
$$
\dbar a \ =\ \phi -\Psi(R)
$$
where $\phi$ is a  smooth $\dbar$-closed $(0,q)$-form and $\Psi(R)$ is the $(0,q)$-component
of a rectifiable cycle $R$ of codimension $q$. These two are related by a larger enveloping 
complex which contains them both.  

Similar remarks apply to the case where one truncates
at level $p>1$.

Other geometrically motivated spark complexes will
be studied in forthcoming papers.

\vfill\eject

\centerline{ 1. HOMOLOGICAL SPARK COMPLEXES}\medskip

We begin with a generalization of the homological algebra
introduced in [19].

\Def{1.1} A {\bf homological spark complex} is a triple of cochain
complexes \linebreak $(F^*, E^*, I^*)$ together with morphisms
$$
I^*\ @>{\Psi}>>\  F^*\  \supset\  E^*
$$
such that:
\roster
\item "(i)"\quad  $\Psi(I^k)\cap E^k\ =\ \{0\}$ \ \ for $k>0$, \medskip
\item "(ii)" \quad  $H^*(E) \cong H^*(F)$, and\medskip 
\item "(iii)" \quad $\Psi: I^o\to F^0$ is injective.
\endroster
\medskip

 Note that   
$$
\Psi(I^0)\cap E^0\ \subset\  Z^0(E) \ =\ H^0(E) \ =\ H^0(F)
$$
since for any $a\in \Psi(I^0)\cap E^0$, we have  $da \in \Psi(I^1)\cap
E^1=\{0\}$

\medskip

\Def{1.2} In a given  spark complex  $(F^*,E^*,I^*)$ a {\bf spark of
degree $k$}  is a pair  
$$
(a,r)\in F^k\oplus I^{k+1}
$$
which satisfies the {\bf spark equation}
\medskip
\roster
\item "(i)"\quad $da \ =\ \vff - \Psi(r) \qquad\text{ for some }
              e\in   E^{k+1}$, and  
\medskip 
\item "(ii)" \quad $dr\ =\ 0$.
\endroster
\medskip
\noindent
The group of sparks of degree $k$ is denoted by
$\cs^k=\cs^k(F^*,E^*,I^*)$. 
Note that by 1.1\,(i)  \medskip

\, (iii)\ \ \quad $d\vff=0$.
\medskip
\noindent

\medskip

\Def{1.3} Two sparks $(a,r), (a'r')\in \cs^k(F^*,E^*,I^*)$ are 
{\bf equivalent} if there exists a pair
$$
(b,s)\in F^{k-1}\oplus I^{k}
$$
\medskip
\roster
\item "(i)" \quad $a -a' \ =\ db + \Psi(s)  $ 
\medskip
\item "(ii)"\quad $r-r'\ =\ -ds$.
\endroster
\medskip\noindent
The set of equivalence classes is called the {\bf group of spark
classes of degree $k$} associated to the given spark complex and will be
denoted by $\hH^k(F^*, E^*, I^*)$ or simply  $\hH^k$ when the
complex in question is evident.
Note that $\hH^{-1} = H^0(I)$.
\medskip

We now derive the fundamental exact sequences associated to a 
homological spark complex $(F^*,E^*,I^*)$.   Let
$Z^k(E)=\{\vff\in E^k : d\vff=0\}$ and set
$$
Z_I^k(E)\ \equiv\ \{\vff \in Z^k(E) \ :\
[\vff]=\Psi_*(\rho)  \quad\text{for some}\ \  \rho \in H^k(I)\}
\tag1.1
$$ 
where $[\vff]$ denotes the class of $\vff$ in $H^k(E)\cong H^k(F)$.

\Lemma{1.4}
{\sl  There exist well-defined surjective homomorphisms:
$$
\hH^k\ @>{\d_1}>>\ Z_I^{k+1}(E)         \and
\hH^k\ @>{\d_2}>>\ H^{k+1}(I)       
\tag1.2
$$
given on any representing spark $(a, r)\in \cs^k$ by
$$
\d_1(a,r) = \vff\and \d_2(a,r)=[r]
$$
 where 
$da =\vff-\Psi(r)$ as in 1.2\,(i).  }

\pf
It is straightforward to see that these maps are well-defined.
To see that $\d_1$ is surjective, consider $e\in  Z_I^{k+1}(E)$. By
definition there exists  $r\in I^{k+1}$ with $dr=0$ such that 
$[e]=\Psi_*[r]\in H^{k+1}(F)\cong H^{k+1}(E)$. Hence, there exists
$a\in F^k$ with $da=e-\Psi(r)$ and $\d_1(a,r)= e$ as desired.
The map $\d_2$ is surjective because if $[r]\in H^{k+1}(I)$, then $\Psi(r)$
represents a class in $H^{k+1}(F) \cong H^{k+1}(E)$.  Picking a
representative $e\in E^{k+1}$ of this class yields a spark $(a,r)$ with
$da=e-\Psi(r)$.\qed

\Lemma{1.5}
{\sl  Let $\hH^k_E$ = $\ker \d_2$. Then}
$$
\hH^k_E\ =\ E^k/Z^k_I(E)
$$
\pf
Suppose $\a\in \hH^k$ is represented by the spark $(a,r)$. Then 
$da=e-\Psi(r)$ with $e\in E^{k+1}$, and $dr=0$. Now $\d_2\a=0$ means that 
$[r]=0\in H^{k+1}(I)$, i.e., $r= -ds$ for some $s\in I^k$. The equivalent
spark $(a-\Psi(s),0)$ satisfies $d(a-\Psi(s))=e$.  Since $H^*(E)\cong
H^*(F)$, we know by  [19, Lemma 1.5] that  we can find $b\in F^{k-1}$
so that $a-\Psi(s)+db\in E^k$.  This proves that each $\a\in\hH^k_E$ has a
representative of the form $(a,0)$ with $a\in E^k$.
If $(a,0)$ is equivalent to 0, then $a=db +\Psi(s)$ for some $b\in F^{k-1}$
and some $s\in I^k$ with $ds=0$.  That is, $a\in Z_I^{k}(E)$.\qed

\Def{1.6}
Associated to any spark complex $(F^*, E^*, I^*)$ is the 
{\bf cone complex} $(G^*, D)$ defined by setting
$$
\aligned
G^k\ \equiv \ F^k &\oplus I^{k+1}\qquad\qquad k\geq -1  \\
D(a,r)\  =\ (da &+ \Psi(r), -dr)
\endaligned
\tag1.3$$
Note that there is a short exact sequence of complexes
$$
0\ \to\ \ F^*\ \to \ G^*\ \to\ I^*(1)\ \to\ 0
\tag1.4$$
where $I^k(1) \equiv I^{k+1}$.  The morphism $\Psi$ defines a
chain map of degree 1: 
$$
F^*\quad \ @<{\Psi}<<\ \quad I^*(1) 
$$
which induces the connecting homomorphisms in the associated long 
exact sequence in cohomology.

\Prop{1.7} {\sl
There are two fundamental short exact sequences:}
$$\aligned
0\ \arr\ H^k(G)\ \arr\ &\hH^k\ @>{\d_1}>>\ \ \ Z_I^{k+1}(E)\ \arr\ \   0 \\
0\ \, \ \arr\ \  \hH^k_E\ \ \arr\ &\hH^k\ @>{\d_2}>>\ \ H^{k+1}(I)\ \arr\ \ 0
\endaligned
\tag1.5
$$
\pf  This follows immediately from Lemmas 1.4, 1.5 and Definition 1.6.
\qed
\medskip

Consider the homomorphism $\Psi_*:H^k(I)\to H^k(F)\cong H^k(E)$, and define
$$
H^k_I(E)\ \equiv \ \text{Image}\{\Psi_*\}
\and
Ker^k(I)\ \equiv \ \text{ker}\{\Psi_*\}
$$ 
The exact sequences above fit into the following $3\times 3$ commutative
grid. 

\Prop{1.8}{\sl Associated to any spark complex $(F^*, E^*, I^*)$ 
is the commutative diagram  
$$
\CD
\ @.  0   @.  0  @.  0 @. \  \\
@. @VVV @VVV @VVV @.  \\
0 @>>> \frac {H^k(E)}{H^k_I(E)}  @>>> \hH^k_E  @>>> dE^k    @>>> 0 \\
@. @VVV @VVV @VVV @.  \\
0 @>>> H^k(G)  @>>> \hH^k  @>{\d_1}>>  Z_I^{k+1}(E)    @>>> 0 \\
@. @VVV @V{\d_2}VV @VVV @.  \\
0 @>>> Ker^{k+1}(I)  @>>> H^{k+1}(I) @>{\Psi_*}>> H^{k+1}_I(E)   @>>> 0
\\ @. @VVV @VVV @VVV @.  \\
\ @.  0   @.  0  @.  0 @. \  
\endCD
\tag1.6
$$

whose rows and columns are exact.}
\ 

\vskip .3in

\centerline{ 2.\  QUASI-ISOMORPHISMS OF SPARK COMPLEXES}

\vskip .1in

In this section we introduce a useful criterion for showing that
two spark complexes have isomorphic groups of spark classes.

\Def{2.1} Two spark complexes $(F^*, E^*, I^*)$ and $
(\overline{F}^*, \overline{E}^*,\overline{I}^*)$ are 
{\bf quasi-isomorphic} if there exists a commutative diagram of
morphisms  
$$
\CD
\overline{I}^* @>{\overline{\Psi}}>> \overline{F}^*  \quad \supset\quad
\overline{E}^*\\ 
@A{\psi}AA   \cup\qquad\qquad \!\! \| \\
{I}^* @>{\Psi}>> {F}^*  \quad \supset\quad
{E}^*
\endCD
$$
inducing an isomorphism
$$
\psi^*:H^*(I)\ @>{\cong}>>\ H^*(\overline{I})
\tag2.1
$$

\Prop{2.2} {\sl A quasi-isomorphism of spark complexes 
$(F^*, E^*, I^*)$ and $
(\overline{F}^*, \overline{E}^*,\overline{I}^*)$
induces an isomorphism 
$$
\hH^*(F^*, E^*, I^*) \ \cong \ 
\hH^*(\overline{F}^*, \overline{E}^*,\overline{I}^*)
$$
of the associated groups of spark classes.  In fact, it induces an
isomorphism of the grids (1.6) associated to the two complexes.}

\pf There is evidently a mapping $\hH^*(F^*, E^*, I^*) \to 
\hH^*(\overline{F}^*, \overline{E}^*,\overline{I}^*)$. To see that it is
onto, consider a spark $(\overline{a},\overline{r})\in 
\overline{\cs}^k$ with $d\overline{r}=0$ and $d\overline{a}={\vff}-
\overline{\Psi}(\overline{r})$  where 
$\vff \in \overline{E}^{k+1} = {E}^{k+1}$.
By the RHS of (2.1) there exists an element $\ov{s}\in \ov{I}^k$
such that $\ov r = \psi (r) - d\ov s$ for some $r\in I^{k+1}$ with
$dr=0$.  Hence, $ d \ov a = \vff  -\ov \Psi (\ov r) 
=\vff -\ov \Psi (\psi( r)) + \ov\Psi (d\ov s)
=\vff -\Psi ( r) + d \ov\Psi (\ov s)$, and we have
$$
d\{\ov a -\ov\Psi (\ov s)\}\ =\ \vff-\Psi(r) \ \in \ F^{k+1}.
$$
It follows (see [19, Lemma 1.5]) that there exists $\ov b\in 
{\ov F}^{k-1}$ such that 
$$
a\ \equiv\ \ov a -\ov\Psi (\ov s) +d \ov b \ \in\ F^k.
$$
Consequently, $(\ov a, \ov r)$ is equivalent to   
$(\ov a -\ov\Psi (\ov s) +d \ov b, \ov r + d\ov s) = (a, \psi(r))$ 
which is the image of the spark $(a,r)\in\cs^k= \cs^k(F^*,E^*, I^*)$,
(note that $da=e-\Psi(r)$).  Hence,  the map is onto as claimed.

We now prove that the mapping is injective.  Suppose that the image
of $(a,r) \in\cs^k$ is equivalent to 0 in $\ov{\cs}^k$. This means that
there exists a pair $(\ov b, \ov s)\in {\ov F}^{k-1}\oplus {\ov I}^k$ such
that
$$
a\ =\ d\ov b - \ov \Psi (\ov s)
\and
\psi(r)\ =\ d\ov s.
\tag2.2$$

Since $\psi_*: H^*( I)\to H^*(\ov I)$ is an isomorphism , the RHS of
(2.2) implies that there exists $s\in I^k$ such that
$
r=ds
$ 
Hence $(a,r)$ is equivalent in $\cs^k$ to $(a',0)$ where $a'=a+\Psi(s)$.
The triviality condition (2.2) now becomes
$$
a'\ =\ d\ov b - \ov \Psi (\ov s')
\and
0\ =\ d\ov s'.
$$
Again since $\psi_*: H^*( I)\to H^*(\ov I)$ is an isomorphism, there
exists  $s'\in I^k$ with $ds'=0$ such that $\psi(s')= \ov s' +d \ov t$
where  $\ov t \in {\ov I}^{k-1}$.  Hence, $\ov\Psi(\ov
s')=\ov\Psi(\psi(s')-d\ov t) = \Psi(s')-d\ov\Psi(\ov t)$, and we conclude
that $$
a' +\Psi(s')\ =\ d({\ov b}')
$$
where $\ov b'=\ov b +\ov\Psi (\ov t) \in {\ov F}^{k-1}$. Since   
$H^*(\ov F)\cong H^*(F)$, there exists an element $b\in F^{k-1}$ with 
$$
a' +\Psi(s')\ =\ db.
$$
Hence, $(a',0)$ is equivalent in $\cs^k$ to $(0,0) = (a'+\Psi(s')-db,
ds')$. \qed

\medskip

\Remark{2.3} Proposition 2.2 can be strengthened by replacing
the inclusion $F^*\subset {\ov F}^*$ with a strong homological 
equivalence $F^*\to {\ov F}^*$ which induces an isomorphism
$E^*\to {\ov E}^*$.  A {\sl strong homological equivalence}
is a chain map whose kernel and cokernel are acyclic.
\medskip

Proposition 2.2  unifies the diverse theories
of sparks, gerbes, differential characters and holonomy maps (see [19]).
In what follows we examine a ``$\dbar$- analogue'' of these objects.

\vskip.3in


\centerline{ 3. \v CECH-DOLBEAULT SPARKS}\vskip.1in

Suppose $X$ is a complex manifold of dimension $n$, and let
$$
0\ \to\ \co {}\ \to\ \ce^{0,0}\ @>{\dbar}>> \ \ce^{0,1}\
@>{\dbar}>> \ \ce^{0,2} \ @>{\dbar}>> \ \dots\ @>{\dbar}>> \ \ce^{0,n}
$$
denote the Dolbeault resolution of the sheaf $\co{}=\co{X}$ of
holomorphic functions on $X$ by smooth  $(0,q)$-forms.  Suppose
$\cu=\{U_\a\}_{\a\in A}$ is a covering of $X$ by Stein open sets
such that all finite intersections $U_{\a_1}\cap\dots\cap U_{\a_N}$
are contractible, and consider the double \v Cech-Dolbeault complex 
$$
\CD
C^{0}(\cu, \,\ce^{0,n}) @>{\d}>> C^{1}(\cu, \,\ce^{0, n}) @>{\d}>> 
C^{2}(\cu, \,\ce^{0,n}) @>{\d}>> \dots @>{\d}>> C^{n}(\cu, \,\ce^{0, n}) \\
@A{\dbar}AA @A{\dbar}AA @A{\dbar}AA @. @A{\dbar}AA  \\
\vdots @. \vdots @. \vdots @. \ @. \vdots @.     \\
@A{\dbar}AA @A{\dbar}AA @A{\dbar}AA @. @A{\dbar}AA  \\
C^{0}(\cu, \,\ce^{0,2}) @>{\d}>> C^{1}(\cu, \,\ce^{0, 2}) @>{\d}>> 
C^{2}(\cu, \,\ce^{0,2}) @>{\d}>> \dots @>{\d}>> C^{n}(\cu, \,\ce^{0, 2}) \\
@A{\dbar}AA @A{\dbar}AA @A{\dbar}AA @. @A{\dbar}AA  \\
C^{0}(\cu, \,\ce^{0,1}) @>{\d}>> C^{1}(\cu, \,\ce^{0, 1}) @>{\d}>> 
C^{2}(\cu, \,\ce^{0,1}) @>{\d}>> \dots @>{\d}>> C^{n}(\cu, \,\ce^{0, 1}) \\
@A{\dbar}AA @A{\dbar}AA @A{\dbar}AA @. @A{\dbar}AA  \\
C^{0}(\cu, \,\ce^{0,0}) @>{\d}>> C^{1}(\cu, \,\ce^{0, 0}) @>{\d}>> 
C^{2}(\cu, \,\ce^{0,0}) @>{\d}>> \dots @>{\d}>> C^{n}(\cu, \,\ce^{0, 0}) \\
\endCD
$$
There are two edge complexes:
$$
\{\ker(\d) \text{  on the left column}\}\ \cong \ Z^0(\cu, \,\ce^{0,*})
\ \cong\   H^0(\cu, \,\ce^{0,*})\  \cong\  \ce^{0,*}(X)
\tag3.1$$
the standard {\sl smooth Dolbeault complex} with differential $\dbar$, and 
$$
\{\ker(\dbar) \text{  on the bottom row}\}\ \cong \ C^{*}(\cu, \,\co{})
\tag3.2$$
the standard {\sl \v Cech complex} with coefficients in the sheaf $\co{}$.
There is also the  {\sl total complex} $(F^*, D)$ where
$$
F^k\ \equiv\ \bigoplus_{p+q=k} C^p(\cu, \ce^{0,q})  \and
D\ \equiv \  (-1)^p\d+\dbar
\tag3.3$$
\Lemma{3.1} {\sl The inclusions of the edge complexes
$$
\ce^{0,*}(X)\ \subset F^*\and C^{*}(\cu, \,\co{})\ \subset F^*
$$
each induce an isomorphism in cohomology.}

\pf
There are two spectral sequences converging to $H^*(F)$ cf. [12, 4.8].
For the first one, $'E_1^{*,*}$ = the  $\d$-cohomology of $C^*(\cu,
\ce^{0,*})$, which, by the fineness of the sheaves $\ce^{0,*}$,
reduces to the edge complex (3.1) in the far left column and
zero elsewhere.

For the second sequence, $''E_1^{*,*}$ = the  $\dbar$-cohomology of $C^*(\cu,
\ce^{0,*})$, which, since each $\ce^{0,q}$ is fine and each
$U_{\a_1}\cap\dots\cap U_{\a_1}$ is Stein, reduces to the edge complex (3.2)
in the bottom row and zero elsewhere. 
\qed
\medskip

Associated to this double complex are many interesting spark complexes.
The most basic is the following.

\Def{3.2} The {\sl \v Cech-Dolbeault spark complex } is the  homological
spark complex $(F^*, E^*,I^*)$ where $F^*$ is defined by (3.3), 
$$
E^*=\ce^{0,*}(X)\and
I^*\equiv C^*(\cu,\bbz)@>{\Psi}>>C^*(\cu,\co{})
$$
and where $\Psi$ is the chain map associated to the inclusion of the
sheaf of locally constant integer-valued functions into $\co{}$.
The  groups of spark classes associated to this complex will be denoted
by $\hH^{0,*}(X)$.
\medskip 

\Prop{3.3} {\sl The  diagram (1.6) for the groups  $\hH^{0,*}(X)$
can be written as
$$
\CD
\ @.  0   @.  0  @.  0 @. \  \\
@. @VVV @VVV @VVV @.  \\
0 @>>> \frac {H^k(X,\co{})}{H^k_{\bbz}(X,\co{})}  @>>> \hH^{0,k}_\infty(X) 
@>>> \dbar\ce^{0,k}(X)    @>>> 0 \\ @. @VVV @VVV @VVV @.  \\
0 @>>> H^k(X, \co{}^{\times})  @>>> \hH^{0,k}(X)  @>{\d_1}>> 
\cz^{0,k+1}_{\bbz}(X)    @>>> 0 \\
 @. @VVV @V{\d_2}VV @VVV @.  \\
0 @>>> \d \{H^k(X, \co{}^{\times})\}  @>>> H^{k+1}(X;\bbz) @>>>
H^{k+1}_{\bbz}(X, \co{})   @>>> 0 \\ @. @VVV @VVV @VVV @.  \\
\ @.  0   @.  0  @.  0 @. \  
\endCD
\tag3.4
$$
where  $H^{k}_{\bbz}(X, \co{})= \Image\{H^{k}(X;\bbz)
\to H^{k}(X,\co{})\}$ and $\d:H^k(X, \co{}^{\times})\to  H^{k+1}(X;\bbz)$ is
the coboundary map  coming from the exponential sequence
$0\to\bbz\to\co{}\to\co{}^{\times}\to0$ and where 
$\cz^{0,k}_{\bbz}(X)$ denotes $\dbar$-closed $(0,k)$-forms representing 
classes in $H^{k}_{\bbz}(X, \co{})\cong H^{0,k}_{\bbz}(X)$.
 The group $\hH^{0,k}_\infty(X)\cong 
\ce^{0,k}(X)/\cz^{0,k}_{\bbz}(X)$ consists exactly of the
spark  classes representable by smooth $(0,k)$-forms.}

 \pf  This is essentially straightforward. The only real point is the
identification of $H^*(G)$  in (1.6) with $H^*(X, \co{}^{\times})$.
This is done explicitly in Theorem 9.1.

\Remark{3.4} If $\cu'$ is another open cover   having the same properties as $\cu$, then the 
associated groups of \v Cech-Dolbeault spark classes and 3$\times$3 diagrams are isomorphic.
This follows from the fact that if $\cu'$ is also a refinement of $\cu$, the associated 
 \v Cech-Dolbeault spark complexes are quasi-isomorphic.

\Remark{3.5} In  the above discussion one can replace 
$\bbz$ with any subring $\Lambda\subset\bbc$.  In this case
the sheaf $\co{}^{\times}$ is replaced by $\co{}/\Lambda$

\vskip.3in

\centerline{ 4.  \  $\co{}^{\times}$-GERBES WITH CONNECTION}
\vskip.1in

To each \v Cech-Dolbeault spark  one can associate  a slightly more
geometric object called an ``$\co{}^{\times}$-grundle''or a
(generalized) ``$\co{}^{\times}$-gerbe with connection''.  Under the
correspondence, spark equivalence becomes a certain ``gauge equivalence'',
so the groups of spark classes discussed in \S 3 become gauge equivalence
classes of  $\co{}^{\times}$-gerbes with connection.
These objects have very nice interpretations in low dimensions. 

An ``$\co{}^{\times}$-grundle'' of degree $k$ is obtained from a \v
Cech-Dolbeault spark $A$ of degree $k$ by replacing the bottom
component $A^{k,0}$ with its exponential
$$
g_{\a_1\dots\a_{k+1}}\ =\ e^{2\pi i A_{\a_1\dots\a_{k+1}}}
$$
Therefore, $g \in C^k(\cu, \ce^{\times})$ where $\ce^{\times}$ denotes the sheaf
of smooth $\bbc^*$-valued functions.  The bottom component of the 
spark equation:
$$
\d A^{k,0}\ =\ r
$$
where $r_{\a_1\dots\a_{k+1}} \in \bbz$, implies that $g$ is a cocycle, that
is 
$$
\d g\ =\ 0.
$$
\Def{4.1} An {\bf $\co{}^{\times}$-grundle of degree $k$} is a pair
$(A,g)$ where $g \in C^k(\cu, \ce^{\times})$ satisfies 
$
\d g=0
$
 and  $A\in \bigoplus_{p+q=k, q>0} C^p(\cu, \ce^{0,q})$ satisfies
$$\aligned
\dbar A^{0,k}\ &=\ \vf \in \ce^{0,k+1}(X) \\
\dbar A^{1,k-1} +(-1)^k \d A^{0,k} \ \   &=\  0  \\
\dbar A^{2,k-2} +(-1)^{k-1}  \d A^{1,k-1} \ &=\ 0  \\
\vdots\qquad\qquad &\ \  \\
\dbar A^{k-1,1} + \d A^{k-2,2} \ &=\ 0  \\
\frac1 {2\pi i}\frac {\dbar g}g - \d A^{k-1,1} \ &=\ 0  
\endaligned
$$

\Def{4.2} Two $\co{}^{\times}$-grundles $(A,g),(\wt A,\wt g)$ of degree $k$
are {\bf gauge equivalent} if there exists a pair
$(B,h)$ where $h \in C^{k-1}(\cu, \ce^{\times})$ satisfies 
$$
\d h = g {\wt g}^{-1}
$$
and  $B\in \bigoplus_{p+q=k-1, q>0} C^p(\cu, \ce^{0,q})$ satisfies
$$\aligned
 A^{0,k} - {\wt A}^{0,k}\ &=\ \dbar B^{0,k-1} \\
 A^{1,k-1} - {\wt A}^{1,k-1}\ &=\ \dbar B^{1,k-2} +(-1)^{k-1}  \d B^{0,k-1}  \\
 A^{2,k-2} - {\wt A}^{2,k-2}\ &=\ \dbar B^{2,k-3} +(-1)^{k-2}  \d B^{1,k-2}  \\
&\ \  \vdots\\
 A^{k-1,1} - {\wt A}^{k-1,1}\ &=\  \frac1 {2\pi i}\frac {\dbar h}h -\d
B^{k-2,1}  \endaligned
$$

The group of gauge equivalence classes of $\co{}^{\times}$-grundles
is denoted by $\hH^{0,k}_{\text{grundle}}(X)$. The following proposition
is straightforward to check.

\Prop{4.3} {\sl The correspondence above induces an isomorphism}
$$
\hH^{0,k}(X)\ \cong\ \hH^{0,k}_{\text{grundle}}(X) 
$$

In the next two sections we examine these groups $\hH^{0,k}(X)$ in low
degrees where they  have interesting interpretations analogous to those
of differential characters [4], [20].

\vskip.3in

\centerline{ 5.  \  THE CASE OF DEGREE 0}
\medskip

From Proposition 4.3 one immediately deduces the isomorphism
$$
\hH^{0,0}(X) \ \cong\ \Map(X, \bbc^\times)
$$ 
with the space of $C^\infty$-maps to $\bbc^\times$.  The two fundamental
exact sequences (1.5) become 
$$
\aligned
0\ \to\ \co{}^{\times}(X)\ \to\ &\Map(X, \bbc^\times)\ @>{\dbar}>>\
\cz^{0,1}_{\bbz}(X) \ \to\ 0 \\
0\ \to\ \Map(X, \bbc)/\bbz\ @>{\exp}>>\ &\Map(X, \bbc^\times)\ @>{}>>\
H^1(X;\,\bbz)\ \to\ 0
\endaligned
$$

\vskip.3in

\centerline{6.  \ THE CASE OF DEGREE 1}
\medskip

From Proposition 4.3 one sees that an element of $\hH^{0,1}(X)$
is represented by a pair $(A,g)$ where $g_{\a \b}:U_{\a}\cap 
U_{\b}\to \bbc^\times$ is a \v Cech 1-cocycle, 
and $A_\a\in \ce^{0,1}(U_\a)$ are (0,1)-forms satisfying
$$
A_\a - A_\b\ =\ \frac 1 {2\pi i} \frac{\dbar g_{\a\b}}{g_{\a\b}}
\qquad\text{ in } U_{\a}\cap  U_{\b}.
$$
Thus $g$ gives the data for a complex line bundle $L$ on $X$ and $A$
represents the (0,1)-part of a connection on $L$. These objects have 
intrinsic interest.

\Def{6.1}  Let $L\to X$ be a smooth complex line bundle on a complex
manifold $X$. Then a {\bf $\dbar$-connection} on $L$ is a linear mapping
$$
\dbar_A:\ce^{0,0}(X, L)\ \arr\ \ce^{0,1}(X,L)
$$
from smooth sections of $L$ to smooth (0,1)-forms with values in $L$    
such that
$$
\dbar_A(f\s)\ =\ (\dbar f)\otimes \s + f\dbar_A(\s)
$$
for all $f\in C^\infty(X)$ and $\s\in \ce^{0,0}(X, L)$. 
    Two  $\dbar$-connections $\dbar_A$, $\dbar_B$ are said to be {\bf gauge
equivalent} if $\dbar_B=g\circ \dbar_A \circ g^{-1}$ for some bundle
isomorphism $g:L\to L$.

\medskip

A $\dbar$-connection extends naturally to define a Dolbeault
sequence
$$
\ce^{0,0}(X, L)\ @>{\dbar_A}>>\ \ce^{0,1}(X, L)\ @>{\dbar_A}>>\ 
\ce^{0,2}(X, L)\ @>{\dbar_A}>>\ \dots\ @>{\dbar_A}>>\  \ce^{0,n}(X, L)
\tag6.1$$ 
for non-holomorphic bundles $L$.  In general $\dbar_A^2$ is not zero. 
However, on a section $\s\in \ce^{0,0}(X, L)$  one has
$$
\dbar_A^2 \s \ =\ \vf_A  \otimes\s
$$
for $\vf_A \in \ce^{0,2}(X)$.  The form $\vf_A$ is called the
{\bf $\dbar$-curvature} of the connection. It has the following properties.

\Prop{6.2} {\sl  
\roster
\item $\vf_A$ depends only on the gauge equivalence
class of the $\dbar$-connection.

\item   $L$ admits a $\dbar$-connection 
with curvature $\vf_A\equiv 0$ if and only if $L$ is (smoothly) equivalent
to a holomorphic line bundle.  In fact each $\dbar_A$ with $\vf_A\equiv0$
determines a unique holomorphic structure on $L$.
 \endroster}

\pf  The first assertion is a calculation.  The second is a standard
consequence of the Newlander-Nirenberg Theorem.

\Prop {6.3} {\sl Let $X$ be a complex manifold.  Then there is a natural
isomorphism: }
$$
\hH^{0,1}(X)\ \cong\ \frac{
\{\text{ smooth complex line bundles with
$\dbar$-connection on } X\} }
{\text{gauge equivalence}}
$$
\pf Exercise.
\medskip

Note that the $\dbar$-spark equation in this case can be written
$$
\dbar (A,g)\ =\ \vf_A - C_1(L)
\tag6.2$$
where $\vf_A$ is the curvature of the $\dbar$-connection and $C_1(L)$
is a \v Cech representative of the first Chern class of the line bundle $L$.
In this case the  two fundamental exact sequences (1.5) become:
$$
\aligned
0\ \to\ H^1(X, \co{}^{\times})\ \to\ &\hH^{0,1}(X)\ @>{\vf_A}>>\
\cz^{0,2}_{\bbz}(X) \ \to\ \  0 \\
0\ \to\ \hH^{0,1}_{\text{triv}}(X)\ @>{}>>\ &\hH^{0,1}(X)\ @>{c_1}>>\
H^2(X;\,\bbz)\ \to\ 0
\endaligned
$$
where $\hH^{0,1}_{\text{triv}}(X)\subset\hH^{0,1}(X)$ is the subgroup
where the line bundles are topologically trivial.

\vskip.3in

\centerline{ 7.  \  $\dbar$-SPARKS}
\medskip

Suppose $X$ is a complex manifold of dimension $n$. Let $\ce^{p,q}(X)$
denote the smooth forms of bidegree $(p,q)$ on $X$ and 
${\cd'}^{p,q}(X)\supset \ce^{p,q}(X)$ the {\sl generalized forms}
or {\sl currents} of bidegree  $(p,q)$ 
on $X$.  (Recall that by definition ${\cd'}^{p,q}(X)$ is the topological 
dual space to $\cd^{n-p,n-q}(X)$ with the $C^\infty$-topology.)\  The
standard Dolbeault decomposition induces a decomposition
$$
{\cd'}^k(X)\ =\ \bigoplus_{p+q=k} {\cd'}^{p,q}(X)
\tag7.1
$$
of the  currents of degree $k$ on $X$.

Let $\cR^k(X)\subset {\cd'}^k(X)$ denote the group of locally rectifiable
currents of degree $k$ (and dimension $2n-k$) on $X$ (cf. [8]). A current
$T\in {\cd'}^k(X)$ is called locally {\bf integral} if 
$
T \in\cR^k(X)
$
and $dT\in\cR^{k+1}(X)$.  The group of locally integral
currents of degree $k$ on $X$ is denoted $\In^k(X)$.
The complex of sheaves $0\to \bbz \to \In^*$ is a fine resolution of 
$\bbz$ and one has a natural isomorphism
$$
H^*(\In)\ \cong\ H^*(X;\,\bbz)
\tag7.2
$$
For compact manifolds $X$ this is a basic result of Federer and Fleming [9], [8, 4.4.5].

\Note{7.1}  In all that follows, one could replace the complex
$\In^*(X)$ with the subcomplex $\In^*_{\text{deR}}(X)$ of de Rham {\sl
integral chain currents}, which are  defined by
integration over (locally finite) $C^\infty$ singular chains 
with $\bbz$-coefficients (cf. [7]).
This  follows   from Proposition 2.2 and the fact
that $  H^*(\In^*_{\text{deR}}(X))\ \cong\ H^*(X;\,\bbz)$.  Readers 
unfamiliar with integral currents may want to replace them 
with integral chain currents.

\Def{7.2}  By the {\bf $\dbar$-spark complex} on $X$ we mean the triple
$(F^*,E^*,I^*)$ where
$$\aligned
F^*\ &\equiv\ {\cd'}^{0,*}(X)\\
E^*\ &\equiv\ \ce^{0,*}(X)\\
I^*\ &\equiv\ \In^{*}(X)
\endaligned$$
with 
$$
\ce^{0,*}(X)\ \subset\ {\cd'}^{0,*}(X)
\and
\Psi:\In^{*}(X)\ @>{}>>\ {\cd'}^{0,*}(X)
$$
where the inclusion that of smooth forms into the space of forms 
with distribution coefficients, and where for $r\in\In^k(X)$, we have 
$
\Psi(r) \equiv r^{0,k},
$
the $(0,k)\th$ component of $r$ in the decomposition (7.1).  
The isomorphism
$H^*(E)\cong H^*(F)$ is standard. The fact that $E^k\cap \Psi(I^k)=\{0\}$
for $k>0$ is proved in Appendix B.

Thus a  {\bf $\dbar$-spark of degree $k$} is a pair
$(a,r)\in{\cd'}^{0,k}(X)\oplus \In^{k+1}(X)$ such that $d r=0$ and
$a$ satisfies the {\bf $\dbar$-spark equation:}
$$
\dbar a \ =\ \vf - r^{0,k+1}
\tag7.3
$$
The set of these will be denoted by $\cs^{0,k}(X)$.  A $\dbar$-spark 
$(a,r) \in \cs^{0,k}(X)$ is equivalent to zero if there exits
$(b,s)\in{\cd'}^{0,k-1}(X)\oplus \In^{k}(X)$ such that 
$$
a\ =\ \dbar b+s^{0,k+1}   \and   r\ =\ -ds
\tag7.4$$
The set of  spark classes is denoted  by
$$
\hH^{0,k}_{\dbar-\text{spark}}(X)\ \equiv \ \cs^{0,k}(X)/\sim.
$$

\vskip.3in

\centerline{ 8.  \  \v CECH-DOLBEAULT HYPERSPARKS}
\medskip

Suppose $X$ is a complex manifold of dimension $n$, and let
$$
0\ \to\ \co {}\ \to\ {\cd'}^{0,0}\ @>{\dbar}>> \ {\cd'}^{0,1}\
@>{\dbar}>> \ {\cd'}^{0,2} \ @>{\dbar}>> \ \dots\ @>{\dbar}>> \ {\cd'}^{0,n}
$$
denote the Serre resolution of the sheaf $\co{}=\co{X}$ by
generalized  $(0,q)$-forms.  Suppose $\cu=\{U_\a\}_{\a\in A}$ is a covering
of $X$ as in \S 3, and consider the double
\v Cech-Serre complex  
$$
\CD
\vdots @.\vdots @.\vdots @. \\
@A{\dbar}AA @A{\dbar}AA @A{\dbar}AA  \\
C^{0}(\cu, \,{\cd'}^{0,2}) @>{\d}>> C^{1}(\cu, \,{\cd'}^{0, 2}) @>{\d}>> 
C^{2}(\cu, \,{\cd'}^{0,2}) @>{\d}>> \dots   \\
@A{\dbar}AA @A{\dbar}AA @A{\dbar}AA   \\
C^{0}(\cu, \,{\cd'}^{0,1}) @>{\d}>> C^{1}(\cu, \,{\cd'}^{0, 1}) @>{\d}>> 
C^{2}(\cu, \,{\cd'}^{0,1}) @>{\d}>>  \\
@A{\dbar}AA @A{\dbar}AA @A{\dbar}AA  \\
C^{0}(\cu, \,{\cd'}^{0,0}) @>{\d}>> C^{1}(\cu, \,{\cd'}^{0, 0}) @>{\d}>> 
C^{2}(\cu, \,{\cd'}^{0,0}) @>{\d}>> \dots  \\
\endCD
$$
There are two edge complexes:
$$\aligned
\{\ker(\d) \text{  on the left column}\}
\ &\cong\   Z^0(\cu, \,{\cd'}^{0,*})\  \cong\  {\cd'}^{0,*}(X)\ \supset\ 
 \ce^{0,*}(X),\\
\{\ker(\dbar) \text{  on the } & \text{bottom row}\}\ \cong \ C^{*}(\cu,
\,\co{})
\endaligned$$
where in the second line we are using the standard regularity  
results for $\dbar$. Consider the  {\sl total  complex} $({\ov F}^*, D)$
where  $$
{\ov F}^k\ \equiv\ \bigoplus_{p+q=k} C^p(\cu, {\cd'}^{0,q})  \and
D\ \equiv \  (-1)^p\d+\dbar
\tag8.1$$
\Lemma{8.1} {\sl The inclusions of the edge complexes
$$
\ce^{0,*}(X)\ \subset {\ov F}^*\and C^{*}(\cu, \,\co{})\ \subset {\ov F}^*
$$
each induce an isomorphism in cohomology.}

\pf The argument is essentially the same as the one given for Lemma 3.1.\qed
\medskip

Let $\In^q$ denote the sheaf of germs of locally integrally flat currents
of degree $q$ on $X$. The standard boundary operator gives a resolution
$$
0\to\bbz\to \In^0 \to \In^1 \to\In^2 \to \In^3\to \dots
$$
of the constant sheaf $\bbz$ (cf. [21, Appendix A]).
Consider the double complex  
$$
\CD
\vdots @.\vdots @.\vdots @. \\
@A{d}AA @A{d}AA @A{d}AA  \\
C^{0}(\cu, \,\In^2) @>{\d}>> C^{1}(\cu, \,\In^2) @>{\d}>> 
C^{2}(\cu, \,\In^2) @>{\d}>> \dots   \\
@A{d}AA @A{d}AA @A{d}AA   \\
C^{0}(\cu, \,\In^1) @>{\d}>> C^{1}(\cu, \,\In^1) @>{\d}>> 
C^{2}(\cu, \,\In^1) @>{\d}>>  \\
@A{d}AA @A{d}AA @A{d}AA  \\
C^{0}(\cu, \,\In^0) @>{\d}>> C^{1}(\cu, \,\In^0) @>{\d}>> 
C^{2}(\cu, \,\In^0) @>{\d}>> \dots  \\
\endCD
$$
There are two edge complexes:
$$\aligned
\{\ker(\d) \text{  on the left column}\}
\ \cong\   Z^0(\cu, \,&\In^*)\  \cong\  \In^*(X),\\
\{\ker(d) \text{  on the }  \text{bottom row}\}\ &\cong \ C^{*}(\cu,
\,\bbz)
\endaligned$$
 Consider the  {\sl total  complex} $({\ov I}^*, D)$
where  $$
{\ov I}^k\ \equiv\ \bigoplus_{p+q=k} C^p(\cu, \In^q)  \and
D\ \equiv \  (-1)^p\d+\dbar
\tag8.2$$
\Lemma{8.2} {\sl The inclusions of the edge complexes
$
\In^*(X)\ \subset {\ov I}^*$ and $C^{*}(\cu, \,\bbz)\ \subset {\ov I}^*
$
 induce  isomorphisms }
$$
H^*(\In^*(X)) \cong H^*(C^{*}(\cu, \,\bbz)) \cong H^*({\ov I}^*) \cong 
H^*(X;\,\bbz).
$$
\pf This follows as before since the sheaves $\In^q$  are
acyclic. (See [21, Appendix A].)\qed
\medskip

Inclusion $\In^q\subset {\cd'}^q$ followed by
projection ${\cd'}^q\to {\cd'}^{0,q}$ gives a morphism of sheaves
$
\In^q \to {\cd'}^{0,q}
$ 
which induces a mapping of double complexes
$$
\Psi:C^p(\cu,\In^q) \to C^p(\cu,{\cd'}^{0,q})
\tag8.3
$$

\Def{8.3}  By the {\sl \v Cech-Dolbeault hyperspark complex }   on $X$ we
mean the  homological
spark complex $({\ov F}^*, {\ov E}^*,{\ov I}^*)$  where
$
{\ov E}^*\ \equiv \ce^{0,*}(X)\subset {\ov F}^*
$
is given as in Lemma 8.1 and 
$
\Psi:{\ov I}^*\to {\ov F}^*
$
is defined by (8.3).  The associated group of spark classes will be denoted
by $\hH^{0,k}_{\text{hyperspark}}(X)$

\medskip

\Theorem{8.4} {\sl The \v Cech-Dolbeault hyperspark complex
is quasi-isomorphic to the (smooth) \v Cech-Dolbeault spark complex 
and also to the (analytic) $\dbar$-spark complex of \S 8.  Hence there are
natural isomorphisms}
$$
\hH^{0,k}(X)\ \cong\ 
\hH^{0,k}_{\dbar-\text{spark}}(X)\ \cong\ 
\hH^{0,k}_{\text{hyperspark}}(X)
$$

\pf  Let $(F^*,E^*,I^*)$ denote the \v Cech-Dolbeault spark complex 
defined in 3.2. The natural inclusion $\ce^{0,*}\subset{\cd'}^{0,*}$
of the smooth forms into the generalized forms gives an inclusion
$F^*\to {\ov F}^*$ which is  the identity on $E^*$. We define 
a mapping $I^*\to {\ov I}^*$ by the inclusions
$$
I^*\equiv C^*(\cu,\bbz)\  \subset\  C^*(\cu,\In^0)\  \subset\ \ 
C^*(\cu,\In^*) \equiv {\ov I}^*.
$$
This commutes with the maps $\Psi$, and by Lemma 8.2 it induces an
isomorphism  $H^*(I)\cong H^*({\ov I})$. This establishes the first
assertion.

Now let $(F^*,E^*,I^*)$ denote the   $\dbar$-spark complex 
defined in 7.2, and consider the inclusion 
$$
F^k={\cd'}^{0,k}(X)\  =\  Z^0(\cu, {\cd'}^{0,k})\ \subset\ 
\bigoplus_{p+q=k}C^p(\cu, {\cd'}^{0,q})\ =\ 
{\ov F}^k
$$ 
which is essentially the identity on $E^*$. One sees that under this
inclusion  $I^*=F^*\cap {\ov I}^*$ is the vertical edge complex of the 
double complex ${\ov I}^*$. Hence, by Lemma 8.2 one has that
 $H^*(I^*)\cong H^*({\ov I}^*)$, and
so this inclusion is a quasi-isomorphism as claimed. \qed

\vskip.3in

\centerline{ 9. A CURRENT RESOLUTION OF THE SHEAF \  $\co{}^\times$}
\medskip
One of the first consequences of Theorem 8.4 is the following  theorem
of deRham-Federer type, giving an isomorphism between 
$H^*(X, \co{}^\times)$ and the cohomology of a certain complex of currents
on a complex manifold  $X$.  

Consider the complex of sheaves on $X$:
$$
0\to\cg^{-1} \to \cg^0 @>{\ov D}>>\cg^1 @>{\ov D}>>\cg^2 @>{\ov
D}>>\dots
\tag9.1
$$
with
$$
\cg^q \ \equiv\ {\cd'}^{0,q}\oplus \In^{q+1}
\and
\ov D (a,r)\ \equiv\ (\dbar a + r^{0,q+1}, -dr).
$$
where ${\cd'}^{0,*}$ and $\In^*$ are as in \S 8.
Let 
$$
G^q(X) \equiv {\cd'}^{0,q}(X)\oplus \In^{q+1}(X)
$$
be the group of global sections and consider the associated complex
$$
G^{-1}(X)@>{\ov D}>> G^0(X)@>{\ov D}>> G^1(X)@>{\ov D}>> 
G^2(X)@>{\ov D}>> \dots 
\tag9.2$$
Note  that $\cg^{-1}=\In^{0}$ is the sheaf of locally integrable 
$\bbz$-valued functions, and $\ov D : \cg^{-1}\to \cg^{0}={\cd'}^0\oplus
\In^1$ is given by $\ov D (f)=(f,-df)$.  We define 
$$
\widetilde{\cg}^0\ \equiv\ \cg^0/\ov D \cg^{-1}
$$
to be the quotient sheaf and claim that the kernel of $\ov D$ on 
$\widetilde{\cg}^0$ is the sheaf $\co{}^\times$. To see this fix a
contractible open set $U$ and consider an element $(a,r)\in
G^0(U)$ with $\ov D(a,r) = (\dbar a+r^{0,1}, -dr)=(0,0)$. Since
$H^1(\ci^*(U))=H^1(U,\bbz)=0$, we have $r=df$ for $f:U\to \bbz$ 
locally integrable. Hence $(a,r)\sim (a+f, r-df) =(\widetilde a,0)$ and
$\dbar \widetilde a = \bar a + [df]^{0,1}=\dbar a+r^{0,1}=0$. Two such
holomorphic functions $(\widetilde a,0)$, $(\widetilde a',0)$ are
equivalent if they differ by a integer constant. Thus we obtain a complex
of sheaves
$$
0\to\co{}^{\times} \to \widetilde{\cg}^0 @>{\ov D}>>\cg^1 @>{\ov D}>>\cg^2
@>{\ov D}>>\dots
\tag9.3
$$

\medskip

\Theorem{9.1} {\sl The complex (9.3) is a soft resolution of the sheaf 
$\co{}^\times$, and there are natural isomorphisms
$$
H^q(G(X))\ \cong\ H^q(X;\co{}^\times)
\tag9.4$$
for all $q\geq0$.}

\pf  Fix $x\in X$ and consider the cofinal family of contractible
Stein neighborhoods of $x$. For any such neighborhood $U$ we have
$H^q(G^*(U))=0$ for all $q>0$. This results from the long exact sequence in
cohomology associated to the short exact sequence of complexes:
$
0\to {\cd'}^{0,*}(U)\to G^*(U)\to I^{*+1}(U)\to 0.
$ 
Exactness at $q=0$ was proved above. Since the sheaves ${\cd'}^{0,*}$ and
$\ci$ are soft, so are the sheaves $\cg^*$ and  $\widetilde{\cg}^0$. This
easily implies (9.4), which can also be deduced from the quasi-equivalence
of the complex  (9.1) with the trivial complex $\co{}^\times$ (in degree
0).\qed

\vskip.3in

\centerline{ 10. RING STRUCTURE}
\medskip

In this section we show that on any complex manifold $X$, the group
$\hH^{0,*}(X)$  carries the structure of a graded commutative ring with the
property that the homomorphisms $\d_1$ and $\d_2$ in the fundamental
exact sequences (3.4) are ring homomorphisms.  To do this we shall use the
$\dbar$-spark representation of $\hH^{0,*}(X)$ and its analogue in the real
case.

Recall from  [20] (cf. [10] and  [14]) that on any $C^\infty$-manifold $X$
the  {\sl differential characters} of degree $k$ can be defined as the spark
classes $\hH^k(X)$ of the spark complex $(F^*,E^* ,I ^*)$ where 
$$
F ^*\equiv {\cd'}^*(X),  \qquad  E ^*\equiv \ce^*(X),
\qquad I ^*\equiv\In^*(X).
\tag10.1$$

There is a ring structure on
$\hH^*(X)$ defined at the spark level as follows. Fix $\a\in
\hH^k(X)$  and $\b\in \hH^\ell(X)$ and choose representatives 
$(a,r)\in\a$ and $(b,s)\in\b$ with
$$
da\ =\ \phi-r\and
db \ =\ \psi - s
$$
where $a\wedge s$, $ r\wedge b$ and $r\wedge s$ are all well defined
and $r\wedge s \in \In^{k+\ell+2}(X)$. It is shown in [20] 
that this can always be done and that 
$$
\a * \b \ \equdef\ [a\wedge \psi +(-1)^{k+1}r\wedge b] \ =\ 
[a\wedge s +(-1)^{k+1}\phi\wedge b] \ \in\  \hH^{k+\ell+1}(X)
\tag10.2$$
defines a graded commutative ring structure on $\hH^*(X)$ which coincides
with the one defined by Cheeger in [3].

\Note{} Given a spark  a $\in{\cd'}^k(X)$ with  da=$\phi$-r
where $\phi$ is smooth and  r is integrally flat, the elements $\phi$ and
r are uniquely determined by a (cf. [20]). For this reason we refer to the
spark  (a,r)  simply by a.

\Remark{10.1}  In discussing differential characters and sparks the
standard references only consider {\sl real} differential forms and
currents. However, if one replaces these with {\sl complex} differential
forms and currents and if one replaces $\bbr/\bbz$ with $\bbc/\bbz \cong
\bbc^\times$, the standard discussion remains valid. We shall assume that
the forms and currents in (10.1) are  complex. \medskip

Let us denote the $\dbar$-spark complex of \S 7 by $(F ^{0,*},E ^{0,*},
I^*)$ where
$$
F ^{0,*}\equiv {\cd'}^{0,*}(X),  \qquad  E ^{0,*}\equiv \ce^{0,*}(X),
\qquad I ^*\equiv\In^*(X)
$$

\Prop{10.2} {\sl The  projection $\pi:{\cd'}^{*}(X)\to
{\cd'}^{0,*}(X)$ determines a morphism of spark complexes $(F^*,E^*,I^*)\
\to\ (F^{0,*},E^{0,*},I^{*})$ which induces a surjective additive
homomorphism $$
\Pi:\hH^*(X)\ \to \hH^{0,*}(X)
\tag10.3$$
whose kernel is an ideal.}

\pf
Note that for any $a\in {\cd'}^{k}(X)$  one has 
$$[da]^{0,k+1}=\dbar(a^{0,k}).$$ 
From this it is straightforward to see that 
$$
F^*\oplus E^*\oplus I^* \ @>{\pi\oplus\pi\oplus\text{Id}}>>\ F^{0,*}\oplus
E^{0,*}\oplus I^* 
$$
is a map of complexes which commutes  with structure
maps: $I^*\to F^*$, $I^*\to F^{0,*}$, etc.  in the definition of a spark
complex. Consequently, the the induced map  $(a,r) \mapsto
(a^{0,q},r)$ on sparks descends to a well-defined homomorphism
$\Pi:\hH^k(X)\to\hH^{0,k}(X)$ as claimed.

To see that $\Pi$ is surjective, consider a $\dbar$-spark $(A,r) \in
F^{0,k}\oplus I^{k+1}$ with $\dbar A = \Phi-r^{0,k+1}$. By theorems of  
de Rham [7] there exist  $\phi_0\in E^{k+1}$ and 
$a_0\in F^k$ such that $da_0=\phi_0-r$. (That is, there exists a smooth
representative of the cohomology class of $r$ in $H^*(F^*)$.) Let
$A_0=[a_0]^{0,k}$ and $\Phi_0=[\phi_0]^{0,k+1}$ and note that $\dbar
A_0=\Phi_0-r^{0,k+1}$. Hence, $\dbar (A-A_0)= \Phi-\Phi_0$ is a smooth form.
It follows that there exist $b\in F^{0,k-1}$ and $\psi \in E^{0,k}$ with 
$A-A_0=\psi+\dbar b$ (cf. [19, Lemma 1.5]). Set $a=a_0+\psi+db$ and note
that $da=(\phi_0+d\psi)-r$.  Hence, $a$ is a spark of degree $k$ and 
$a^{0,k}=A_0+\psi +\dbar b=A$. Hence, the mapping $\Pi$ is surjective as
claimed.

We now invoke the following.

\Lemma{10.3} {\sl On $\hH^k(X)$, one has that}
$$
\ker\,(\Pi)\ =\ \{\a\in \hH^k(X) : \exists (a,0)\in\a \text { where $a$ is
smooth and } a^{0,k}=0   \} 
$$
\pf
Clearly if $\a=[(a,0)]$ where $a^{0,k}=0$, then $\Pi(\a)=0$.
Conversely, suppose $\a\in\ker\,(\Pi)$ and choose any spark $(a,r)\in \a$.
Then  $\Pi(\a)=0$ means that there exist $b\in{\cd'}^{0,k-1}(X)$ and $s\in
\In^k(X)$ with $r=-ds$ and $a^{0,k}= -\dbar b+s^{0,k} =(-db+s)^{0,k}$.
Replace $(a,r)$ by the equivalent element $(\wt a,0)\equiv (a +db-s,r+ds)$
and note that $\wt a^{0,k}=0$. Repeated application of [19, Lemma 1.5]
now shows that after modification by some $d\wt b$ we may assume that
$\wt a$ is smooth.  Specifically, since $\dbar {\wt a}^{1,k-1}$ is smooth,
there exists ${\wt b}^{1,k-2}$ so that ${\wt a}^{1,k-1} +\dbar{\wt
b}^{1,k-2}$  is smooth. Replacing $\wt a$ by $\wt a + d\wt b$ we can assume
that  ${\wt a}^{1,k-1}$ is smooth and $d\wt a =0$.  Thus $\dbar {\wt
a}^{2,k-2} =-\partial {\wt a}^{1,k-1}$ is smooth, and by [19, Lemma 1.5]
there exists ${\wt b}^{1,k-2}$ with ${\wt a}^{2,k-2}+\dbar {\wt b}^{1,k-2}$
smooth. Continuing inductively completes the proof.
 \qed \medskip

It remains to show that $\ker\,(\Pi)$ is an ideal. Fix $\a\in \ker\,(\Pi)$
and choose $(a,0)\in \a$ with $a^{0,k}=0$ by Lemma 10.3.
Let $\b\in \hH^\ell(X)$ be any class, and $(b,s)\in\b$ any representative
with spark equation $db=\psi-s$.
Then by formula (10.2) the class $\a*\b$ is represented by         
$(a\wedge \psi,0)$ and $(a\wedge \psi)^{0,k+\ell+1} = a^{0,k}\wedge 
\psi^{0,\ell+1} =0$. Hence, $\a*\b\in \ker\,(\Pi)$, and so $\ker\,(\Pi)$ is a
ideal as claimed.\qed

\medskip
The projection $\Pi$ is well behaved with respect to the fundamental
exact sequences.

\Prop {10.4}  {\sl There are commutative diagrams
$$
\CD  
\hH^k(X) @>{\d_1}>>   \cz^{k+1}_{\bbz}(X)  \\
@V{\Pi}VV  @VV{\pi}V     \\
\hH^{0,k}(X) @>{\d_1}>>   \cz^{0,k+1}_{\bbz}(X)
\endCD
\qquad\ \qquad
\CD  
\hH^k(X) @>{\d_2}>>   H^{k+1}(X;\,\bbz)  \\
@V{\Pi}VV    @VV{\text{\rm =}}V      \\
\hH^{0,k}(X) @>{\d_2}>>    H^{k+1}(X;\,\bbz)
\endCD
$$
where $\pi$ denotes the projection.  In fact the surjective homomorphism 
$\Pi$ expands to a morphism of the fundamental diagram for $\hH^k(X)$
([19]): 
$$
\CD
\ @.  0   @.  0  @.  0 @. \  \\
@. @VVV @VVV @VVV @.  \\
0 @>>> \frac {H^k(X,\bbc)}{H^k_{\text{\rm free}}(X;\bbz)}  @>>>
\hH^{k}_\infty(X)  @>>>  d\ce^{k}(X)    @>>> 0 \\ 
@. @VVV @VVV @VVV @.  \\ 
0 @>>> H^k(X, \bbc^{\times})  @>>> \hH^{k}(X) 
@>{\d_1}>>  \cz^{k+1}_{\bbz}(X)    @>>> 0 \\ @. @VVV @V{\d_2}VV @VVV @.  \\
0 @>>>  H^k_{\text{\rm tor}}(X, \bbz)  @>>> H^{k+1}(X;\bbz) @>>>
H^{k+1}_{\text{\rm free}}(X)   @>>> 0 \\ 
@. @VVV @VVV @VVV @.  \\
\ @.  0   @.  0  @.  0 @. \  
\endCD
\tag10.4
$$
to the diagram (3.4).}

\pf
Fix $\a\in\hH^{0,k}(X)$ and  $\a'\in\hH^k(X)$ with $\Pi(\a')=\a$.
Choose a representative $(a,r)\in \a'$ with $da=\phi -r$ where
$\phi\in \cz^{k+1}(X)$. Then $(\pi(a),r) = (a^{0,k},r)$ represents
$\a$ and we see that $\dbar a^{0,k} = [da]^{k+1} =\phi^{k+1}-r^{k+1}$.
Thus $\d_1\a=\phi^{k+1}=\pi(\phi)=\pi(\d_1\a')$. The argument for $\d_2$ is
similar. The assertion concerning the diagrams follows directly.
\qed

This brings us to the  main result.

\Theorem{10.5} {\sl  For any complex manifold $X$ the $\dbar$-spark classes
$\hH^{0,*}(X)$ carry a biadditive product 
$$
*:\hH^{0,k}(X) \times \hH^{0,\ell}(X)\ \arr\ \hH^{0,k+\ell+1}(X)
$$
which makes $\hH^{0,*}(X)$ into a graded commutative ring (commutative in
the sense that $$\a*\b=(-1)^{k+\ell+1}\b*\a$$ for $\a\in\hH^{0,k}(X)$
and $\b \in \hH^{0,\ell}(X)$.)  With respect to this ring structure
 the fundamental maps $\d_1$ and $\d_2$ are ring homomorphisms.
}

\pf
The first assertions follow immediately from  Proposition 10.2.
For the last assertion we apply Proposition 10.4
as follows.  Suppose $\a\in\hH^{0,k}(X)$ and $\b \in \hH^{0,\ell}(X)$ are
given and choose $\a'\in\hH^{k}(X)$ and $\b' \in \hH^{\ell}(X)$ with
$\Pi(\a')=\a$ and $\Pi(\b')=\b$. If $\d_1\a' = \phi$ and $\d_1\b' = \psi$,
then $\d_1(\a'*\b')=\phi\wedge\psi$ and so $\d_1(\a*\b)=
\d_1\Pi(\a'*\b')=\pi \d_1(\a'*\b')= \pi(\phi\wedge\psi)=
\pi(\phi)\wedge\pi(\psi)= \d_1(\a)\wedge \d_1(\b)$. 
The argument for $\d_2$ is similar.
 \qed
\medskip

\Remark{10.5}  The product on $\hH^{0,*}(X)$ can be defined at the spark
level as follows.  Let $(a,r)\in\a\in\hH^{0,k}(X)$ and 
$(b,s)\in\b\in\hH^{0,\ell}(X)$ be sparks with $\dbar a=\phi-r$ and
$\dbar b=\psi-s$. Assume that the currents $a\wedge s^{0,\ell+1}$, $r^{0,k+1}
\wedge b$ and $r\wedge s$ are well  defined and that $r\wedge s$ is
integrally flat. (The intersection theory for currents in [20]  shows
that such representatives exist,)  Then the product $\a*\b$ is represented by
each of the following: 
$$
(a\wedge \psi +(-1)^{k+1}r^{0,k+1}\wedge b,\ r\wedge s)
\and
(a\wedge s^{0,\ell+1} +(-1)^{k+1}\phi\wedge b,\ r\wedge s).
$$

\Cor{10.6} {The subgroup $H^*(X,\co{}^\times)=\ker \d_1$ is an ideal in
$\hH^{0,*}(X)$.  In particular, $H^*(X,\co{}^\times)$ has a  ring
structure compatible with the spark product. In this product the coboundary
map
$$
\d:H^*(X,\co{}^\times)\ \arr\ H^*(X;\,\bbz)
$$
from the sheaf sequence $0\to\bbz\to\co{}\to\co{}^\times\to 0$
is a ring homomorphism.}

\vskip.3in

\centerline{ 11.  FUNCTORIALITY}
\medskip

The main result of this section is the following.

\Theorem{11.1}  {\sl 
Any holomorphic map $f:Y\to X$ between
complex manifolds  induces a graded ring homomorphism
$$
f^*:\hH^{0,*}(X) \to \hH^{0,*}(Y) 
$$
with the property that  if $g:Z\to Y$ is holomorphic,
then $(f\circ g)^*=g^*\circ f^*$. In other words, $\hH^{0,*}(\bullet)$
is a graded ring functor on the category of complex manifolds and
holomorphic maps. }

\pf We know that $f$ induces a ring homomorphism $f^*:\hH^*(X)\to\hH^*(Y)$
with the asserted property.  We need only show that $f^* \ker\,(\Pi)
\subset \ker\,(\Pi)$. This follows directly from the Lemma 10.3.\qed

\vskip.3in

\centerline{ 12.  REFINED CHERN CLASSES FOR HOLOMORPHIC BUNDLES}
\medskip

In this section we construct refined Chern classes ${\wh d}_k(E)\in
H^{2k-1}(X,\co{}^\times)$ for holomorphic bundles $E$ which extend the
tautological case $k=1$. These classes possess the usual properties and map
to the integral Chern classes under the coboundary map 
$\d:H^{2k-1}(X,\co{}^\times)\to  H^{2k}(X,\bbz)$. In fact they descend 
to holomorphic $K$-theory in the sense of Grothendieck.

These classes can also be accessed through Deligne cohomology,
and our construction could be considered {\sl a Chern-Weil approach 
to Deligne cohomology} in lowest degree. In [13]  this is extended to the full Deligne theory.

Our point of departure is the fundamental work of Cheeger and Simons [4]
who showed that for a smooth complex vector bundle  $E\to X$ with unitary
connection $\nabla$ there exist refined Chern classes ${\wh c}_k(E,\nabla)
\in \hH^{2k-1}(X)$ with
$$
\d_1({\wh c}_k(E,\nabla))\ =\ c_k(\Omega^{\nabla})
\and
\d_2({\wh c}_k(E,\nabla))\ =\ c_k(E)
\tag12.1$$
where $c_k(E)$ is the integral $k\th$ Chern class and $c_k(\Omega^{\nabla})$
is the Chern-Weil form representing $c_k(E)\otimes \bbr$ in the curvature of
$\nabla$ (cf. (10.4)). Setting $\wh c(E) =  1+{\wh c}_1+{\wh c}_2 
+\dots  $ they show
$$
{\wh c}(E\oplus E',\nabla\oplus \nabla' )\ =\ 
{\wh c}(E,\nabla)*{\wh c}( E', \nabla' )
\tag12.2$$

Whenever $X$ is a complex manifold we can take the projections
$$
{\wh d}_k(E,\nabla)\ \equiv\ \Pi \{{\wh c}_k(E,\nabla)\} \in \hH^{0,2k-1}(X)
$$
and note that by (12.1), (12.2) and Proposition 10.4
$$
\d_1({\wh d}_k(E,\nabla))\ =\ c_k(\Omega^{\nabla})^{0,2k}
\and
\d_2({\wh d}_k(E,\nabla))\ =\ c_k(E)
\tag12.3$$
and with
${\wh d}(E,\nabla)\equiv 1+{\wh d}_1(E,\nabla)+{\wh d}_2(E,\nabla)+\dots$,
$$
{\wh d}(E\oplus E',\nabla\oplus \nabla' )\ =\ 
{\wh d}(E,\nabla)*{\wh d}( E', \nabla' ).
\tag12.4$$

Suppose now that $E$ is holomorphic and is provided with a hermitian
metric $h$.  Let $\nabla$ be the associated canonical hermitian connection.
Then ${\wh c}_k(E,\nabla)$ is of type $(k,k)$ and so by (12.3) we have 
$$
{\wh d}_k(E,\nabla) \ \in\ \ker(\d_1)\ =\ H^{2k-1}(X,\co{}^{\times})
\tag12.5$$

\Prop{12.1} {\sl The class  in
(12.5) is independent of the choice of hermitian metric.}

\pf
Let $h_0,h_1$ be hermitian metrics on $E$ with canonical connections
$\nabla^0,\nabla^1$ respectively. Then (see [4])
$$
{\wh c}_k(E,\nabla^1) - {\wh c}_k(E,\nabla^0)\ =\ 
[T]
$$
where $[T]$ is the differential character represented by the smooth
transgression form
$$
T\ =\ T(\nabla^1,\nabla^0)\ \equiv\  k \int_0^1
C_k(\nabla^1-\nabla^0, \Omega_t,...,\Omega_t)\,dt
$$
where $C_k(X_1,...,X_k) $ is the polarization of the $k\th$ elementary
symmetric function and where $\Omega_t$ is the curvature of the connection
$\nabla^t \equiv t\nabla^1-(1-t)\nabla^0$. Fix a local
holomorphic frame field for $E$ and let $H_j$ be the hermitian matrix
representing the metric $h_j$ with respect to this trivialization.
Then
$$
\nabla^1-\nabla^0 = \omega_1-\omega_0\qquad\qquad\text{where} 
\qquad \omega_j \equiv H_j^{-1}\partial H_j = \partial \log H_j.
$$
In this framing,   $\nabla^t=d + \omega_t$ where
$\omega_t=t\omega_1-(1-t)\omega_0$ and so its curvature
$\Omega_t=d\omega_t-\omega_t\wedge \omega_t$ only has Hodge components of 
type $(1,1)$ and $(2,0)$. It follows that the Hodge components  
$$
T^{p,q} \ =\ 0  \qquad\text{ for }\ \  p<q.
\tag12.6$$
In particular, $T^{0,2k-1}=0$ and so $\Pi[T]=
{\wh d}_k(E,\nabla^1) - {\wh d}_k(E,\nabla^0)=0$. \qed

\medskip

By Proposition 12.1 each holomorphic vector bundle $E$ of rank $k$
has a well defined total refined Chern class
$$
{\wh d}(E)\ =\ 1+{\wh d}_1(E)+\dots+{\wh d}_k(E)\ \in \
H^{\text{*}}(X,\co{}^{\times})
$$
Denote by $\cv^k(X)$ the set of isomorphism classes of holomorphic vector
bundles of rank $k$ on $X$, and note that
$\cv(X)=\coprod_{k\geq0}\cv^k(X)$ is an additive monoid under Whitney sum.

\Theorem{12.3} {\sl On any complex manifold there is a natural
transformation of functors
$$
{\wh d} : \cv(X) \ \arr\ H^*(X,\co{}^{\times})
$$
with the property that: \medskip

(i)\ \ \ \ ${\wh d}(E\oplus F)\ =\ {\wh d}(E)*{\wh d}(F)$,  \medskip

(ii)\ \ \ ${\wh d}:\cv^1(X)\ @>{\approx}>> \ 1+H^1(X,\co{}^{\times})$
is an isomorphism, and \medskip

(iii)\ \ under the coboundary map $\d:H^\ell(X,\co{}^{\times})\to
H^{\ell+1}(X,\,\bbz)$,
$$
\d\circ {\wh d}\ =\ c \ \ \text{(the total integral Chern class)}.
$$
}

\pf Property (i) follows from (12.4), Property (ii) is classical, and
Property (iii) follows from (12.3). \qed

\medskip

Property (i) implies that ${\wh d}$ extends to the additive group
completion 
$$
\cv^+(X)\equiv \cv(X)\times\cv(X)/\sim
\tag12.7
$$ 
where 
$(E,F)\sim (E',F')$ iff there exists $G\in \cv(X)$ with $E\oplus F'\oplus G
\cong E'\oplus F\oplus G$. It is natural to ask is whether ${\wh d}$
then descends to the Grothendieck quotient. It does.

\Theorem{12.4}  {\sl For any short exact sequence of holomorphic vector
bundles on $X$
$$
0\ \arr\ E'\ \arr\ E\ \arr\ E''\ \arr\ 0
$$
one has }
$$
{\wh d}(E) \ =\ {\wh d}(E')*{\wh d}(E'')
$$

\pf
Choose hermitian metrics $h'$ and $h''$ for $E'$ and $E''$ respectively,
and  choose a $C^\infty$-splitting
$$
0\ \arr\ E'\ @>{i}>>\ E\ \underset {@<<{\s}<}\to {@>{\pi}>>}
\ E''\ \arr\ 0.
$$
Define a hermitian metric $h=h'\oplus h''$ on $E$ via the smooth isomorphism
$$
(i,\s):E'\oplus E'' \ \arr\ E.
$$

We must show that the $\dbar$-spark classes for the direct sum of the
canonical hermitian connections on $E'\oplus E''$ agree with the spark
classes for the canonical hermitian connection on $E$. For this it will
suffice to work over an open set $U\subset X$ on which there exist
holomorphic framings
$$
(e_1',...,e_n') \ \  \text{ for }\ \  E'
\and
(e_1'',...,e_m'') \ \  \text{ for }\ \  E''
$$
We now choose two framings form $E$ over $U$.

\noindent
{\bf Framing 1 (holomorphic):} \ \ 
$$
(e_1,...,e_{n+m})\ =\ (e_1',...,e_n',{\wh e}_1'',..., {\wh e}_m'')
$$
where each ${\wh e}_k''$ is a holomorphic lift of ${ e}_k''$.

\noindent
{\bf Framing 2 (smooth and direct-sum compatible):} \ \ 
$$
({\wt e}_1,...,{\wt e}_{n+m})\ =\ (e_1',...,e_n',\s  e_1'',..., \s e_m'').
$$

Let $H'$ and $H''$ be the smooth hermitian-matrix-valued functions
representing the metric $h'$ and $h''$ in their respective holomorphic
frames over $U$.  Then the connection 1-form for the direct sum connection
on $E$ in framing 2 is 
$$
\omega \ =\ \omega'\oplus\omega'' \qquad\ \ \text{where}\ \ \ 
\omega'\ =\ \partial H'\cdot (H')^{-1}\ \ \text{and}\ \ 
  \omega''\ =\ \partial H''\cdot (H'')^{-1}.
$$
In particular in framing 2 this connection form is of type (1,0).

Now let $H$ be the smooth hermitian-matrix-valued function
representing the metric $h$ in the holomorphic framing 1 of $E$ over $U$.
In this framing the connection 1-form of the canonical hermitian connection
on $E$ is 
$$
\theta\ =\ \partial H\cdot H^{-1}.
$$
We want to compute the connection 1-form ${\wt \theta}$ for this connection
in  the second framing $({\wt e}_1,...,{\wt e}_{n+m})$.   For this we
consider the change of framing
${\wt e}_k = \sum_{\ell=1}^{n+m} g_{k\ell} e_{\ell}$ and recall 
(cf. [11, p. 72]) that
$$
{\wt \theta}\ =\ dg\cdot g^{-1} + g\cdot \theta\cdot  g^{-1}
\tag12.8$$
where $g=(g_{k\ell})$.  We now observe that $g_{k\ell}=\d_{k\ell}$ for
$1\leq k\leq n$ and so $g$ has the form
$$
g\ =\ \left(\matrix I&0\\A&I
\endmatrix\right)
\and 
g^{-1}\ =\ \left(\matrix I&0\\-A&I
\endmatrix\right).
$$
In particular we have that
$$
dg\cdot g^{-1}\ =\ dg\cdot g\ =\ g\cdot dg\ =\ 
dg\ =\ \left(\matrix 0&0\\dA&0
\endmatrix\right)
\tag12.9$$

Suppose now that $\Phi$ is any Ad-invariant symmetric $k$-multilinear
function on the Lie algebra of $GL_{n+m}(\bbc)$. Then the two given
connections on $E$ give rise to two Cheeger-Simons differential characters,
say $\wh\Phi_0$ and $\wh\Phi_1$, and the difference  
$$
\wh\Phi_1 - \wh\Phi_0\ =\ [T]
$$
where $[T]$ is the character associated to the smooth differential form
$$
T\ =\ k \int_0^1 \Phi(\wt \theta -\omega, \Omega_t,...,\Omega_t)\,dt.
$$
where $\Omega_t$ is the curvature 2-form of the connection $\omega_t
\equiv (1-t)\wt\theta +t\omega$.
 To compute the corresponding difference in $\dbar$-characters we take
the projection onto $(0,2k-1)$-forms
$$\aligned
T^{0,2k-1}\ &=\ k \int_0^1 \Phi(\wt \theta -\omega,
\Omega_t,...,\Omega_t)^{0,2k-1}\,dt  \\
&=  k \int_0^1 \Phi(dg\cdot g^{-1},
\Omega_t,...,\Omega_t)^{0,2k-1}\,dt \\
&=  k \int_0^1 \Phi(dg,
\Omega_t,...,\Omega_t)^{0,2k-1}\,dt
\endaligned\tag12.10$$
where the second line follows from the fact that $\omega$ and $\theta$
are of type (1,0) and the last line follows from (12.9).

Observe now that
$$
\Phi(dg, \Omega_t,...,\Omega_t)^{0,2k-1}\ =\ 
\Phi(\dbar g, \Omega_t^{0,2},...,\Omega_t^{0,2})
$$
and from type considerations one computes that
$$
\Omega_t^{0,2} =  \left\{d\omega_t-\omega_t\wedge\omega_t\right\}^{0,2} = 
(1-t)^2\dbar g\wedge\dbar g=(1-t)^2 [\dbar g,\dbar g].
$$
From the Adjoint-invariance of $\Phi$ we conclude that 
$$
\Phi(\dbar g, [\dbar g,\dbar g],...,[\dbar g,\dbar g])\ =\ 0.
$$
 Indeed for any matrix-valued 1-form $w$, invariance implies that
$$\aligned
\Phi([w,w],w,[w,w],...,&[w,w])+\Phi(w,[w,w],...,[w,w])\\
&\ +\Phi(w,w,[w,[w,w]],[w,w]...,[w,w]) +\dots\\\
&\qquad\qquad =\  2\Phi(w,[w,w],...,[w,w])\ =\ 0
\endaligned$$
since $[w,[w,w]]=0$ by the Jacobi Identity. We conclude that
$T^{0,2k-1}=0$ and the proof is complete.\qed

\medskip

Consider the natural transformation
$$
{\wh d}: \cv(X)^+\ \arr H^*(X,\co{}^{\times})
$$
where $\cv(X)^+$ is the group completion of $\cv(X)$ (cf. (12.7)).
Following Grothendieck we define the holomorphic $K$-theory of $X$
to be the quotient
$$
K_{\text{hol}}(X)  \ \equiv\ \cv(X)^+/\sim
$$
where $\sim$ is the equivalence relation generated by setting 
$[E]\sim [E'\oplus E'']$ whenever there is a short exact sequence
of holomorphic bundles $0\to E'\to E\to E''\to 0$ on $X$.

\Cor{12.5} {\sl  The natural transformation\  \  $\wh d$ defined above
descends to a natural transformation 
$$
\wh d : K_{\text{hol}}(X)\ \arr\  H^*(X,\co{}^{\times})
$$
such that properties (i), (ii) and (iii) of Theorem 12.3 continue to hold.
In particular for algebraic manifolds this gives a total Chern class map
$$
\wh d:CH(X)\ \arr\  H^*(X,\co{}^{\times})
$$
from the group of algebraic cycles modulo rational equivalence.}

\vskip.3in

\centerline{ 13.  BOTT VANISHING FOR HOLOMORPHIC FOLIATIONS}
\medskip

In 1969 R. Bott constructed a family of connections on the normal bundle
of any smooth foliation of a manifold. Using these classes he established
the vanishing of characteristic classes of the normal bundle in
sufficiently high degrees [1].  Cheeger and Simons then showed that 
for these classes which vanish, the corresponding differential characters
are well defined (independent of the choice of Bott connection) and
represent secondary invariants of the foliation [4, \S 7].  These
invariants are highly non-trivial and can vary continuously as
$\bbr^k$-valued objects.

Suppose now that $N$ is the normal bundle to a holomorphic foliation of
codimension-$q$ on a complex manifold $X$.  Then there are two natural
families of connections to consider on $N$: the family of Bott connections
and the family of canonical hermitian connections.

\Prop{13.1} {\sl Let  $P(c_1,...,c_q)$ be a polynomial in  Chern classes
which is of pure cohomology degree $2k$ with $k>2q$.  Then the
$\dbar$-character $P({\wh c}_1,...,{\wh c}_q)^{0,2k-1}$ for the Bott
connections agrees with the $\dbar$-character $P({\wh d}_1,...,{\wh
d}_q)$ for the  canonical hermitian connections.
}

\Note{} The polynomial $P$ has {\sl pure cohomology degree} $2k$ 
if it satisfies the weighted homogeneity condition: \ $P(tc_1, t^2
c_2,...,t^qc_q)=t^kP(c_1,...,c_q)$ for all  $t\in\bbr$.

\Theorem{13.2}  {\sl Let $N$ be a holomorphic bundle of rank $q$ on a
complex manifold $X$.  If $N$ is (isomorphic to) the normal bundle of a
holomorphic foliation of $X$, then for every  polynomial $P$ of pure
cohomology degree $k>2q$, the associated refined Chern class satisfies
$$
P({\wh d}_1(N),...,{\wh d}_q(N)) \ \in\ H^{2k-1}(X;\,\bbc^{\times}) \subset
H^{2k-1}(X;\,{\co{}}^{\times})
$$
 } 

\pf
This is an immediate consequence of Proposition 13.1 which we shall now
prove.

Suppose $N$ is the normal bundle to a holomorphic foliation $\cf$ of
codimension-$q$. Then $X$ has a distinguished atlas of coordinate charts
$(z_\a,w_\a):U_\a\to\bbc^{n-q}\times\bbc^q$  such that in the open subset
$U_\a\subset X$, $\cf$ is defined by the equation $w_\a=$ constant. 
Under the change of such coordinates one has
$$
\frac{\partial  w_\a}{\partial z_\b}\  =\
0,
\qquad\text{ that is }\qquad w_\a\ =\ w_\a(w_\b)
\tag13.1$$
depends on $w_\b$ alone. Note that therefore the operator
$$
\partial_w\ \equiv\ \sum d w^j \wedge \frac{\partial}{\partial w^j}
$$
is independent of the choice of distinguished coordinates $(z,w)$ and
therefore globally defined on $X$.

 Suppose now that a hermitian metric $h$ is given for the normal
bundle $N=\text{span}\{dw\}$ and let $H_\a$ be the hermitian matrix
representing $h$ in the holomorphic frame $dw^1_{\a},...,dw^q_{\a}$ 
 We define a connection 1-form $\theta_\a$
for $N$ in this frame by setting 
$$
\theta_\a\ \equiv\ \partial_w H_\a\cdot H_{\a}^{-1}
$$
The transition functions for the  $N=$ are
given by the jacobian matrix 
$$
g\ =\ g_{\a,\b}\ =\ \frac{\partial (w_\a)}{\partial (w_\b)}.
$$
A straightforward calculation shows that
$$
\theta_\a \ =\ g\cdot \theta_\b\cdot g^{-1} +\partial_w g\cdot g^{-1}\ 
=\ g\cdot \theta_\b\cdot g^{-1} +dg\cdot g^{-1},
$$
and so these 1-forms assemble to give a well-defined Bott connection
$\nabla$ on $N$.

Recall that in the local frame $dw^1_{\a},...,dw^q_{\a}$ for $N$ the
canonical hermitian connection $\wt \nabla$is given by the 1-form
$$
{\wt \theta}_\a\  \equiv\ \partial H_\a\cdot H_{\a}^{-1}.
$$

We can now explicitly compute the transgression term in any 
distinguished coordinate system $(z,w)$. (We shall drop
the $\a$'s for convenience.)  Given $P$, let $\Phi(X_1,...,X_k)$ be the
Ad-invariant $k$-multilinear symmetric function on the Lie algebra 
$\gl_q(\bbc)$ such that $P(\s_1(X),...,\s_q(X))=\Phi(X,...,X)$ where 
$\s_j(X)$ is the $j\th$ elementary symmetric function of the eigenvalues
of $X$.  Then the difference between the $\dbar$-differential characters
associated to $P$ for the two connections $\wt\nabla$ and $\nabla$ is the
character associated to the smooth form 
$$
T^{0,2k-1}\ =\ k\int_0^1\Phi(\wt\theta-\theta,
\Omega_t,...,\Omega_t)^{0,2k-1} \, dt.
$$
 where 
$$
\wt\theta-\theta\ =\ \partial_z H\cdot H^{-1}
$$
and 
$
\Omega_t\ =\ d\theta_t - \theta_t\wedge\theta_t
$ with
$$
\theta_t\ \equiv\ (1-t)\wt\theta + t\theta
\ =\ \partial_w H\cdot H^{-1}-t(\partial_z H\cdot H^{-1}).
$$
Since $\theta_t$ is of type $1,0$ we have that $\Omega_t$ is of type
$1,1$ plus $2,0$. Hence, $T^{0,2k-1}=0$.\qed

\medskip
\Note{13.3} The above calculation shows that 
$$
T^{p,q}\ =\ k\int_0^1\Phi(\wt\theta-\theta,
\Omega_t,...,\Omega_t)^{p,q} \, dt\ =\ 0\qquad\text{ for all } p<q.
$$
This will allow us to generalize Theorem 13.2 to all Deligne cohomology.

\vskip.3in
\vfill\eject

\centerline{ 14.  GENERALIZATIONS AND THE RELATION TO DELIGNE COHOMOLOGY}
\medskip

Most of the discussion above can be easily generalized to other truncations
of the de Rham complex.  This leads to spark complexes with multiplicative
structure on the associated spark classes. The role of 
$H^*(X,\co{}^{\times})$ in the above is now played by more general Deligne
cohomology groups.

Fix an integer $p > 0$ and consider the truncated de Rham complex
$(\Dp *,\ddp)$ with
$$
\Dp k\ \equiv 
\underset {r< p}\to {\bigoplus_{r+s=k}}{\cd'}^{r,s}(X)
\and
\ddp\equiv \Psi\circ d
$$
where 
$$
\Psi:{\cd'}^k(X)\ \arr\ \Dp k
$$
is the projection $\Psi(a) = a^{0,k}+ a^{1,k-1}+\dots + a^{p-1,k-p+1}$.
 Note the subcomplex
$$
\Ep k\ \equiv 
\underset {r < p}\to {\bigoplus_{r+s=k}}{\ce}^{r,s}(X)
$$
of smooth forms with projection $\Psi:\ce^*(X)\to \Ep *$ whose
kernel is a {\bf $d$-closed ideal.}

\Def{14.1}  By the {\bf $\ddp$-spark complex of
level $p$} we mean the triple $(F^*,E^*, I^*)$ where
$$\aligned
F^k\ &\equiv\ \Dp k \\
E^k\ &\equiv\ \Ep k \\
I^k\ &\equiv\  \In^k(X)
\endaligned$$
with maps
$$
E^*\ \subset F^*\and \Psi:I^*\ \arr\ F^*
$$
given by the inclusion and projection above. The group of associated spark
classes in degree $k$ will  be denoted by $\Hp k$.
\medskip

Note that a spark in this complex is a pair 
$(a,r)\in \Dp k \times \In^{k+1}(X)$ such that $dr=0$,
 and $a$ satisfies the $\ddp$-{\bf spark equation}
$$
\ddp a \ =\ \phi - \Psi(r)
\tag14.1$$
for some  smooth $(k+1)$-form $\phi$ on $X$. Note that
$\phi^{\ell,k+1-\ell}=0$ for all  $\ell\geq p$ and that
$\ddp \phi = 0$, i.e., $d\phi \equiv 0 \ \ (\text{mod} \ker \Psi)$.

The $\dbar$-spark complexes of level $1$ are exactly the $\dbar$-sparks
discussed in \S 7. 

Of course we have not yet established that the triple  $(F^*,E^*,I^*)$ in 
Definition 14.1 is a spark complex. The fact that $E^k\cap I^k=\{0\}$ for
$k>0$ is proved in Proposition B.1.  To show that the inclusion
$E^*\subset F^*$ induces an isomorphism in cohomology 
$$
H^*(E)\ \cong\ H^*(F)\ \equiv\ \Hop * {},
$$
consider $E^*$ and $F^*$ as  double complexes and note
that the inclusion induces an isomorphism on vertical 
$\dbar$-cohomology (and hence in total cohomology by a standard spectral
sequence argument  (cf. [2, Lemma 1.2.5]). Note that $\Hop * {}$ is the 
hypercohomology of the  complex of sheaves: 
$0\to \Omega^0_X\to\Omega^1_X\to\Omega^2_X\to\dots\to\Omega^{p-1}_X\to 0$.

To analyze the groups of spark classes $\Hp k$ we   recall the
following (cf. [2] or [24]).

\Def {14.2} By the {\bf Deligne complex of level $p$} we mean the
complex of sheaves
$$
\Dcx p\,: \qquad\qquad 0\to \bbz \to \Omega^0_X \to
\Omega^1_X\to\Omega^2_X\to\dots\to\Omega^{p-1}_X\to0 \qquad\qquad\qquad
$$
where $\Omega^k_X$, the sheaf of holomorphic $k$-forms on $X$, is
considered to live in degree $k+1$.
(Typically the constant sheaf $\bbz$ is embedded into $\co X = \Omega^0_X$
by  $m\mapsto (2\pi i)^pm$. We will not adopt this convention here.)
By the {\bf Deligne cohomology } of $X$ in level $p$  we mean the
hypercohomology of this complex:  
$$
\HD * p \ \equiv H^*(X, \Dcx p)
$$ 
\medskip

Note that when $p=1$ the complex of sheaves $\Dcx 1$ is
quasi-equivalent to ${\Cal O}_X^{\times}[-1]$, the sheaf              
${\Cal O}_X^{\times}$  shifted to the right by 1 so its formal degree is 1
and not 0. Thus we have that 
$$
\HD k 1\ \cong\ H^{k-1}(X,\co{X}^\times),
$$
and so in the fundamental exact grid (3.4), the left-middle term could
be replaced by $\HD k 1$.  In this form the picture generalizes to
all levels.

\Prop{14.3} {\sl The first fundamental exact sequence (1.5) for the 
group $\Hp k$ is
$$
0\ \arr\ \HD {k+1} p\ \arr\ \Hp k\ @>{\d_1}>>\ {\Zp {k+1}} \ \arr\ 0
$$
where ${\Zp {k+1}}$ is the set of $\ddp$-closed forms in $\Ep {k+1}$ 
which represent classes in 
$
\Hop {k+1}{\bbz} \equiv \Image\{\Psi_*:H^{k+1}(X;\,\bbz) \arr \Hop
{k+1}{}\}. 
$ }

\pf The identification of Image$(\d_1)$ is straightforward. To identify the
kernel consider the acyclic resolution
$$
\CD
\In^*_X @>{\Psi}>>  {\cd'}^{*,*}_X  \\
@AAA  @AAA  \\
\bbz @>{\psi}>>  \Omega_X^*
\endCD
$$
of the Deligne complex.   Proposition A.3 then gives the following.

\Prop{14.4} {\sl There is a natural isomorphism}
$$
\HD {*} p \ \cong\ H^*(\text{Cone}\,(\In^*(X) @>{\Psi}>> \Dp *))
$$
This cone complex is exactly the cone complex $G^*$ associated 
to our spark complex as in (1.5). \qed
\medskip

 The diagram (1.6) for the groups $\Hp {k-1}$ can be written as
$$
\CD
\ @.  0   @.  0  @.  0 @. \  \\
@. @VVV @VVV @VVV @.  \\
0 @>>> \frac {\Hop {k-1} {}}{\Hop {k-1} {\bbz}}  @>>> 
\hH^{k-1}_{\infty}(X,p) @>>> \ddp\Ep {k-1}    @>>> 0 \\ @. @VVV @VVV @VVV
@.  \\ 0 @>>> \HD {k} p @>>> \Hp {k-1}  @>{\d_1}>> 
{\Zp {k}}   @>>> 0 \\ @. @VVV @V{\d_2}VV @VVV @.  \\
0 @>>> \ker(\Psi_*)  @>>> H^{k}(X;\bbz) @>{\Psi_*}>>
\Hop {k} \bbz   @>>> 0 \\ @. @VVV @VVV @VVV @.  \\
\ @.  0   @.  0  @.  0 @. \  
\endCD
$$
where as usual $ \hH^{k-1}_{\infty}(X,p)$ denotes the spark classes  representable by
smooth forms.

Note that when $X$ is Kaehler,  we have
$$\aligned
\ker(\Psi_*) \ &=\ H^{k}(X;\bbz)\cap \bigoplus_{j\geq p} H^{j,k-j}(X) \\
&=\ H^{k}(X;\bbz)\cap \underset{r+s=k}\to{\bigoplus_{|r-s|\leq k-2p}}
H^{r,s}(X) \endaligned
$$
where the second line is deduced from the reality of $H^k(X;\,\bbz)$.
In particular when $k=2p$ we deduce
$$
\CD
\ @.  0   @.  0  @.  0 @. \  \\
@. @VVV @VVV @VVV @.  \\
0 @>>> \cj_p(X)  @>>>  \hH^{2p-1}_{\infty}(X,p)
@>>> \ddp\Ep {2p-1}    @>>> 0 \\ @. @VVV @VVV @VVV @.  \\
0 @>>> \HD {2p} p @>>> \Hp {2p-1}  @>{\d_1}>> 
{\Zp {2p}}   @>>> 0 \\ @. @VVV @V{\d_2}VV @VVV @.  \\
0 @>>> \text{Hdg}^{p,p}(X)  @>>> H^{2p}(X;\bbz) @>{\Psi_*}>>
\Hop {2p} \bbz   @>>> 0 \\ @. @VVV @VVV @VVV @.  \\
\ @.  0   @.  0  @.  0 @. \  
\endCD
$$
where $\cj_p(X)$ denotes the Griffiths' $p\th$ intermediate Jacobian and 
$\text{Hdg}^{p,p}(X)\subset H^{2p}(X;\bbz)$ are the {\sl Hodge classes},
i.e., the integral classes representable over $\bbr$ by   closed
$(p,p)$-forms.  The left vertical sequence, which is classical for Deligne
cohomology, is deduced directly from this theory.
 
\medskip
\noindent
{\bf Remark 14.5. (The Deligne class of a holomorphic chain.)} A {\sl
holomorphic k-chain} on $X$ is an integral current which
can be written as a  locally finite sum 
$
Z=\sum_j \,n_jV_j
$
where for each $j$, $n_j\in \bbz$ and $V_j$ is an irreducible complex
 subvariety of dimension-$k$. Any such chain determines a spark
$(0,Z)\in \cs^{2p-1}(X)$ where $p=n-k$  (since $d0=0-\Psi(Z)=0$) and therefore
a class in $\hH^{2p-1}(X)$.   This is clearly in the
kernel of $\d_1$ and we obtain a class
$$
[(0,Z)] \in \HD {2p} p.
$$
This retrieves Griffiths' generalized Abel-Jacobi mapping when $X$ is
compact Kaehler and $Z$ is homologous to zero. The existence of this class
goes back to Deligne [6].

\medskip
\noindent
{\bf Remark 14.6. (The Deligne class of a maximally complex cycle.)}
The above construction clearly generalizes to any integral cycle
of restricted Hodge type. An interesting case is the following.
An integral cycle $M\in\ci^{2p+1}(X)$ with Dolbeault decomposition 
$M=M^{p,p+1}+M^{p+1,p}$ is called {\sl maximally complex} (cf. [15]).
As above, any such a cycle determines a class in Deligne cohomology
$$
[(0,M)] \in \HD {2p+1} p.
$$
As in [15] we shall say that $M$ is the {\sl boundary of a holomorphic
chain} if $M=dZ$ where $Z$ is a current defined by a holomorphic 
$(n-p)$-chain in $X-M$.

\Prop{14.6} {\sl If $M$ is a maximally complex cycle which bounds a 
holomorphic chain, then its Deligne class in $\HD {2p+1} p$ is 
zero.}

\pf If $M=dZ$ where $Z\in \ci^{2p}(X)\cap{\cd'}^{p,p}(X)$, then
$[(0,M)]= [(\Psi(Z),M-dZ)]=[(0,0)]=0$ in $\HD {2p+1} p$.\qed

\medskip
The above remarks also apply to cycles with compact support with classes
in {\sl  compactly supported Deligne cohomology}.  In this case one
retrieves the moment conditions characterizing boundaries of holomorphic
chains in [15].

\vskip .3in

\centerline{ 15.\ \  THEOREMS IN THE GENERAL SETTING} \medskip

Most of the results proved for $\dbar$-sparks carry over to 
the general $\ddp$-sparks of \S 14.

\Theorem {15.1}  {\sl  The projection $\Psi$ induces a morphism from the de
Rham-Federer spark complex to the $\ddp$-spark complex in 14.1. This
induces a surjective homomorphism 
$$
\Pi:\hH^*(X) \ \arr\ \Hp *
$$
of abelian groups whose kernel is an ideal. Hence, $\Hp *$ carries a ring
structure given at the spark level by explicit formulas.
} 

\Theorem {15.2}  {\sl  Holomorphic vector bundles $E$ have   
Chern classes $d_p(E) \in \HD {2p}p$ defined as the image of the 
Cheeger-Simons classes ${\wh c}_p(E,\nabla)\in \hH^{2p}(X)$
where $\nabla$ is any canonical hermitian connection on $E$.
}
\medskip

This constitutes a form of Chern-Weil Theory for Deligne characteristic
classes.

\medskip

Proofs of these theorems appear in  [13] where Ning Hao also establishes the
following. \medskip


$\bullet$  Analytic formulas for the full product in Deligne cohomology.\medskip

$\bullet$ A \v Cech-deRham spark complex equivalent to the one given above.

\medskip

$\bullet$  A strengthened  Bott vanishing theorem 
and the construction of secondary classes in 
 
\ \ \  Deligne cohomology for holomorphic foliations.

\vskip .3in

\centerline{ APPENDIX A.\ \  HYPERCOHOMOLOGY AND CONE COMPLEXES} \medskip
In this section we present some elementary homological algebra relevant to
Deligne cohomology and the following  cone construction. 

\Def{A.1} Let $\Psi:\ci^*\to \cj^*$ be a degree-0 mapping of cochain
complexes. The associated {\bf cone complex} $C^*\equiv
\text{Cone}\,(\ci^* @>{}>>\cj^*)$ is defined by
$$
C^k \equiv \ci^{k+1}\oplus\cj^k
$$
with differential $D:C^k\to C^{k+1}$ given by
$$
D(a,b) \equiv (-da, db + \Psi(a)).
$$\medskip

Consider now a two-step complex of sheaves
$$
\ca \ @>{\psi}>>\ \cb
$$
on a manifold $X$ (with $\ca$ in degree 0), and an acyclic resolution
$$
\CD
I^*  @>{\Psi}>>  J^*\\
@AAA  @AAA  \\
\ca @>{\psi}>>  \cb
\endCD
$$

\Prop{A.2} {\sl There is a natural isomorphism
$$
H^*(X, \ca\to\cb) \ \cong\ H^*(\text{Cone}\,(\Gamma I^* @>{}>> \Gamma J^*))
$$
where $\Gamma$ denotes the sections functor.}

\Note{} This is the a relative version of the classical isomorphisms:
$$
H^*(X, \ca)\ \cong\ H^*(\Gamma I^*)
\and
H^*(X, \cb)\ \cong\ H^*(\Gamma J^*).
$$

\pf  By definition the hypercohomology $H^*(X, \ca\to\cb)$ is the total
cohomology of the double complex
$$
\CD
\vdots @. \vdots \\
@A{d}AA  @AA{d}A  \\
\Gamma I^2 @>{\Psi}>> \Gamma  J^2  \\
@A{d}AA  @AA{d}A  \\
\Gamma I^1 @>{\Psi}>> \Gamma  J^1  \\
@A{d}AA  @AA{d}A \\
\Gamma I^0 @>{\Psi}>> \Gamma  J^0  \\
\endCD
$$
whose total differential 
$$
 \Gamma I^p\oplus \Gamma J^{p-1} \ @>{D}>>\
\Gamma I^{p+1}\oplus \Gamma J^{p} 
$$
is given by $D(a,b) \equiv (-da, \Psi(a)+db)$. This is exactly the cone
complex as asserted.\qed
\medskip

The simple fact asserted in Proposition A.2 has a useful generalization.
Consider a complex of sheaves on $X$:
$$
\ca\ @>{\psi}>>\ \cb^{0}\to \cb^1\to \cb^2\to \dots\to \cb^N
$$
where $\cb^k$ has formal degree $k+1$, and suppose we are given an
acyclic resolution
$$
\CD
I^* @>{\Psi}>>  J^{**}  \\
@AAA  @AAA  \\
\ca @>{\psi}>>  \cb^*
\endCD
$$
so that $H^*(X,\ca) \cong H^*(\Gamma I^*)$ and $H^*(X,\cb^*)\cong 
H^*(\Gamma  J^{**})$. Arguing exactly as above proves the following.

\Prop{A.3} {\sl  There is a natural isomorphism
$$
H^*(X, \ca\to\cb^*) \ \cong\ 
H^*(\text{Cone}\,(\Gamma I^* @>{}>> \Gamma J^{**}))
$$
}

\vfill \eject

\vskip.3in

\centerline{ APPENDIX B.\ \  NON-SMOOTHNESS OF INTEGRAL   CURRENTS}
\medskip

Recall that a current $T$ on a manifold is called {\sl locally integral} if
both  $T$ and $dT$ are locally rectifiable. The main result of this section
is the following.

\Prop{B.1} {\sl Let $X$ be a complex manifold and for fixed integers
$k,p>0$ consider the projection
$$
\Psi:{\cd'}^k(X)\ \arr\ \Dp k 
$$
defined in \S 14.  Then for any locally rectifiable  current $T\in \In^k(X)$,
$$
\Psi(T)\ \ \text{is smooth}\qquad\Rightarrow\qquad \Psi(T)=0.
$$
}

\pf Consider an open subset $U\subset \bbr^N$.
 For any current $S\in {\cd'}^k(U)$ and any $x\in U$ we have the {\sl upper
$m$-density}:
$$
\overline{\Theta}^m(x,\|S\|)\ \equiv\ \overline{\lim_{r\to 0}}
\frac{\|S\|(B(x,r))}{\a_mr^m}
$$ 
where $\|S\|$ is the {\sl total variation measure} of $S$ (cf. [8]),
$m=N-k$ is the dimension of $S$, $B(x,r)$ is the ball of radius $r$
centered at $x$ and $\a_m$ is the volume of the unit ball in $\bbr^m$. 

We now observe that 
$$
\text{$S$ is smooth}\qquad\Rightarrow\qquad
\overline{\Theta}^m(x,\|S\|)\ =\ 0 
\tag{B.1}$$
for all $x\in U$ and all $m<N$.  To see this, suppose $S$ is smooth and fix
$\delta>0$ with $B(x,\d)\subset\subset U$. Let
$c=\sup_{|y-x|<\delta}\|S_y\|$. Then $\|S\|(B(x,r)))\leq cr^{N}$ for $r<\d$,
and so $\overline{\Theta}^m(x,\|S\|)\leq cr^{N-m}$
for $r<\d$, which proves (B.1).

Suppose now that $T$ is a locally rectifiable current of dimension $m<N$
(i.e., of degree $k=N-m>0$). Let $\ch^m$ denote Hausdorff measure in
dimension $m$.  Then $T=\eta\ch^m\bigr|_B$ for some  $(\ch^m,m)$ locally
rectifiable subset $B$ and some $\bigwedge^m\bbr^N$-valued function $\eta$
which is $L^1_{\text{loc}}$   with respect
to $\ch^m\bigl|_B$. In fact $\eta$ is a simple $m$-vector of with $|\eta|
\in \bbz^+$,  $\ch^m\bigl|_B$-a.e. The corresponding total variation
measure is just $\|T\|=|\eta|\ch^m\bigl|_B$.

Suppose now that $X$ is an open subset of $U\subset \bbc^n$ and that $T$
is a locally rectifiable current on $X$  with $\|T\|=\eta\ch^m\bigl|_B$
as above.  Then
$$
S\equiv \Psi(T) = \Psi(\eta)\ch^m\bigl|_B
$$
where $\Psi(\eta)$ denotes the pointwise projection of $\eta$ onto
$\bigoplus_{\ell<p}\bigwedge^{\ell, k-\ell}$.  

Now assume that $S$ is smooth.  Then by (B.1) we have that 
$$\overline{\Theta}^m(x,\|S\|)\
 =\ \overline{\Theta}^m(x,|\Psi(\eta)|\ch^m\bigl|_B)\ =\ 0 
$$
for all $x\in U$. However, this implies that $\Psi(\eta)=0$,
$\ch^m\bigl|_B$-a.e., and so $S=0$ as claimed.\qed

\vfill\eject
\centerline{\bf References}

\vskip .3in

\nobreak

\ref \key  1 \by R. Bott \paper On a topological obstruction to 
integrability  
\pages 127-131 \jour Proc. of Symp. in Pure Math. \vol 16 \yr 1970 
\endref

\ref\key  2 \by  J.-L. Brylinski\book Loop Spaces, Characteristic Classes
and Geometric Quantization  \publ 
 Birkhauser\publaddr Boston \yr 1993 \endref

\ref\key  3  \by J. Cheeger\paper Multiplication of 
Differential
Characters \jour Instituto Nazionale di Alta Mathematica, 
Symposia Mathematica \vol XI \yr 1973 \pages 441--445
\endref

\ref\key  4 \by J. Cheeger and J. Simons\paper
Differential Characters and Geometric Invariants
\inbook Geometry and Topology \publ Lect. Notes in Math. no. 1167,
Springer--Verlag \publaddr New York\yr 1985 \pages 50--80\endref

\ref\key  5  \by S. S.  Chern and J. Simons\paper Characteristic forms and geometric invariants \jour Ann. of Math. \vol 99 \yr 1974 \pages 48-69
\endref

\ref\key  6  \by P. Deligne\paper Th\'eorie de Hodge, II
 \jour Publ. I. H. E.S.\vol 40 \yr 1971 \pages 5-58
\endref

\ref\key  7 \by  G. de Rham \book Vari\'et\'es Diff\'erentiables, formes,
courants, formes harmoniques\publ 
 Hermann\publaddr Paris \yr 1955\endref

\ref\key  8 \by   H. Federer\book Geometric Measure 
Theory\publ 
 Springer--Verlag\publaddr New York \yr 1969\endref

\ref\key  9  \by  H. Federer and W. Fleming \paper Normal and Integral
currents \jour Annals of Math. \vol 72 \yr 1960  \pages 458-520\endref

\ref \key  10  \by H. Gillet and C. Soul\'e \paper  
Arithmetic chow groups and differential characters
 \pages 30-68 \inbook Algebraic K-theory;  Connections with Geometry and
Topology   \publ Jardine and Snaith (eds.), Kluwer Academic Publishers \yr
1989 \endref

\ref\key  11 \by   P. Griffiths and J. Harris\book Principles of
Algebraic Geometry\publ John Wiley and Sons\publaddr New York \yr
1978\endref

\ref\key  12 \by   R. Godement\book Th\'eorie des faisceaux\publ 
Hermann\publaddr Paris \yr 1964\endref

\ref\key  13 \by  N. Hao \paper 
Ph.D. Thesis\jour Stony Brook University, 2007.
\endref

\ref \key  14  \by B. Harris \paper  
Differential characters and the Abel-Jacobi map
 \pages 69-86 \inbook Algebraic K-theory;  Connections with Geometry and
Topology   \publ Jardine and Snaith (eds.), Kluwer Academic Publishers \yr
1989 \endref

\ref\key  15  \by R. Harvey and H. B. Lawson, Jr. \paper
  On boundaries of complex analytic varieties, I and II
\jour Annals of Mathematics \vol 102 and 106 \yr 1975 and 1977 \pages
223-290 and 213-238
\endref

\ref\key 16 \bysame \paper
A Theory of Characteristic Currents Associated with a Singular  Connection
\jour Ast\'erisque  \vol 213 \publ Soci\'et\'e Math. de France, Paris, 1993\endref

\ref\key 17 \bysame \paper 
Geometric residue theorems
\jour Amer. J. Math.  \vol 117 \yr 1995\pages 829-873   \endref

\ref\key 18 \bysame \paper 
Finite volume flows and Morse Theory
\jour Ann. of Math.  \vol 153 \yr 2001\pages 1-25.\ \  arXiv.math.DG/0101268   \endref

\ref\key 19 \bysame \paper 
From sparks to grundles -- differential characters
\jour Comm. in Analysis and Geometry  \vol 14 \yr 2006\pages 25-58.\ \  arXiv.math.DG/0306193    \endref

\ref\key 20 \by R. Harvey, H. B. Lawson, Jr. and J. Zweck \paper 
The de Rham-Federer theory of differential characters and character 
duality \jour Amer. J. Math. \vol 125 \yr 2003 \pages  791-847
\endref

\ref\key  21  \by R. Harvey and J. Zweck\paper Stiefel--Whitney
Currents\jour J. Geometric Analysis
\vol 8 No 5 \yr 1998 \pages  805--840
\endref

\ref\key  22  \bysame \paper Divisors and Euler sparks of atomic sections\jour Indiana Univ. Math. J.\vol 50 \yr 2001 \pages  243-298
\endref

\ref \key  23  \by J. Simons \paper Characteristic forms and
transgression: characters associated to a connection \jour Stony
Brook preprint,   \yr 1974 \endref

\ref \key  24  \by C. Voisin \paper Th\'eorie de Hodge et G\'eometrie
Alg\'ebrique Complexe \publ Soci\'et\'e   Math\'ematique de France\publaddr Paris
    2002 \endref

\end